\documentclass{amsart}
\usepackage{graphicx}
\usepackage{amsmath}
\usepackage{amssymb}
\usepackage{amsfonts}
\usepackage{amscd}
\usepackage{latexsym}
\usepackage{graphics}
\usepackage[all]{xy}

 \newtheorem{theorem}{Theorem}[section]
  \newtheorem{prop}[theorem]{Proposition}
  \newtheorem{cor}[theorem]{Corollary}
  \newtheorem{lemma}[theorem]{Lemma}

  \newtheorem{thm}[theorem]{Theorem}

\def\deg{\operatorname{deg}}
\def\det{\operatorname{det}}

\def\deg{\operatorname{deg}}

\newcommand{\cal}{\mathcal}
\renewcommand{\u}{\underline}

\newcommand{\lra}{\longrightarrow}

\newcommand{\mt}{\mapsto}
\newcommand{\h}{\simeq}

\newcommand{\mb}{\mathbb}
\newcommand{\de}{\delta}
\renewcommand{\t}{\theta}
\newcommand{\D}{\Delta}
\newcommand{\f}{\frac}
\renewcommand{\l}{\left}
\renewcommand{\r}{\right}
\renewcommand{\lg}{\langle}
\newcommand{\rg}{\rangle}
\newcommand{\be}{\begin{equation}}
\newcommand{\ee}{\end{equation}}
\newcommand{\bce}{\begin{center}}
\newcommand{\ece}{\end{center}}
\newcommand{\bq}{\begin{eqnarray}}
\newcommand{\eq}{\end{eqnarray}}
\newcommand{\n}{\nonumber}
\renewcommand{\v}{\varepsilon}
\newcommand{\ra}{\rightarrow}

\renewcommand{\a}{\alpha}
\renewcommand{\b}{\beta}
\renewcommand{\rm}{\textrm}
\newcommand{\q}{\quad}

\newcommand{\la}{\lambda}
\newcommand{\La}{\Lambda}
\newcommand{\g}{\gamma}
\newcommand{\G}{\Gamma}
\newcommand{\s}{\sum}
\renewcommand{\i}{\infty}
\renewcommand{\c}{\cdot}
\newcommand{\cs}{\cdots}
\newcommand{\ba}{\begin{array}}
\newcommand{\ea}{\end{array}}
\renewcommand{\o}{\overline}
\newcommand{\om}{\omega}
\newcommand{\Om}{\Omega}
\newcommand{\mf}{\mathfrak}

\renewcommand{\s}{\sigma}

\renewcommand{\ss}{\subset}
\newcommand{\bpm}{\begin{pmatrix}}
\newcommand{\epm}{\end{pmatrix}}
\newcommand{\bbm}{\begin{bmatrix}}
\newcommand{\ebm}{\end{bmatrix}}
\newcommand{\ld}{\ldots}
\newcommand{\tl}{\tilde}
\newcommand{\op}{\oplus}
\newcommand{\vp}{\varphi}
\newcommand{\ot}{\otimes}
\newcommand{\hr}{\hookrightarrow}

\renewcommand{\k}{\kappa}
\newcommand{\df}{\dfrac}
\newcommand{\bop}{\bigoplus}
\renewcommand{\bot}{\bigotimes}
\newcommand{\ti}{\times}
\newcommand{\wtl}{\widetilde}

\newcommand{\wh}{\widehat}

\newcommand{\vpi}{\varpi}
\newcommand{\bw}{\bigwedge}
\newcommand{\bs}{\backslash}
\newcommand{\is}{\stackrel{\sim}{\lra}}

\numberwithin{equation}{section}

\begin{document}

\title{Eisenstein Series on Loop Groups}

\author{Dongwen Liu}
\address{Department of Mathematics, University of Connecticut, Storrs, CT 06269, USA}
\email{dongwen.liu@uconn.edu}


\date{\today}


\begin{abstract}
Based on Garland's work, in this paper we construct the Eisenstein series on the adelic loop groups over a number field, induced from either a cusp form or a quasi-character which is assumed to be unramified. We compute the constant terms, prove their absolute and uniform convergence under the affine analog of Godement's criterion. For the case of quasi-characters the resulting formula is an affine Gindikin-Karpelevich formula. Then we prove the convergence of Eisenstein series themselves in certain analogs of Siegel subsets.
\end{abstract}

\maketitle

\tableofcontents

\section{Introduction}

One of the most important tools to study automorphic forms is the theory of Eisenstein series. In the fundamental work of R. Langlands \cite{Lan}, he showed how to get automorphic $L$-functions from the constant terms of Eisenstein series. This method, which was further developed by F. Shahidi and known as the Langlands-Shahidi method, has been applied to Ramanujan conjecture and Langlands functoriality \cite{kim1, KS, sha1}. On the other hand, H. Garland \cite{Gar3, Gar4, Gar5} has made very important generalizations to loop groups. He considered the Eisenstein series induced from a character, and proved the absolute convergence of constant terms first and then the Eisenstein series itself, under certain affine analog of Godement's criterion. His work lays the foundation of this field, and gives the first example of automorphic forms on infinite dimensional groups.

Based on the methods of Garland, in this paper we study the Eisenstein series defined on adelic loop groups over a number field, induced from either  a cusp form or a quasi-character which is unramified. We prove the absolute and uniform convergence of these series, and analyze their constant terms and Fourier coefficients.

Given an untwisted affine Kac-Moody Lie algebra $\wtl{\mf{g}}$ associated to
a complex simple Lie algebra $\mf{g}$, we have the affine root system $\wtl{\Phi}$ and the set of simple roots $\wtl{\D}=\{\a_0,\a_1,\ld,\a_n\}$ such that $\D=\{\a_1,\ld,\a_n\}$ is the set of simple roots of $\mf{g}.$ The affine Weyl group $\wtl{W}$ is isomorphic to the semidirect product $W\ltimes Q^\vee$ where $W$ is the Weyl group of $\mf{g}$ and $Q^\vee$ is the coroot lattice. Associated to $\wtl{\mf{g}}$ we first construct the central extension
\[
1\ra F^\times \ra \wh{G}(F((t)))\ra G(F((t)))\ra 1,
\]
then form the semi-direct product
\[
\wtl{G}(F((t)))=\wh{G}(F((t)))\rtimes \s(F^\times),
\]
where $F$ is any field, and $\s(\mf{q}),$ $\mf{q}\in F^\times$ acts on $F((t))$ as the automorphism $t\mt \mf{q}t.$ There are two methods to construct the central extension. One is via the tame symbol (\ref{tame}) on $F((t))^\times$, and the other is to use a rational representation of $G$ and the method of determinant bundles (see \cite{A}). In Theorem \ref{new} we give the explicit relation between  these two constructions. More precisely, we obtain a homomorphism between loop groups, which is identity after modulo the center, and when restricted on the center is to the power of the Dynkin index of the rational representation.

For a number field $F$ with adele ring $\mb{A}$ and idele group $\mb{I}$, we may form the adelic loop group $\wtl{G}(\mb{A}\lg t\rg)= \wh{G}(\mb{A}\lg t\rg)\rtimes\s(\mb{I})$,
where $\mb{A}\lg t\rg=\prod\limits_v'F_v((t))$ is the restricted product with respect to ${\cal O}_v((t))$ for all finite places $v.$ The ``$F$-rational points'' of the loop group is  $\wtl{G}(F\lg t\rg)$ where $F\lg t\rg=F((t))\cap\mb{A}\lg t\rg.$ We have defined the subgroups $\wh{B}_v$ and $\wh{K}_v$ of $\wh{G}(F_v((t)))$ for each place $v$, which are analogues of Borel subgroup and maximal compact subgroup respectively. More concretely, $\wh{B}_v$ is the preimage of the Borel subgroup of $G(F_v)$ under the map
\[
\wh{G}(F_v[[ t]])\lra G(F_v [[t]])\stackrel{t=0}{\lra} G(F_v).
\]
We have the Iwasawa decomposition $\wh{G}(F_v((t)))=\wh{B}_v\wh{K}_v.$ The group $\wh{G}(F_v[[t]])$ can be interpreted as the maximal parabolic subgroup of $\wh{G}(F_v((t)))$ corresponding to $\D.$ It can be shown that the central extension splits over $G(F_v[[t]]),$ i.e. we may realize $G(F_v[[t]])$ as a subgroup of $\wh{G}(F_v((t)))$ canonically. The corresponding results for the adelic groups are formulated in an obvious way. There is also the Bruhat decomposition $\wh{G}(F_v((t)))=\wh{B}_v\wtl{W}\wh{B}_v$. These results can be proved by using standard theory of Tits systems.

Fix $\mf{q}\in\mb{I}.$ If $f$ is an unramified cusp form on $G(\mb{A})$, $s\in\mb{C}$, we define a function $\wtl{f}_s$ on $\wh{G}(\mb{A}\lg t\rg)\rtimes\s(\mf{q})$ by
\[
\wtl{f}_s(g)= |c|^s f(p_0),
\]
where $g=c p \s(\mf{q}) k$ is the Iwasawa decomposition such that $c\in\mb{I},$ $p\in G(\mb{A}\lg t\rg_+)$ with $\mathbb{A}\lg t\rg_+ =\mathbb{A}\lg t\rg\cap\mathbb{A}[[t]]$, $k\in\wh{K}$, and $p_0$ is the image of $p$ under the projection $G(\mb{A}\lg t\rg_+)\stackrel{t=0}{\lra}G(\mb{A}).$ This function is well-defined and we construct the Eisenstein series defined on $\wh{G}(\mb{A}\lg t\rg)\rtimes\s(\mf{q})$ as
\[
E(s,f,g)=\sum_{\g\in \wh{G}(F\lg t\rg_+)\bs \wh{G}(F\lg t\rg)}\wtl{f}_s(\g g),
\]
where $F\lg t\rg_+=F\lg t\rg\cap F[[t]].$ The Eisenstein series is left invariant under $\wh{G}(F\lg t\rg)$ and right invariant under $\wh{K}.$ Similar construction applies for an unramified quasi-character $\chi_{\wh{T}}$ on $\wh{T}(\mb{A})/\wh{T}(F)$ where $\wh{T}$ is the maximal torus of $\wh{B}$.

The unipotent radical $\wh{U}$ of $\wh{B}$ is the subgroup corresponding to the set of all the positive roots of $\wtl{\Phi}.$ It can be proved that $\wh{U}(F)\bs \wh{U}(\mb{A})$ is compact and inherits the product measure from that of $\mb{A}/F.$ We define the constant term of $E(s,f,g)$ along $\wh{B}$ by
\[E_{\wh{B}}(s,f,g)=\int_{\wh{U}(F)\bs\wh{U}(\mb{A})}E(s,f,ug)du.
\]
The following theorem generalizes Garland's results in \cite{Gar4}.

\begin{thm}
\label{convint}
$($i$)$ Suppose that $g\in\wh{G}(\mb{A}\lg t\rg)\rtimes \s(\mf{q})$ with $\mf{q}\in\mb{I}$ and $|\mf{q}|>1,$ $s\in H=\{z\in\mb{C}|\mf{Re }z>h+h^\vee\},$ where $h$ $($resp. $h^\vee)$ is the Coxeter $($resp. dual Coxeter$)$ number. Then $E(s,f,ug)$, as a function on $\wh{U}(F)\bs\wh{U}(\mb{A})$,
converges absolutely outside a subset of measure zero and is
measurable.
\\
$($ii$)$ For any $\v, \eta>0,$ let $H_\v=\{z\in\mb{C}|\mf{Re}z>h+h^\vee+\v\},$ $\s_\eta=\{\s(\mf{q})|\mf{q}\in\mb{I}, |\mf{q}|>1+\eta\}.$ The integral defining
$E_{\wh{B}}(s,f,g)$ converges absolutely and uniformly for $s\in H_\v,$ $g\in \wh{U}(\mb{A})\Om \s_\eta\wh{K},$ where $\Om$ is a compact
subset of $T(\mb{A}).$
\\
$($iii$)$ Replace $\wtl{f}_s$ by the height function $h_s=|c|^s$, and denote the resulting series by $E(s,h,g),$ then for $a\in\wh{T}(\mb{A})$,
\be\label{zetaint} E_{\wh{B}}(s,h,a\s(\mf{q}))=\sum_{w\in
W\bs \wtl{W}}(a\s(\mf{q}))^{\tl{\rho}-w^{-1}\tl{\rho}+w^{-1}sL}c_w(s),
\ee where the summation is taken over representatives of minimal length of the cosets $W\bs\wtl{W},$ $L$ is the fundamental weight corresponding to $\a_0$, $\tl{\rho}\in\tl{\mf{h}}^\ast$ satisfies
$\lg\tilde{\rho},\a_i^\vee\rg=1,$ $i=0,1,\ld,n$ and
\be\label{cfunint} c_w(s)=\prod_{\b\in \wtl{\Phi}_+\cap
w\wtl{\Phi}_-}\f{\La_F(\lg sL-\tl{\rho},\b^\vee\rg)}{\La_F(\lg
sL-\tl{\rho},\b^\vee\rg+1)}, \ee with $\La_F$ the normalized Dedekind zeta
function.
\end{thm}

The formula (\ref{zetaint}) is an affine analogue of Gindikin-Karpelevich formula. The condition $\mf{Re}s>h+h^\vee$ is an affine analogue of Godement's criterion. Similarly if we consider the Eisenstein series $E(\chi_{\wh{T}},g)$ induced from an unramified quasi-character $\chi_{\wh{T}}$ on $\wh{T}(\mb{A})/\wh{T}(F)$, then under the condition $|\mf{q}|>1$ and $\mf{Re}(\chi_{\wh{T}}\a_i^\vee)>2$, $i=0,1,\ld,n,$ the constant term $E_{\wh{B}}(\chi_{\wh{T}}, a\s(\mf{q}))$, $a\in\wh{T}$, is given by
\be
E_{\wh{B}}(\chi_{\wh{T}},a\s(\mf{q}))=\sum_{w\in\wtl{W}}(a\s(\mf{q}))^{\wtl{\rho}+w^{-1}(\chi_{\wh{T}}-\wtl{\rho})}c_{w}(\chi_{\wh{T}})
\ee
where \be c_{w}(\chi_{\wh{T}})=\prod_{\b\in \wtl{\Phi}_+\cap
w\wtl{\Phi}_-}|\D_F|^{-\f{1}{2}}\f{L(-\lg\wtl{\rho},\b^\vee\rg, \chi_{\wh{T}} \b^\vee)}{L(1-\lg\wtl{\rho},\b^\vee\rg,\chi_{\wh{T}}\b^\vee)}.
\ee
Here $\D_F$ is the discriminant of $F$ and $L(s,\chi)$ is the Hecke $L$-function.

We have also considered the Fourier coefficients of our Eisenstein series. To obtain a general formula would be quite difficult and non-trivial. But at least for $\wtl{SL}_2$ we have computed everything explicitly, and the formulas are given in Section 4.4.

Following Garland's approach in \cite{Gar5}, we also prove some results on the absolute convergence of Eisenstein series themselves instead of the constant terms. For example we establish uniform convergence over certain analogues of Siegel sets.  The proof is technical and involves systematic use of Demazure modules together with estimations of some norms for both archimedean and non-archimedean cases. Let us only state the main results along this direction.

\begin{thm}
Fix $\mf{q}\in\mb{I},$ $|\mf{q}|>1.$ There exists a constant $c_\mf{q}>0$
depending on $\mf{q}$, such that for any $\v>0$ and compact subset
$\Om$ of $T(\mb{A})$, $E(s,f,g)$ and $E(s,h,g)$ converge absolutely
and uniformly for $s\in \{z\in\mb{C}|\mf{Re}z>\max(h+h^\vee+\v,
c_\mf{q})\}$ and $g\in \wh{U}(\mb{A})\Om \s(\mf{q}) \wh{K}.$
\end{thm}

\begin{thm}\label{convaeint}
There exist constants $c_1, c_2>0$ which depend on the number field $F$, such that for any $\v>0$ and compact subset $\Om$ of $T(\mb{A})$, $E(s,f,g)$
and $E(s,h,g)$ converge absolutely and uniformly for
$s\in\{z\in\mb{C}|\mf{Re}z>\max(h+h^\vee+\v, c_1h^\vee)\}$ and $g\in
\wh{U}(\mb{A})\Om\s_{c_2}\wh{K}$.
\end{thm}

We conjecture that Theorem \ref{convaeint} is true for $c_1=1$ (in which case the first condition reads $s\in H_\v$)
and arbitrary $c_2>0.$ In other words we conjecture that the domain of uniform convergence for the constant term $E_{\wh{B}}(s,f,g)$ in Theorem \ref{convint} also applies for $E(s,f,g)$ itself. We again interpret this as the analogue of Godement's criterion. We have proved the conjecture for $F=\mb{Q}.$ For the geometric analogue we know that the conjecture is true for $F=\mb{F}_q(T),$ the function field of $\mb{P}^1_{\mb{F}_q }.$

The theory of Eisenstein series on infinite dimensional groups is far from complete. Besides above conjecture, let us propose some other related open problems.

({\bf A}) Build the foundations of representation theory and harmonic analysis for infinite dimensional algebraic groups. Since we are dealing with groups which are not locally compact, we do not have Haar measures.  One should also concern induction from ramified representations, in contrast to what we do in this paper where we only consider induction from unramified cusp forms or quasi-characters.

({\bf B}) Generalize the theory of Eisenstein series further to all Kac-Moody groups (see \cite{Kumar} for the theory of Kac-Moody groups), and also non-split infinite dimensional groups. Compute the constant and non-constant coefficients to see if there are any new $L$-functions \cite{Gar6, sha}.

({\bf C}) Establish the Maass-Selberg relations \cite{MS1, MS2, MS3, MS4} and as applications prove the analytic continuation and functional equations for the Eisenstein series. This project, together with ({\bf B}), would be crucial for the generalization of Langlands-Shahidi method.

({\bf D}) In his thesis \cite{Pa}, M. Patnaik investigated the geometric meaning of Eisenstein series on loop groups over a function field, where he used the concept of ribbons \cite{Ka}. It would be interesting to consider this problem for the number field case.

\section{Affine Kac-Moody Lie Algebras}

In this section we review the theory of affine Kac-Moody Lie algebras. The basic references are \cite{FK, Kac, Zhu2}.

\subsection{Definition}

Let $\mf{g}$ be a complex simple finite dimensional Lie algebra. Let $(,)$ denote an invariant symmetric bilinear form on $\mf{g}$, normalized such that the square length of a long root is equal to 2. Following \cite{FK} we call it the standard bilinear form. The affine Lie algebra $\wtl{\mf{g}}$ is a complex infinite dimensional Lie algebra constructed as follows.

Let $\mb{C}[t,t^{-1}]$ be the algebra of Laurent polynomials in the indeterminate $t$ over $\mb{C}.$ For a Laurent polynomial $P=\sum c_i t^i$ the residue is defined by $\rm{Res }P=c_{-1}.$ Consider the complex infinite dimensional Lie algebra $\overline{\mf{g}}=\mb{C}[t,t^{-1}]\ot_\mb{C}\mf{g}.$ The invariant form on $\mf{g}$ can be extended naturally to a bilinear $\mb{C}[t,t^{-1}]$-valued form on $\overline{\mf{g}}$, which we again denote by $\left(,\right).$ Any derivation $D$ of $\mb{C}[t,t^{-1}]$ can be extended to a derivation of $\overline{\mf{g}}$ by $D(P\ot g)=D(P)\ot g.$

Define a $\mb{C}$-valued bilinear form $\psi$ on $\overline{\mf{g}}$ by
\[
\psi(x,y)=\rm{Res }\left( \f{dx}{dt},y\right).
\]
$\psi$ satisfies the properties

(i) $\psi(x,y)=-\psi(y,x)$ and

(ii) $\psi([x,y],z)+\psi([y,z],x)+\psi([z,x],y)=0.$

Then $\psi$ is a $2$-cocycle and we define $\hat{\mf{g}}$ to be the corresponding 1-dimensional central extension of $\overline{\mf{g}}$. The affine algebra $\wtl{\mf{g}}$ is obtained by adding to $\hat{\mf{g}}$ a derivation $d$ which acts on $\overline{\mf{g}}$ as $t\dfrac{d}{dt}$ and acts on the center as 0.

More precisely, $\wtl{\mf{g}}$ is the complex vector space
\[
\wtl{\mf{g}}=(\mb{C}[t,t^{-1}]\ot_\mb{C}\mf{g})\op\mb{C}c\op\mb{C}d
\]
with the Lie bracket:
\[
[x_1\op\a_1 c\op\b_1d,x_2\op \a_2c\op\b_2 d]=\l([x_1,x_2]+\b_1t\f{dx_2}{dt}-\b_2t\f{dx_1}{dt}\r)\op\psi(x_1,x_2)c.
\]
Here $x_i\in\overline{\mf{g}}$, $[x_1,x_2]$ is the bracket in the Lie algebra $\overline{\mf{g}}$ and $\a_i, \b_i\in\mb{C}.$

We introduce a $\mb{C}$-valued bilinear form $\left(,\right)$ on $\wtl{\mf{g}}$ by
\[
\left( x_1\op \a_1 c\op \b_1d, x_2\op \a_2c\op \b_2d\right)=\rm{Res }(t^{-1}\left( x_1, x_2\right))+\a_1\b_2+\a_2\b_1.
\]
It is easy to check that this bilinear form is symmetric, non-degenerate and invariant. Note that the restriction of the form $\left(,\right)$ to the subalgebra $\mf{g}\ss\wtl{\mf{g}}$ induces the standard bilinear form on $\mf{g}.$ Following \cite{FK} we also call the form $\left(,\right)$ on $\wtl{\mf{g}}$ the standard bilinear form.

\subsection{Root system of $\wtl{\mf{g}}$ and subalgebras in $\wtl{\mf{g}}$}

Let $\mf{h}$ denote a Cartan subalgebra of $\mf{g}.$ Let
$\mf{g}=\mf{h}\op\sum\limits_{\a\in\Phi}\mf{g}_\a$ be the root space
decomposition of $\mf{g}$ with respect to $\mf{h};$ here $\Phi\ss
\mf{h}^\ast$ is the system of roots. We fix a choice of positive
roots $\Phi_+\ss\Phi$; let $\D=\{\a_1,\ld,\a_n\}$ be the subset of
simple roots and let $\wtl{\a}$ be the highest root.

Define the following subalgebra in $\wtl{\mf{g}}$:
\[
\wtl{\mf{h}}=\mf{h}\op\mb{C}c\op\mb{C}d.
\]
This is a maximal abelian diagonalizable subalgebra in $\wtl{\mf{g}}$, and is called a Cartan subalgebra of $\wtl{\mf{g}}.$ For $\a\in\wtl{\mf{h}}^\ast$ the attached root space is
\[
\wtl{\mf{g}}_\a=\{x\in\wtl{\mf{g}}|[h,x]=\a(h)x, h\in\wtl{\mf{h}}\},
\]
and $\a$ is called a root if $\wtl{\mf{g}}_\a\neq 0.$ We extend any linear function $\la\in\mf{h}^\ast$ to a linear function on $\wtl{\mf{h}}$, which we still denote by $\la$, by setting $\la(c)=\la(d)=0.$ Let $\de\in\wtl{\mf{h}}^\ast$ be defined by $\de|_{\mf{h}+\mb{C}c}=0,$ $\de(d)=1.$ Similarly, define $L\in\wtl{\mf{h}}^\ast$ by $L|_{\mf{h}+\mb{C}d}=0,$
$L(c)=1.$

The decomposition of $\wtl{\mf{g}}$ into a sum of root spaces with respect to $\wtl{\mf{h}}$ is
\[
\wtl{\mf{g}}=\wtl{\mf{h}}\op \sum_{\a\in \Phi,
i\in\mb{Z}}(t^i\ot_\mb{C}\mf{g}_\a)\op\sum_{i\in\mb{Z}\bs\{0\}}(t^i\ot_\mb{C}\mf{h}).
\]
Therefore the root system of $\wtl{\mf{g}}$ with respect to $\wtl{\mf{h}}$ is
\[
\wtl{\Phi}=\{\a+i\de|\a\in\Phi,
i\in\mb{Z}\}\cup\{i\de|i\in\mb{Z}\bs\{0\}\}.
\]
A root $\b=\a+i\de$ with $\a\in\Phi$ is called real and a root
$\b=i\de,$ $i\in\mb{Z}\bs\{0\}$ is called imaginary. The
multiplicity $\dim\wtl{\mf{g}}_\b$ of a root $\b\in\wtl{\Phi}$ is $1$
if $\b$ is real and is $n$ otherwise.

The following properties of the standard form on $\wtl{\mf{g}}$ can be deduced from the corresponding properties of the standard form on $\mf{g}:$
\bq
&&\left(,\right)|_{\wtl{\mf{h}}}\rm{ is non-degenerate;}\n\\
&&\left(,\right)|_{\wtl{\mf{g}}_\b\op\wtl{\mf{g}}_{-\b}}\rm{ is non-degenerate;}\n\\
&&\left(\wtl{\mf{g}}_\b,\wtl{\mf{g}}_{\g}\right)=0\rm{ if }\b+\g\neq0.\n \eq
Let $\nu: \wtl{\mf{h}}\h\wtl{\mf{h}}^\ast$ be the isomorphism induced from the standard bilinear form, and we still write $(,)$ for the induced bilinear form on $\wtl{\mf{h}}^\ast$. Moreover we denote by $\lg,\rg$ the canonical pairing $\wtl{\mf{h}}^\ast\times\wtl{\mf{h}}\ra\mb{C}.$ Note that
\[
\left(\a_1+i_1\de+j_1L,\a_2+i_2\de+j_2L\right)=\left(
\a_1,\a_2\right)+i_1j_2+i_2j_1,
\]
and that a root $\a\in\wtl{\Phi}$ is real if and only
if $\left( \a,\a\right)\neq0,$ in which case $\left(\a,\a\right)>0.$ Let us write
\[
\wtl{\Phi}=\wtl{\Phi}_{re}\cup\wtl{\Phi}_{im}
\]
for the decomposition of $\wtl{\Phi}$ into real roots and imaginary roots.

Define a subsystem of positive roots $\wtl{\Phi}_+$ by
\[
\wtl{\Phi}_+=\{\a+i\de|\rm{either }i>0,\rm{ or }i=0, \a\in\Phi_+\}.
\]
Then $\wtl{\Phi}=\wtl{\Phi}_+\cup(-\wtl{\Phi}_+)$, and the
corresponding system of simple roots is
\[
\wtl{\D}=\{\a_0=\de-\wtl{\a}, \a_1,\ld,\a_n\}.
\]
The following subalgebras of $\wtl{\mf{g}}$ are the analogues of the maximal nilpotent and the Borel subalgebras of $\mf{g}$:
\[
\wtl{\mf{n}}_+=\bop_{\b\in\wtl{\Phi}_+}\wtl{\mf{g}}_\b,\q
\wtl{\mf{n}}_-=\bop_{\b\in\wtl{\Phi}_+}\wtl{\mf{g}}_{-\b}, \q
\wtl{\mf{b}}=\wtl{\mf{h}}\op\wtl{\mf{n}}_+.
\]

\subsection{Affine Weyl group of $\wtl{\mf{g}}$}

For a real root $\b\in\wtl{\Phi}_{re}$, let $\b^\vee\in\wtl{\mf{h}}$ be
the coroot.

\begin{lemma}\label{coroot}
Let $\a^\vee\in\mf{h}$ be the coroot of the root $\a\in\Phi,$ then
the coroot of $\b=i\de+\a$ is $\b^\vee=\df{i}{2}\left(
\a^\vee,\a^\vee\right) c+\a^\vee.$ In particular, the coroot of $\a_0$
is $c-\wtl{\a}^\vee.$
\end{lemma}
\textit{Proof.} Let $x_\a, x_{-\a}$ and $\a^\vee$ be a standard basis
of $\mf{g}_\a+\mf{g}_{-\a}+\mb{C}\a^\vee\h \mf{sl}_2,$ then \bq
\b^\vee &=&[t^i\ot x_\a,t^{-i}\ot x_{-\a}]\n\\
&=&[x_\a,x_{-\a}]+i\left( x_\a,x_{-\a}\right) c\n\\
&=&\a^\vee+\f{i}{2}\left( [\a^\vee,x_\a],x_{-\a}\right) c\n\\
&=&\a^\vee+\f{i}{2}\left( \a^\vee, [x_\a, x_{-\a}]\right) c\n\\
&=&\a^\vee +\f{i}{2}\left(\a^\vee,\a^\vee\right) c.\n
\eq
One can check easily that
\[
[\b^\vee, t^{i}\ot x_\a]=2t^{i}\ot x_\a,\q [\b^\vee, t^{-i}\ot x_{-\a}]=-2t^{-i}\ot x_{-\a}.
\]
\hfill$\Box$

Let $\rho\in\mf{h}^\ast$ be the half sum of all the positive roots in $\Phi,$ then $\lg\rho,\a_i^\vee\rg=1$ for $i=1,\ld,n.$ Let $\wtl{\rho}=\rho+(1+\lg\rho,\wtl{\a}^\vee\rg)L\in \wtl{\mf{h}}^\ast$, then $\lg\wtl{\rho},\a_i^\vee\rg=1$ for $i=0,1,\ld,n.$ The number $h^\vee:=1+\lg \rho, \wtl{\a}^\vee\rg$ is called the dual Coxeter number of the root system $\Phi.$ Therefore $\wtl{\rho}=\rho+h^\vee L.$

For a real root $\b\in\wtl{\Phi}_{re}$, let $r_\b$ be the reflection
whose action on $\wtl{\mf{h}}^\ast$ is given by
\[
r_\b(\la)=\la-\lg\la,\b^\vee\rg\b,\q \la\in \wtl{\mf{h}}^\ast,
\]
and whose action on $\wtl{\mf{h}}$ is given by
\[
r_\b(h)=h-\lg \b,h\rg\b^\vee,\q h\in \wtl{\mf{h}}.
\]
The two actions are dual to each other:
\[
\lg r_\b(\la),r_\b(h)\rg=\lg \la,h\rg,\q h\in\wtl{\mf{h}}, \la\in \wtl{\mf{h}}^\ast.
\]
The group $\wtl{W}\ss GL(\wtl{\mf{h}}^\ast)$ generated by $r_\b$'s
over all real roots $\b\in \wtl{\Phi}_{re}$ is called the affine Weyl
group of $\wtl{\mf{g}}.$ The form $\left(,\right)|_{\wtl{\mf{h}}^\ast}$ is
$\wtl{W}$-invariant. Note that any real root is a
$\wtl{W}$-conjugate of a simple root and the line $\mb{C}\de$ is the
fixed point set for $\wtl{W}.$ Write $r_i$ instead of $r_{\a_i}$,
$i=0,1,\ld,n.$ Then the group $\wtl{W}$ is generated by $r_i$'s.

Let $W$ be the Weyl group of $\mf{g}$, which can be identified with the subgroup of $\wtl{W}$ generated by the reflections $r_1,\ld,r_n.$

Let $Q$ be the root lattice of $\mf{g},$ i.e., the $\mb{Z}$-lattice
generated by $\D,$ and let $Q^\vee$ denote the coroot lattice, i.e.,
the lattice generated by $\a_i^\vee, i=1,2,\ld,n.$ It is known that
$\a^\vee\in Q^\vee$ for all $\a\in\Phi.$

\begin{thm}
The affine Weyl group $\wtl{W}$ is isomorphic to $W\ltimes Q^\vee.$
Write $T_h$ for $h\in Q^\vee$ as an element in $W\ltimes Q^\vee.$
Then the isomorphism is given by, for $\a\in\Phi,$
\[
r_\a\mt r_\a,\q r_{i\de-\a}r_\a\mt T_{i\a^\vee}.
\]
\end{thm}
See \cite{FK, Kac, Kumar}. The following two lemmas give the explicit action of $\wtl{W}$ on $\wtl{\mf{h}}$ and $\wtl{\mf{h}}^\ast.$

\begin{lemma}\label{weyl1}
The $\wtl{W}$-action on $\wtl{\mf{h}}$ fixes $c$. The action is given
by the formula: for $\a\in\Phi,$ $\g\in Q^\vee,$ $h\in\mf{h}$, $i\in\mb{Z}$, \bq
r_\a(h+id)&=&r_\a(h)+id,\n\\
T_{\g}(h+id)&=&h+id+i\g-\l(\left( h,\g\right)+\f{i}{2}\left(\g,\g\right)\r)c,\n
\eq
or equivalently,
\[
T_{\g}(x)=x+\left( x,c\right)\g-\l(\left( h,\g\right)+\f{1}{2}\left( \g,\g\right)\left( x,c\right)\r)c,\q \forall x\in\wtl{\mf{h}}.
\]
\end{lemma}
\textit{Proof.} Since $\left( c,\b\right)=0$ for every real root $\b,$ $c$ is
fixed by $\wtl{W}.$ We prove the second formula for $\g=\a^\vee$
with $\a\in\Phi.$ The general case can be reduced to this one. \bq
T_{\a^\vee}(h+kd)&=&r_{\de-\a^\vee}r_\a(h+id)\n\\
&=&r_{\de-\a}(r_\a(h)+id)\n\\
&=&r_\a(h)+id-\left( r_\a(h)+id, \de-\a\right)\l(\f{(\a^\vee,\a^\vee)}{2}c-\a^\vee\r)\n\\
&=&r_\a(h)+id+(i+\left( h,\a\right))\l(\a^\vee-\f{(\a^\vee,\a^\vee)}{2}c\r)\n\\
&=&h+id+i\a^\vee-\l(\left( h,\a^\vee\right)+\f{i}{2}\left(\a^\vee,\a^\vee\right)\r)c.\n
\eq
\hfill$\Box$

\begin{lemma} \label{weyl2}
The $\wtl{W}$-action on $\wtl{\mf{h}}^\ast$ fixes $\de$. The action
is given by the formula: for $\a\in\Phi,$ $\g\in Q^\vee,$
$\la\in\mf{h}^\ast$, $i\in\mb{Z}$,\bq
r_\a(\la+iL)&=&r_\a(\la)+iL,\n\\
T_\g(\la+iL)&=&\la+iL+i\nu(\g)-\l(\lg\la,\g\rg+\f{i}{2}\left( \g,\g\right)\r)\de,
\n
\eq
or equivalently,
\[
T_\g(x)=x+\left( x, \de\right)\nu(\g)-\l(\lg x,\g\rg+\f{1}{2}\left( \g,\g\right)\left( x,\de\right)\r)\de,\q \forall x\in \wtl{\mf{h}}^\ast.
\]
In particular,
\[
T_\g(\b)=\b-\lg\b,\g\rg\de,\q \b\in\wtl{\Phi}.
\]
\end{lemma}

The proof of this lemma is similar to that of Lemma \ref{weyl1}, and
can be reduced to the case $\g=\a^\vee$ with $\a\in\Phi.$

\section{Constructions of Loop Groups}

We shall construct the loop groups associated with complex
simple Lie algebras, and obtain central extensions of loop groups by using Tame
symbols. Then we discuss the
highest weight representations of loop groups. Another construction of loop groups starting from a linear algebraic
group and a rational representation of this group will be given, and we will see the relationships between these two constructions.
We also construct adelic loop groups and review some
fundamental results of H. Garland \cite{Gar1} on arithmetic
quotients.

\subsection{First construction of loop groups}

We first recall the definition of Chevalley groups. The main references are \cite{Ma, S}.
Let $\mf{g}$ be a complex simple Lie algebra, and we use the same
notations in section 2. Fix a Chevalley basis of $\mf{g}.$ The universal Chevalley
group associated to $\mf{g}$ is a simply connected affine group scheme $G$ over $\mb{Z}$, and for any field $F$ the
 $F$-rational points $G(F)$ of $G$ is generated by the elements $x_\a(u),$
 $\a\in\Phi,$ $u\in F$ subject to relations (\ref{R1})-(\ref{R3}) if rank $\mf{g}\geq2,$ or relations (\ref{R1}), (\ref{R3})
 and (\ref{R2'}) if $\mf{g}=\mf{sl}_2.$

 For $\a\in\Phi,$ $u, v\in F,$
 \be \label{R1}
x_\a(u)x_\a(v)=x_\a(u+v).
 \ee
 For $\a,$ $\b\in\Phi$, $\a\neq-\b,$ $u,$ $v\in F$,
 \be \label{R2}
 x_\a(u)x_\b(v)x_\a(u)^{-1}x_\b(v)^{-1}=\prod_{i,j\in\mb{Z}^+,i\a+j\b\in\Phi}x_{i\a+j\b}(c_{ij}^{\a\b}u^iv^j),
 \ee
where the order of the right-hand-side is given by some fixed order, and the coefficients $c_{ij}^{\a\b}$ are integers
which depend on this order and the Chevalley basis of $\mf{g}$, but not on the field $F$ or on $u,$ $v.$ For $\a\in\Phi$, $u\in F^\times$ we set
\[
w_\a(u)=x_\a(u)x_{-\a}(-u^{-1})x_\a(u),\q
h_\a(u)=w_\a(u)w_\a(1)^{-1}.
\]
Then for $u,$ $v\in F^\times,$
\be\label{R3} h_{\a}(u)h_{\a}(v)=h_{\a}(uv).\ee
If $\mf{g}=\mf{sl}_2,$ there are only two roots $\pm\a$, and the relation (\ref{R2}) above is
replaced by \be\label{R2'} w_\a(u)x_\a(v)w_\a(-u)=x_{-\a}(-u^2v),\q u\in F^\times, v\in F.\ee

The universal Steinberg group $G'(F)$ is generated by
$\wtl{x}_\a(u),$ $\a\in \Phi$, $u\in F$ subject to relations (\ref{R1}) and (\ref{R2}) if rank $\mf{g}\geq 2$, or
(\ref{R1}) and (\ref{R2'}) if $\mf{g}=\mf{sl}_2.$ Here for $\a\in\Phi$, $u\in F^\times$ we define
\[
\wtl{w}_\a(u)=\wtl{x}_\a(u)\wtl{x}_{-\a}(-u^{-1})\wtl{x}_\a(u),\q
\wtl{h}_\a(u)=\wtl{w}_\a(u)\wtl{w}_\a(1)^{-1}.
\]
Let $\pi: G'(F)\ra G(F)$ be the homomorphism defined by
$\pi(\wtl{x}_\a(u))=x_\a(u)$ for all $\a\in\Phi,$ $u\in F.$
Steinberg (\cite[p.78]{S} Thoeorem 10) proved that, if $|F|>4$, and $|F|\neq 9$ when $\mf{g}=\mf{sl}_2,$ then
$(\pi,G')$ is a universal central extension of $G$. Recall from \cite[p.74]{S} that a central extension
$(\pi, E)$ of a group $G$ is universal if for any central extension $(\pi',E')$ of $G$ there exists a
unique homomorphism $\vp: E\ra E'$ such that $\pi'\vp=\pi,$ i.e. the following diagram is
commutative:
\begin{displaymath}
\xymatrix{ E \ar[rr]^{\vp} \ar[dr]_\pi & & E' \ar[dl]^{\pi'}\\
& G & }
\end{displaymath}

Let $C=\rm{Ker }\pi$. Matsumoto \cite{Ma} and Moore \cite{Mo}  (cf. \cite[pp.86-87]{S} Theorem 12) proved that if $|F|>4$,
then $C$ is isomorphic to the abstract group generated by the symbols $c(u,v)$ ($u, v\in F^\times$) subject to the relations
\bq
 \label{a} && c(u,v)c(uv,w)=c(u,vw)c(v,w),\q c(1,u)=c(u,1)=1,\\
 \label{b} && c(u,v)c(u,-v^{-1})=c(u,-1),\\
 \label{c} && c(u,v)=c(v^{-1},u),\\
 \label{d} && c(u,v)=c(u,-uv),\\
 \label{e} && c(u,v)=c(u,(1-u)v)
\eq
and in the case $\Phi$ is not of type $C_n$ ($n\geq1$) the additional relation
\be \label{f}
c\rm{ is bimultiplicative.}\q\q\q\q\q\q\q\q\q\q\q\q\q
\ee
In this case relations (\ref{a})-(\ref{e}) may be replaced by (\ref{f}) and
\bq
 \label{c'} && c\rm{ is skew,}\\
 \label{d'} && c(u,-u)=1.\\
 \label{e'} && c(u,1-u)=1.
\eq
The isomorphism is given by $c(u,v)\mt \wtl{h}_\a(u)\wtl{h}_\a(v)\wtl{h}_\a(uv)^{-1},$
where $\a$ is a fixed long root. For a field $F$ and an abelian group $A$, a map $c:
F^\times\times F^\times\ra A$ is called a Steinberg symbol on $F^\times\times F^\times$ with values in $A$ if it
satisfies the relations (\ref{a})-(\ref{e}), and is said to be bilinear if it also satisfies
(\ref{f}).

In the Steinberg group $G'$ let $c_\a(u,v)=\wtl{h}_\a(u)\wtl{h}_\a(v)\wtl{h}_\a(uv)^{-1}$, $\a\in\Phi$,
$u, v\in F^\times.$

\begin{lemma}\label{symb} $($\cite[p.26]{Ma}~Lemma~$5.4)$\\
$(a)$ $c_\a(u,v)=c_{-\a}(v,u)^{-1}$, $\forall\a\in\Phi,$\\
$(b)$ If there exists $w\in W$ such that $\b=w\a,$ then $c_\b$ equals $c_\a$ or $c_{-\a},$\\
$(c)$ For $\a, \b\in\Phi$,
\[
\wtl{h}_\a(u)\wtl{h}_\b(v)\wtl{h}_\a(u)^{-1}\wtl{h}_\b(v)^{-1}=c_\a(u,v^{\lg\a,\b^\vee\rg})=c_\b(v,u^{\lg\b,\a^\vee\rg})^{-1},
\]
$(d)$ The Steinberg symbol $c_\a$ is blinear except for the case $G$ is symplectic and $\a$ is a long root.
\end{lemma}

Suppose that $c: F^\times\times F^\times\ra A$ is a Steinberg symbol. By \cite[p.30]{Ma} Th\'{e}or\`{e}m 5.10, there exists a central extension
of $G(F)$ by $A$ such that $c_\a=c$ for any long root $\a$, if either $c$
is bilinear or $G$ is symplectic. In fact, the symbol gives a homomorphism of abelian groups $\phi: C\ra A$. We may assume that $\phi$ is surjective. Then from the universal central
 extension $1\ra C\ra G'(F)\ra G(F)\ra 1$ we obtain
  \[
  1\lra \f{C}{\rm{Ker}\phi}\lra \f{G'(F)}{\rm{Ker}\phi}\lra G(F)\lra 1.
  \]
  Then $G'(F)/\rm{Ker}\phi$ is the required central extension.
 If $c$ is bilinear, then $c_\a=c^{\f{2}{\left(\a,\a\right)}}$, $\forall\a\in\Phi.$ This
can be proved by using Lemma \ref{symb} (c) and checking the Dynkin diagrams. Recall that the square length of a long root equals 2.

The following lemma follows from \cite{S2}, see \cite[pp.23-24]{Ma} Lemme 5.1 and Lemme 5.2.

\begin{lemma}\label{ma}
In a central extension of $G$ by a Steinberg symbol we have the following relations for $\a,\b\in\Phi:$\\
$(a)$ $\wtl{w}_\a(u)\wtl{x}_\b(v)\wtl{w}_\a(u)^{-1}=\wtl{x}_{r_\a\b}(\eta_{\a,\b}u^{-\lg\b,\a^\vee\rg}v)$, where $\eta_{\a,\b}$ are integers equal to $\pm1$ given by \cite{Ma} Lemme 5.1 $(c)$,\\
$(b)$ $\wtl{w}_\a(u)\wtl{h}_\b(v)\wtl{w}_\a(u)^{-1}=\wtl{h}_{r_\a \b}(\eta_{\a,\b}u^{-\lg\b,\a^\vee\rg}v)\wtl{h}_{r_\a\b}(\eta_{\a,\b}u^{-\lg\b,\a^\vee\rg})^{-1}$, i.e.
\[
\wtl{w}_\a(u)\wtl{h}_\b(v)\wtl{w}_\a(u)^{-1}=c_{r_\a\b}(v,\eta_{\a,\b}u^{-\lg\b,\a^\vee\rg})^{-1}\wtl{h}_{r_a\b}(v),
\]
$(c)$ $\wtl{h}_\a(u)\wtl{x}_\b(v)\wtl{h}_\a(u)^{-1}=\wtl{x}_\b(u^{\lg\b,\a^\vee\rg}v),$ \\
$(d)$ $\wtl{h}_\a(u)\wtl{w}_\b(v)\wtl{h}_\a(u)^{-1}=\wtl{w}_\b(u^{\lg\b,\a^\vee\rg}v),$\\
$(e)$ $\wtl{w}_\a(1)\wtl{h}_\b(u)\wtl{w}_\a(1)^{-1}=\wtl{h}_\b(u)\wtl{h}_\a(u^{-\lg\a,\b^\vee\rg}),$\\
$(f)$ $\wtl{w}_\a(u)=\wtl{w}_{-\a}(-u^{-1})$, $\wtl{h}_\a(u)=\wtl{h}_{-\a}(u)^{-1}$, $\wtl{w}_\a(1)\wtl{h}_\a(u)\wtl{w}_\a(1)^{-1}=\wtl{h}_\a(u^{-1}),$\\
$(g)$ $\wtl{w}_\a(1)^{-1}\wtl{x}_\a(u)\wtl{w}_\a(1)=\wtl{x}_{-\a}(-u)=\wtl{x}_\a(-u^{-1})\wtl{w}_\a(u^{-1})\wtl{x}_\a(-u^{-1}),$ $u\neq 0.$
\end{lemma}

The tame symbol defined for the field of Laurent power series
$F((t))$ is the map $(,)_{\rm{tame}}: F((t))^\times\times
F((t))^\times\ra F^\times$ given by
\be\label{tame}
(x,y)_{tame}=(-1)^{v(x)v(y)}\l.\f{x^{v(y)}}{y^{v(x)}}\r|_{t=0},
\ee
where $v$ is the valuation on $F((t))$ normalized such that
$v(t^i)=i.$ Note that tame symbol is trivial on $F^\times\times
F^\times.$

Since the tame symbol is a
bilinear Steinberg symbol, we obtain a central extension of
$G(F((t)))$ by $F^\times,$ assocaited to the inverse of the tame symbol. Let us denote this
central extension by $\wh{G}(F((t))).$ It is generated by
$\wtl{x}_\a(u)$ with $\a\in\Phi,$ $u\in F((t))$ and $F^\times,$
subject to relations (\ref{R1}), (\ref{R2}) and (\ref{RC}) below if rank $\mf{g}\geq 2,$
or relations (\ref{R1}), (\ref{R2'}) and (\ref{RC}) if $\mf{g}=\mf{sl}_2.$  By previous remarks, for each $\a\in\Phi,$
\be\label{RC}
\wtl{h}_\a(x)\wtl{h}_\a(y)\wtl{h}_\a(xy)^{-1}=(x,y)_{tame}^{-\f{2}{\left(\a,\a\right)}},\q x,y\in F((t))^\times.
\ee

Then we have the following exact sequence for
$G(F((t)))$:
\[
1\lra F^\times\lra \wh{G}(F((t)))\stackrel{\pi}{\lra} G(F((t)))\lra 1,
\]
where $\pi$ is given by $\wtl{x}_\a(u)\mt x_\a(u).$

For a real root $\b=\a+i\de\in\wtl{\Phi}_{re}$, and $u\in F,$ $v\in
F^\times,$ we define for  $G(F((t)))$,
\be\label{aff}
\l\{\ba{l}x_\b(u)=x_\a(ut^i),\\
  w_\b(v)=x_\b(v)x_{-\b}(-v^{-1})x_\b(v)=w_\a(vt^i),\\
 h_\b(v)=w_\b(v)w_\b(1)^{-1}=h_\a(v). \ea\r.\ee
 For $\wh{G}(F((t)))$ we can
define the elements by the same formula with $x, w, h$ replaced by
$\wtl{x}, \wtl{w}, \wtl{h}$. From the definition we have \bq
\label{2h}\wtl{h}_\b(u)&=&\wtl{w}_\a(ut^i)\wtl{w}_\a(t^i)^{-1}=\wtl{w}_\a(ut^i)\wtl{w}_\a(1)^{-1}(\wtl{w}_\a(t^i)\wtl{w}_\a(1)^{-1})^{-1}\\
&=&\wtl{h}_\a(ut^i)\wtl{h}_\a(t^i)^{-1}= (u,t^i)_{tame}^{\f{2}{\left(\a,\a\right)}}\wtl{h}_\a(u)\n\\
&=& u^{\f{2i}{\left(\a,\a\right)}}\wtl{h}_\a(u)= u^{\f{2i}{\left(\b,\b\right)}}\wtl{h}_\a(u). \n \eq

It is clear that $\{x_\b(u)|u\in F\}$ (resp. $\{\wtl{x}_\b(u)|u\in
F\}$) forms a subgroup isomorphic to the additive group
$\mb{G}_a(F)$. We call it the root subgroup associated to $\b,$ and
denote by $U_\b$ (resp. $\wtl{U}_\b$).

For a positive imaginary root $\b=i\de\in\wtl{\Phi}_{im+},$ $i\in\mb{N},$ we define the root subgroup $U_\b$ (resp. $\wtl{U}_\b$) as follows, which is isomorphic to $\mb{G}_a^n$. The map
\[
\exp: tF[[t]]\lra 1+tF[[t]]
\]
is a bijection with inverse map $\log.$ The root subgroup $U_\b$ is given by
\be\label{imroot}
U_\b=\{h_{\a_1}(\exp(u_1t^i))\cs h_{\a_n}(\exp(u_n t^i))| u_1,\ld, u_n\in F\}.
\ee
We define $\wtl{U}_\b$ similarly, with $h$ replaced by $\wtl{h}.$

\begin{lemma}\label{sl2}
For each real root $\b=\a+i\de\in\wtl{\Phi}_{re},$ there is a unique
group homomorphism $\vp_\b: SL_2(F)\ra G(F((t)))$ $(\rm{resp. }
\wh{G}(F((t))))$ such that
\[
\bpm 1&u\\0&1\epm\mt x_\b(u)~(\rm{resp. } \wtl{x}_\b(u)),\q\bpm
1&0\\u&1\epm\mt x_{-\b}(u)~(\rm{resp. }\wtl{x}_{-\b}(u)).
\]
\end{lemma}
\textit{Proof.} For $G(F((t)))$ we have
\[
w_\b(u)=x_\a(ut^i)x_{-\a}(-u^{-1}t^{-i})x_{-\a}(ut^i)=w_\a(ut^i),
\]
\[
h_\b(u)=h_\a(ut^i)h_\a(t^i)^{-1}=h_\a(u),
\]
and therefore
\[
h_\b(u)h_\b(v)h_\b(uv)^{-1}=h_\a(u)h_\a(v)h_\a(uv)^{-1}=1.
\]
For $\wh{G}(F((t)))$, we only need to verify the last equation above
since others are similar. However from (\ref{2h}) it follows that
\[
\wtl{h}_\b(u)\wtl{h}_\b(v)\wtl{h}_\b(uv)^{-1}=\wtl{h}_\a(u)\wtl{h}_\a(v)\wtl{h}_\a(uv)^{-1}=1.
\]
This verifies that $x_\b(u)$ (resp. $\wtl{x}_\b(u))$ and $x_{-\b}(u)$
(resp. $\wtl{x}_{-\b}(u)$) satisfy the relations of $SL_2(F).$ \hfill
$\Box$

\

Apply (\ref{2h}) and use properties of the tame symbol, we can translate Lemma \ref{symb} and Lemma \ref{ma} into the data of affine root system $\wtl{\Phi}.$ Assume that $\a=\a_0+i\de,$ $\b=\b_0+j\de\in\wtl{\Phi}_{re},$ where $\a_0$, $\b_0\in\Phi,$ $i, j\in\mb{Z}.$ One should not confuse $\a_0$ with the simple root $\a_0=\de-\wtl{\a}.$ Let $\eta_{\a,\b}=\eta_{\a_0,\b_0}.$

\begin{cor}\label{r}
We have the following relations in $\wh{G}(F((t)))$ for $\a, \b\in\wtl{\Phi}_{re}:$ \\
$(a)$ $\wtl{x}_\a(u)\wtl{x}_\b(v)\wtl{x}_\a(u)^{-1}\wtl{x}_\b(v)^{-1}=\prod\limits_{m,n\in\mb{Z}^+, m\a+n\b\in\wtl{\Phi}_{re}}\wtl{x}_{m\a+n\b}(c^{\a_0\b_0}_{mn}u^m v^n)$, $u, v\in F$,\\
$(b)$ $\wtl{h}_\a(u)\wtl{h}_\a(v)=\wtl{h}_\a(uv)$, $u, v\in F^\times$,\\
$(c)$ $\wtl{h}_\a(u)\wtl{h}_\b(v)=\wtl{h}_\b(v)\wtl{h}_\a(u)$, $u, v\in F^\times$,\\
$(d)$ $\wtl{w}_\a(u)\wtl{x}_\b(v)\wtl{w}_\a(u)^{-1}=\wtl{x}_{r_\a\b}(\eta_{\a,\b}u^{-\lg\b,\a^\vee\rg}v)$, $u\in F^\times$, $v\in F$,\\
$(e)$ $\wtl{w}_\a(u)\wtl{h}_\b(v)\wtl{w}_\a(u)^{-1}=\wtl{h}_{r_\a \b}(v)$, $u, v\in F^\times$,\\
$(f)$ $\wtl{h}_\a(u)\wtl{x}_\b(v)\wtl{h}_\a(u)^{-1}=\wtl{x}_\b(u^{\lg\b,\a^\vee\rg}v),$ $u\in F^\times$, $v\in F$, \\
$(g)$ $\wtl{h}_\a(u)\wtl{w}_\b(v)\wtl{h}_\a(u)^{-1}=\wtl{w}_\b(u^{\lg\b,\a^\vee\rg}v),$ $u, v\in F^\times$,\\
$(h)$ $\wtl{w}_\a(1)\wtl{h}_\b(u)\wtl{w}_\a(1)^{-1}=\wtl{h}_\b(u)\wtl{h}_\a(u^{-\lg\a,\b^\vee\rg}),$ $u\in F^\times$,\\
$(i)$ $\wtl{w}_\a(u)=\wtl{w}_{-\a}(-u^{-1})$, $\wtl{h}_\a(u)^{-1}=\wtl{h}_{\a}(u^{-1})$, $\wtl{w}_\a(1)\wtl{h}_\a(u)\wtl{w}_\a(1)^{-1}=\wtl{h}_\a(u^{-1}),$ $u\in F^\times$,\\
$(j)$ $\wtl{w}_\a(1)^{-1}\wtl{x}_\a(u)\wtl{w}_\a(1)=\wtl{x}_{-\a}(-u)=\wtl{x}_\a(-u^{-1})\wtl{w}_\a(u^{-1})\wtl{x}_\a(-u^{-1}),$ $u\in F^\times.$
\end{cor}
\textit{Proof.} ($a$) By (\ref{R2}),
\bq
&&\wtl{x}_\a(u)\wtl{x}_\b(v)\wtl{x}_\a(u)^{-1}\wtl{x}_\b(v)^{-1}\n\\
&=& \wtl{x}_{\a_0}(ut^i)\wtl{x}_{\b_0}(vt^j)\wtl{x}_{\a_0}(ut^i)^{-1}\wtl{x}_{\b_0}(vt^j)^{-1}\n\\
&=& \prod_{m,n\in\mb{Z}^+, m\a_0+n\b_0\in\Phi}\wtl{x}_{m\a_0+n\b_0}(c^{\a_0\b_0}_{mn}(ut^i)^m (vt^j)^n)\n\\&=& \prod_{m,n\in\mb{Z}^+, m\a+n\b\in\wtl{\Phi}_{re}}\wtl{x}_{m\a+n\b}(c^{\a_0\b_0}_{mn}u^m v^n).\n
\eq
($b$): Since the tame symbol is trivial on $F^\times\times F^\times$,
\bq \wtl{h}_\a(u)\wtl{h}_\a(v)&=& u^{\f{2i}{\left(\a,\a\right)}}v^{\f{2i}{\left(\a,\a\right)}}\wtl{h}_{\a_0}(u)\wtl{h}_{\a_0}(v)\n\\
&=&(uv)^{\f{2i}{\left(\a,\a\right)}}\wtl{h}_{\a_0}(uv)=\wtl{h}_\a(uv).\n \eq
($c$): By Lemma \ref{symb} ($c$), \bq
\wtl{h}_\a(u)\wtl{h}_\b(v)&=& u^{\f{2i}{\left(\a,\a\right)}}v^{\f{2j}{\left(\b,\b\right)}}\wtl{h}_{\a_0}(u)\wtl{h}_{\b_0}(v)\n\\
&=&u^{\f{2i}{\left(\a,\a\right)}}v^{\f{2j}{\left(\b,\b\right)}}\wtl{h}_{\b_0}(v)\wtl{h}_{\a_0}(u)\n\\
&=&\wtl{h}_\b(v)\wtl{h}_\a(u).\n \eq
($d$): By Lemma \ref{ma} ($a$), \bq
\wtl{w}_\a(u)\wtl{x}_\b(v)\wtl{w}_\a(u)^{-1}&=&\wtl{w}_{\a_0}(ut^i)\wtl{x}_{\b_0}(vt^j)\wtl{w}_{\a_0}(ut^i)^{-1}\n\\
&=&\wtl{x}_{r_{\a_0}\b_0}((ut^i)^{-\lg\b_0,\a_0^\vee\rg}vt^j)\n\\
&=& \wtl{x}_{r_\a\b}(u^{-\lg\b,\a^\vee\rg}v),\n
\eq
where the last equality follows from the formulas
\[
\lg\b,\a^\vee\rg=\lg\b_0,\a_0^\vee\rg,\q r_\a\b=r_{\a_0}\b_0+(j-\lg\b_0,\a_0^\vee\rg i)\de.
\]
($e$): By Lemma \ref{ma} ($b$), \bq
&& \wtl{w}_\a(u)\wtl{h}_\b(v)\wtl{w}_\a(u)^{-1}\n\\
&=& v^{\f{2j}{\left(\b,\b\right)}}\wtl{w}_{\a_0}(ut^i)h_{\b_0}(v)\wtl{w}_{\a_0}(ut^i)^{-1}\n\\
&=& v^{\f{2j}{\left(\b,\b\right)}}\l(v,\eta_{\a,\b}(ut^i)^{-\lg\b,\a^\vee\rg}\r)_{tame}^{\f{2}{\left(\b,\b\right)}}\wtl{h}_{r_{\a_0}\b_0}(v)\n\\
&=& v^{\f{2(j-\lg\b,\a^\vee\rg i)}{\left(\b,\b\right)}}\wtl{h}_{r_{\a_0}\b_0}(v)\n\\
&=& \wtl{h}_{r_\a\b}(v).\n
\eq
($f$) and ($g$) are easy consequences of (\ref{aff}) and (\ref{2h}). By ($d$) the left-hand-side of ($h$) equals $\wtl{h}_{r_\a\b}(u).$ By
Lemma \ref{ma} ($e$), the right-hand-side of ($h$) equals
\bq
&& u^{\f{2j}{\left(\b,\b\right)}-\f{2i\lg\a,\b^\vee\rg}{\left(\a,\a\right)}}\wtl{h}_{\b_0}(u)\wtl{h}_{\a_0}(u^{-\lg\a,\b^\vee\rg})\n\\
&=&u^{\f{2(j-\lg\b,\a^\vee\rg i)}{\left(\b,\b\right)}}\wtl{w}_{\a_0}(1)\wtl{h}_{\b_0}(u)\wtl{w}_{\a_0}(1)^{-1}\n\\
&=&u^{\f{2(j-\lg\b,\a^\vee\rg i)}{\left(\b,\b\right)}}\wtl{h}_{r_{\a_0}\b_0}(u)\n\\
&=&\wtl{h}_{r_\a\b}(u).\n
\eq
This proves ($h$). ($i$) follows from ($a$) and ($g$). ($j$) follows from ($d$) and the facts $\eta_{\a,\a}=\eta_{\a,-\a}=1$,
see \cite[p.24]{Ma} Lemme 5.1 (c).
\hfill$\Box$

\

For each $\a\in\Phi_+,$ the subgroup of $G(F)$ generated by $x_\a(u)$ and
$x_{-\a}(u)$ ($u\in F$) is isomorphic to $SL_2(F)$ by the map
\be\label{cpt} x_\a(u)\mt\bpm 1&u\\0&1\epm,\q x_{-\a}(u)\mt\bpm
1&0\\u&1\epm.\ee The subgroup $B(F)$ generated by $x_\a(u)$
($\a\in\Phi_+$) and $h_\a(u)$ is a Borel subgroup, and the subgroup
generated by $x_\a(u)$ ($\a\in\Phi_+$) is the unipotent radical of
$B(F).$ When $F$ is a local field, $G(F)$ is a locally compact
topological group. We choose a maximal compact subgroup $K$ of
$G(F)$ as follows. We first choose for $SL_2(F)$ a maximal compact
subgroup. If $F=\mb{R}$ or $\mb{C}$ we choose $SO_2(\mb{R})$ or
$SU_2(\mb{C})$. If $F$ is non-archemedean, we choose $SL_2({\cal
O}_F)$ where ${\cal O}_F$ is the ring of integers of $F$. Using
(\ref{cpt}) we obtain a maximal compact subgroup in the $SL_2(F)$
corresponding to each positive root $\a.$ Let $K$ be the subgroup
generated by these subgroups. Then we have the Iwasawa decomposition
$G(F)=B(F)K.$

Let $\wh{B}_0$ be the preimage of $B(F)$ of the canonical
projection $G(F[[t]])\ra G(F).$ It is easy to prove that $\wh{B}_0$
is generated by the elements $x_\a(u)$ where either $\a\in\Phi_+$, $u\in
F[[t]]$ or $\a\in \Phi_-, u\in tF[[t]],$ and the elements $h_\a(u)$, $\a\in\Phi,$
$u\in F[[t]]^\times.$ The subgroup $\wh{B}_0$
plays the role of Borel subgroup for $G(F((t)))$. Let $\wh{N}_0$ be
the group generated by $w_\a(u)$ with $\a\in\Phi,$ $u\in
F((t))^\times.$

\begin{lemma}\label{bn1}
The subgroup $\wh{H}_0=\wh{B}_0\cap \wh{N}_0$ is generated by
elements $h_\a(u)$ where $\a\in\Phi,$ $u\in F[[t]]^\times$, and is
normal in $\wh{N}_0.$
\end{lemma}

Let $w_\a=w_\a(1),$ $\wtl{w}_\a=\wtl{w}_\a(1).$ Let $\wh{S}_0=\{w_{\a_0}\wh{H}_0,\ld, w_{\a_n}\wh{H}_0\}\ss
\wh{N}_0/\wh{H}_0.$

\begin{thm}\label{tits1}
$(G(F((t))), \wh{B}_0, \wh{N}_0, \wh{S}_0)$ is a Tits system, and
its Weyl group is isomorphic to $\wtl{W}$. Moreover
$w_{\a_i}\wh{H}_0\mt r_{i}$ gives an isomorphism.
\end{thm}

Theorem \ref{tits1} follows from \cite[pp.37-38]{IM} Theorem 2.22 and 2.24.
Now consider the central extension $\wh{G}(F((t))).$

\begin{thm}\label{lift}
There exists a lifting $G(F[[t]])\ra \wh{G}(F((t)))$ given by
$x_\a(u)\mt \wtl{x}_\a(u)$, where $\a\in\Phi, u\in F[[t]].$
\end{thm}

The proof of Theorem \ref{lift} requires the theory of highest
weight representations of loop groups, and will be given in the next
section. Assume its validity at the moment, we may regard
$G(F[[t]])$ and its subgroups as subgroups of $\wh{G}(F((t))).$ In particular we may identify
$\wtl{U}_\b$ with $U_\b$ for each $\b\in\wtl{\Phi}_+.$ Let
$\wh{B}$ be the preimage of $\wh{B}_0$ in $\wh{G}(F((t)))$ under the
canonical map $\wh{G}(F((t)))\ra G(F((t))).$ Then $\wh{B}=\wh{B}_0\times
F^\times$ is a
subgroup of $G(F[[t]])\times F^\times$. Let $\wh{N}$ be the subgroup of $\wh{G}(F((t)))$
generated by the center $F^\times$ and the elements $\wtl{w}_\a(u)$
with $\a\in\Phi,$ $u\in F((t))^\times.$

\begin{lemma}\label{bn2}
The subgroup $\wh{H}=\wh{B}\cap\wh{N}$ is generated by the center
$F^\times$ and the elements $\wtl{h}_\a(u)$ with $\a\in\Phi,$ $u\in
F[[t]]^\times$, and is normal in $\wh{N}.$
\end{lemma}

Let $\wh{S}=\{\wtl{w}_{\a_0}\wh{H},\ld, \wtl{w}_{\a_n}\wh{H}\}\ss
\wh{N}/\wh{H}.$

\begin{thm}\label{tits2}
$(\wh{G}(F((t))),\wh{B},\wh{N},\wh{S})$ is a Tits system, and its
Weyl group is isomorphic to $\wtl{W}$ under the isomorphism given by
$\wtl{w}_{\a_i}\wh{H}\mt r_{i}.$
\end{thm}

Lemma \ref{bn2} and Theorem \ref{tits2} are immediate consequences
of Lemma \ref{bn1} and Theorem \ref{tits1}. We shall always identify
the quotient $\wh{N}/\wh{B}\cap\wh{N}$ with the affine Weyl group
$\wtl{W}$ using the isomorphism in the theorem.

\begin{lemma}\label{rep}
For every $\a\in\Phi,$ $\wtl{h}_\a(t^i)\in\wh{N}$ and it maps to
$T_{-i\a^\vee}\in\wtl{W}.$
\end{lemma}
\textit{Proof.} $\wtl{h}_\a(t^i)=\wtl{w}_\a(t^i)\wtl{w}_\a^{-1}\in
\wh{N}.$ The isomorphism in Theorem \ref{tits2} maps
$\wtl{w}_{\de-\a}\wtl{w}_\a$ to $r_{\de-\a}r_\a=T_{\a^\vee}.$ On the
other hand, \bq
\wtl{w}_{\de-\a}\wtl{w}_\a&=& \wtl{w}_{-\a}(t)\wtl{w}_\a= \wtl{w}_\a(-t^{-1})\wtl{w}_\a\n\\
&=& \wtl{h}_\a(-t^{-1})\wtl{w}_\a^2= \wtl{h}_\a(-t^{-1}) \wtl{h}_\a(-1)= \wtl{h}_\a(t^{-1}),\n \eq where the 2nd equality used Lemma \ref{ma} ($f$), the 2nd last equality
used \cite[pp.34-35]{Ma} Th\'{e}or\`{e}me 6.3 ($b$), and the last equality used Corollary \ref{r}
($b$). Therefore $\wtl{h}_\a(t^{-1})$ corresponds to
$T_{\a^\vee}\in \wtl{W}.$ \hfill$\Box$

\

The standard results about Tits system implies the Bruhat
decomposition \be\label{bru} \wh{G}(F((t)))=\bigcup_{w\in
\wtl{W}}\wh{B}w\wh{B}. \ee The general Bruhat decomposition with
respect to parabolic subgroups also applies to the loop groups,
where the notion of parabolic subgroups are explained below.

For any $\t\ss\wtl{\D}=\{\a_0,\a_1,\ld,\a_n\},$ there corresponds to
a parabolic subgroup $P_\t$ of $\wh{G}(F((t)))$ such that
$P_{\t_1}\ss P_{\t_2}$ if and only if $\t_1\ss\t_2.$ Let $P_\t=M_\t
N_\t$ be the Levi decomposition where $M_\t$ is the Levi subgroup,
and $N_\t$ is the unipotent radical.

For example, we have
$P_{\wtl{\D}}=\wh{G}(F((t)))$, $P_\emptyset=\wh{B},$
$M_{\emptyset}=T\times F^\times$ where $T\h\mb{G}_m^n$ is generated by
$\wtl{h}_\a(u)$ with $\a\in\Phi$ and $u\in F^\times$. An important example is the maximal parabolic subgroup
$P_{\D}=G(F[[t]])\times F^\times$ with $M_{\D}=G(F)\times F^\times.$
In fact these are the subgroups from which we induce the Eisenstein
series in section 4.

Let $\wh{U}=N_\emptyset$ be the unipotent radical of
$\wh{B}=P_\emptyset$, i.e. $\wh{U}$ is generated by the elements
$\wtl{x}_\a(u)$ where either $\a\in\Phi_+$, $u\in F[[t]]$ or
$\a\in\Phi_-,$ $u\in tF[[t]].$ Let $U^+$ be the subgroup
generated by the elements $\wtl{x}_\a(u)$ where $\a\in\Phi_+$ and
$u\in F[[t]]$, and $U^-$ be the subgroup generated by the
elements $\wtl{x}_\a(u)$ where $\a\in\Phi_-$ and $u\in tF[[t]],$ and
$D=\wh{B}\cap \wh{N}$ be the subgroup generated by the center $F^\times$ and the elements
$\wtl{h}_\a(u)$ where $\a\in\Phi$, $u\in F[[t]]^\times$,  and
$D^1$ be the subgroup of $D$ generated by $\wtl{h}_\a(u)$
where $\a\in\Phi$ and $u\in 1+tF[[t]]$. Let $\wh{T}=T\times
F^\times$, then both $\wh{T}$ and $D^1$ are stable under the
conjugation of $\wtl{W}.$

\begin{lemma} $($\cite[p.29]{IM} Proposition 2.1$)$
We have unique factorizations
\[
D^1=\prod_{\b\in \wtl{\Phi}_{im+}}U_\b,\q U^-U^+=\prod_{\b\in\wtl{\Phi}_{re+}}U_\b,\q D=\wh{T}D^1,
\]
\[
\wh{B}=\wh{T}\wh{U}=U^-DU^+,\q
\wh{U}=U^-D^1U^+=\prod_{\b\in\wtl{\Phi}_+}U_\b.
\]
\end{lemma}

In general let $N^{\pm}_\t=U^\pm\cap N_\t,$ then
\be\label{rs}
N_\t=N^-_\t D^1 N^+_\t.
\ee
Note that if $\Phi_\t$ is the subsystem of $\wtl{\Phi}$ generated by
$\t,$ then
\be\label{rs1}
N_\t=\prod_{\a\in \wtl{\Phi}_+-\Phi_\t}U_\a.
\ee
Let $W_\t$ be the subgroup of $\wtl{W}$ generated by $\{r_i|\a_i\in
\t\}$. The following result is also standard, see \cite{cas}.

\begin{thm}\label{bruhat}
For $\t_1,$ $\t_2\ss\wtl{\D},$ there is the Bruhat decomposition into disjoint unions
\[
\wh{G}(F((t)))=\bigcup_{W_{\t_1}\bs
\wtl{W}/W_{\t_2}}P_{\t_1}wP_{\t_2},
\]
where $w$ runs over a set of representatives of the double cosets in
$W_{\t_1}\bs \wtl{W}/W_{\t_2}$. The following is such a set of
double coset representatives:
\[
W(\t_1,\t_2)=\{w\in\wtl{W}| w^{-1}\t_1\ss\wtl{\Phi}_+,
w\t_2\ss\wtl{\Phi}_+\}.
\]
In the case $\theta_1=\theta$ and $\theta_2=\emptyset$, for each $w\in W(\t,\emptyset)$ there is a bijection
\[
P_{\t}wP_\emptyset\h P_{\t}\times \{w\}\times U_w,
\]
where
\[
U_w=\prod_{\a>0,
w\a<0}U_\a.
\]
\end{thm}

Assume now $F$ is a local field. For a real root
$\b\in\wtl{\Phi}_{re}$, we denote $K_\b$ the image of the standard
maximal compact subgroup of $SL_2(F)$ under the map $\vp_\b$ in
Lemma \ref{sl2}. Let $\wh{K}$ denote the subgroup of
$\wh{G}(F((t)))$ generated as follows: \be\label{K}
\wh{K}=\l\{\ba{ll} \lg K_\b, \b\in \wtl{\Phi}_{re}, \pm 1\in \mb{R}^\times\rg, & \rm{ if }F=\mb{R},\\
\lg K_\b, \b\in \wtl{\Phi}_{re}, S^1\ss\mb{C}^\times\rg, & \rm{ if }F=\mb{C},\\
\lg K_\b, \b\in \wtl{\Phi}_{re}, {\cal O}_F^\times, G({\cal
O}_F[[t]])\rg, & \rm{ if }F\rm{ is }p\rm{-adic}. \ea\r. \ee The
standard method using Tits system shows that there is the Iwasawa
decomposition
\[
\wh{G}(F((t)))=\wh{B}\wh{K}.
\]
For all local field $F$ and all real root $\b,$
$\wtl{w}_\b=\wtl{w}_\b(1)\in \wh{K},$ therefore $\wtl{W}$ has a set of
representatives in $\wh{K}.$ Denote the image of $\wh{K}$ in
$G(F((t)))$ by $\wh{K}_0.$

Now let us construct the full loop group $\wtl{G}(F((t)))$. The
reparametrization group of $F((t))$ is
\[
\rm{Aut}_F F((t))=\l\{\sum^\i_{i=1}u_it^i\in F[[t]]|u_1\neq 0\r\},
\]
where $\s(t)\in \rm{Aut}_F F((t))$ acts on $F((t))$ by $u(t)\mt
u(\s(t)),$ and the group law is $(\s_1\ast\s_2)(t)=\s_2(\s_1(t)).$
This induces an action of $\rm{Aut}_FF((t))$ on $G(F((t)))$ as
automorphisms. It is easy to check that the action of
$\rm{Aut}_FF((t))$ preserves the tame symbol, therefore it acts on
$\wh{G}(F((t)))$ as automorphisms. More precisely, we have
\[
\s(t)\c \wtl{x}_\a(u(t))=\wtl{x}_\a(u(\s(t))),
\]
and the action on the center $F^\times$ is trivial. It is also clear
that the subgroup $G(F)$ is fixed under this action. We have the
semi-direct product group
\[
\wh{G}(F((t)))\rtimes \rm{Aut}_FF((t))
\]
on which there is the standard relation
\[
\s(t)g\s^{-1}(t)=\s(t)\c g.
\]
We shall only consider the subgroup $\s(F^\times )\ss
\rm{Aut}_FF((t))$ which consists of the elements $\s(\mf{q})=\mf{q} t$ ($\mf{q}\in
F^\times$). It is clear that $\s(F^\times )$ is isomorphic to
$\mb{G}_m(F).$ We form the semi-direct product group \be
\label{loop} \wtl{G}(F((t)))=\wh{G}(F((t)))\rtimes \s(F^\times). \ee

Consider the tori in $\wtl{G}(F((t))),$ \be\label{tori} T\hr
\wh{T}\hr \wtl{T}, \ee where $\wtl{T}=\wh{T}\times \s(F^\times).$  The torus
$\wtl{T}\h\mb{G}_m^{n+2}$ will play the role of a maximal torus for
$\wtl{G}(F((t))).$ Then we have the cocharacter lattices
\begin{displaymath}
\xymatrix { X_\ast(T)\ar[r]^\ss \ar[d]^\h & X_\ast(\wh{T})\ar[r]^\ss \ar[d]^\h & X_\ast(\wtl{T}) \ar[d]^\h\\
Q^\vee \ar[r]^\ss & Q^\vee\op \mb{Z}c \ar[r]^\ss &
Q^\vee_{\rm{aff}}}
\end{displaymath}
where $Q^\vee_{\rm{aff}}=Q^\vee\op\mb{Z}c\op\mb{Z}d$ is called the
affine coroot lattice. It has a basis
$\{\a_1^\vee,\ld,\a_n^\vee,c,d\}$, and
$\{\a_0^\vee,\a_1^\vee,\ld,\a^\vee_n,d\}$ is also a basis. The
identification of  $Q^\vee_{\rm{aff}}$ with $X_{\ast}(\wtl{T})$ is,
for $\la=\a^\vee+ic+jd\in Q^\vee_{\rm{aff}},$ where $\a\in\Phi,$
$i,j\in\mb{Z}$, the corresponding cocharacter is $\la: \mb{G}_m\ra
\wtl{T}$ given by
\[
\la(u)=\wtl{h}_\a(u)u^i\s(u^j).
\]

It is clear that $\wh{N}$ normalizes $\wtl{T}$,
and therefore $\wtl{W}$ acts on $\wtl{T}.$ On the other hand $\wtl{W}$ acts on $\wtl{\mf{h}}$ by the formula in Lemma
\ref{weyl1}, and the lattice $Q^\vee_{\rm{aff}}$ is stable under the
action. We have

\begin{lemma}\label{equiv}
The cocharacter map
\[
Q^\vee_{\mathrm{aff}}\times \mb{G}_m\ra \wtl{T},\q (\la,u)\mt \la(u)
\]
is $\wtl{W}$-equivariant.
\end{lemma}
\textit{Proof.} The lemma is equivalent to that, for every
$w\in\wtl{W}$, let $\wtl{w}\in\wh{N}$ be a representative, then
\be\label{equi} \wtl{w}\la(u)\wtl{w}^{-1}=(w\c\la)(u). \ee It is clear
that (\ref{equi}) is true for $w\in W.$ We now prove it for
$w=T_{-\a^\vee}.$ By Lemma \ref{rep}, $\wtl{w}=\wtl{h}_\a(t)$ is a
lifting of $w.$ For $\la=\b^\vee\in Q^\vee,$ (\ref{equi}) is a
special case of the following identity
\[
\wtl{h}_\a(t^i)\wtl{h}_\b(u)\wtl{h}_\a(t^i)^{-1}=u^{i\left(\a^\vee,\b^\vee\right)}\wtl{h}_\b(u),
\]
which follows from Lemma \ref{symb} ($c$). It remains to prove
(\ref{equi}) for $\la=d,$ for which the left-hand-side of
(\ref{equi}) is \bq
\label{lhs}\wtl{h}_\a(t)\s(u)\wtl{h}_\a(t)^{-1}&=& \s(u)\wtl{h}_\a(u^{-1}t)\wtl{h}_\a(t)^{-1}\\
&=&\s(u)u^{-\f{2}{\left(\a,\a\right)}}\wtl{h}_\a(u^{-1}),\n \eq
where we have used (\ref{RC}). By Lemma \ref{weyl1},
\[
T_{-\a^\vee}(d)=-\a^\vee-\f{\left(\a^\vee,\a^\vee\right)}{2}c+d=-\a^\vee-\f{2}{\left(\a,\a\right)}c+d.
\]
Therefore $T_{-\a^\vee}(d)(u)$ is equal to the right-hand-side of
(\ref{lhs}). \hfill$\Box$

\

Finally, define $\wtl{B}=\wh{B}\rtimes\s(F^\times).$ If $F$ is a
local field, we also define \[\wtl{K}=\wh{K}\rtimes\s({\cal M}_F),
\] where ${\cal M}_F$ is the maximal compact subgroup of
$F^\times,$ i.e.
\begin{equation}\label{mf}
{\cal M}_F=\l\{\ba{ll} \{\pm 1\}, & \rm{if }F=\mb{R},\\
S^1, & \rm{if }F=\mb{C},\\
{\cal O}_F^\times, & \rm{if }F\rm{ is }p\rm{-adic}.\ea\r.
\end{equation}

\subsection{Highest weight representations of loop groups}

Let $\la\in\wtl{\mf{h}}^\ast$ be a dominant integral weight, i.e.
$\lg\la,\a_i^\vee\rg\in\mb{Z}_{\geq 0}$, $i=0,\ld,n$ and $\lg\la,d\rg\in\mb{Z}$.
Let $V_\la$ be the corresponding irreducible highest weight
representation of $\wtl{\mf{g}}$, and $v_\la$ be a highest weight
vector. A vector $v\in V_\la$ is said to be homogeneous of weight
$\mu$ if it lies in a weight space $V_{\la,\mu}.$ Every vector $v\in
V_\la$ is a sum of homogeneous elements (called components of $v$).
There is a lattice $V_{\la,\mb{Z}}\ss V_\la$ which is preserved by
the action of $\df{1}{j!}(X_\a\ot t^i)^j\in {\cal U}(\wtl{\mf{g}})$
for every positive integer $j$ and basis vector $X_\a\in\mf{g}_\a$
in the Chevalley basis of $\mf{g}.$ Moreover,
\[
V_{\la,\mb{Z}}=\bop_\mu V_{\la,\mu,\mb{Z}}
\]
where $\mu$ runs over all the weights of $V_\la,$ and
$V_{\la,\mu,\mb{Z}}=V_{\la,\mb{Z}}\cap V_{\la,\mu}.$ Assume that
$V_{\la,\la,\mb{Z}}=\mb{Z}v_\la.$

For any local field $F$, $V_{\la,F}=V_{\la,\mb{Z}}\ot_\mb{Z}F$ is a
representation of $\wtl{G}(F((t)))$, with the action of
$\wtl{x}_\a(ut^i)$ given by
\[
\wtl{x}_\a(ut^i)v=\sum^\i_{j=0}\f{1}{j!}u^j(X_\a\ot t^i)^j v
\]
for every $v\in V_{\la,F}.$ Since $(X_\a\ot t^i)^jv=0$ for $j$ large
enough, the sum above is finite. Since the operators $X_\a\ot t^i$
($i\in\mb{Z}$) are commutative, and $(X_\a\ot t^i)v=0$ for $i$ large
enough, for $u=\sum\limits_{i=N}^\i u_i t^i\in F((t))$ the product
$\prod\limits^\i_{i=N}\wtl{x}_\a(u_it^i)v$ is finite, and we define
the action of $\wtl{x}_\a(u)$ as
$\prod\limits^\i_{i=N}\wtl{x}_\a(u_it^i).$ The action can be extended
to an action of $\wtl{G}(F((t)))$ by setting
\be\label{sac}
\s(\mf{q})v=\mf{q}^{\lg\mu,d\rg}v
\ee
for each $v\in V_{\la,\mu,F}.$

\begin{thm}\label{repn}$($Garland \cite{Gar1}$)$
There is an action of $\wtl{G}(F((t)))$ on $V_{\la,F}$ defined as
above. The action of $u\in F^\times$ on $V_{\la,F}$ is the scalar
$u^{\lg\la,c\rg},$ and the action of $\wtl{h}_\a(u)$ $(\a\in\Phi)$ on
$V_{\la,\mu,F}$ is $u^{\lg\mu,\a^\vee\rg}.$
\end{thm}

If $F=\mb{R}$ or $\mb{C}$, there is a hermitian inner product $(,)$
on $V_{\la,F}$ such that
\\
($i$) $(v_\la,v_\la)=1,$
\\
($ii$) homogeneous vectors with different weights are orthogonal,
\\
($iii$) there is a homogeneous orthonormal basis contained in
$V_{\la,\mb{Z}}$,
\\
($iv$) elements of $\wh{K}$ act as unitary operators.
\\
($v$) $X_\a\ot t^i$ and $X_{-\a}\ot t^{-i}$ are adjoint operators.
\\
In particular, the norm $\|v\|=(v,v)^{\f{1}{2}}\geq 1$ for all $v\in V_{\la,\mb{Z}}$, $v\neq
0.$ If $F$ is a $p$-adic field with ${\cal O}_F$ the ring of
integers and $\pi\in {\cal O}_F$ a uniformizer, we let
\[
V_{\la,{\cal O}_F}=V_{\la,\mb{Z}}\ot {\cal O}_F
\]
and define a norm on
$V_{\la,F}$ by setting $\|0\|=0$ and $\|v\|=|\pi^l|$ for $v\neq 0$, where $l$ is the
largest integer such that $v\in \pi^l V_{\la,{\cal O}_F}.$ Recall that the normalized absolute value
on $F$ is defined by $|\pi|=q^{-1}$ where $q$ is the cardinality of the residue field ${\cal O}_F/\pi{\cal O}_F$. Since the
action of $\wh{K}$ preserves $V_{\la,{\cal O}_F}$, this norm is
preserved by $\wh{K}.$ We also have $\|xv\|=|x|\|v\|$ for $x\in F,$
$x\in V_{\la,F}$, and $\|v_1+v_2\|\leq\max(\|v_1\|,\|v_2\|).$

Now we are ready to prove Theorem \ref{lift}. Let ${\cal H}$ be the
subgroup of $\wh{G}(F((t)))$ generated by the elements
$\wtl{x}_\a(u)$ with $\a\in\Phi,$ $u\in F[[t]].$ Write $\pi$ for the
projection $\wh{G}(F((t)))\ra G(F((t))).$ It suffices to prove the
following lemma.

\begin{lemma}\label{onto}
$\pi: {\cal H}\ra G(F[[t]])$ is an isomorphism.
\end{lemma}
\textit{Proof.} Using ${\cal U}(\wtl{\mf{g}})={\cal
U}_-(\wtl{\mf{g}}){\cal U}_+(\wtl{\mf{g}}),$ we have a decomposition
\[
V_\la= V_\la(0)\op V_\la(1)\op\cs,
\]
where
\[
V_\la(d)=\rm{Span}\{(X_{\a_1}\ot t^{-d_1})\cs (X_{\a_l}\ot
t^{-d_l})v_\la| d_i\geq0, \sum^l_{i=1}d_i=d\}.
\]
It is clear that $V_\la(0)$ is the highest weight module of $\mf{g}$
with highest weight $\la|_{\mf{h}}.$ Then $V_\la(0)$ is a
representation of $G(F)$ and becomes an ${\cal H}$-module via the
following diagram
\begin{displaymath}
\xymatrix{ {\cal H}\ar[r]^\pi \ar[dr]  & G(F[[t]]) \ar[d]^{t=0} &
\\ & G(F) \ar[r] & GL(V_\la(0)) }
\end{displaymath}
Assume $u\in\rm{Ker}(\pi|_{\cal H})\ss F^\times,$ then $u$ acts on
$V_{\la}(0)$ trivially. On the other hand by Theorem \ref{repn} $u$
acts by the scalar $u^{\lg\la,c\rg}.$ Therefore $u^{\lg\la,c\rg}=1$ for any dominant integral weight $\la.$ It follows that
$u=1.$ \hfill$\Box$

\

Using Theorem \ref{repn} we can also prove the following lemma for
$p$-adic loop groups.

\begin{lemma}
If $F$ is $p$-adic, then
\[
\mathrm{Ker}(\pi|_{\wh{K}})={\cal O}_F^\times.
\]
\end{lemma}
\textit{Proof.} If $u\in \rm{Ker}(\pi|_{\wh{K}})\ss F^\times,$ since
$\wh{K}$ preserves $V_{\la,{\cal O}_F}$ we obtain $u^{\lg\la,c\rg}V_{\la,{\cal O}_F}=V_{\la, {\cal O}_F}$ by Theorem
\ref{repn}. It follows that $u^{\lg\la,c\rg}\in{\cal O}_F^\times$
for any dominant integral weight $\la.$ This implies $u\in {\cal
O}_F^\times.$\hfill$\Box$

\subsection{Second construction of loop groups}

We start from the example $G=GL_n.$ A \textit{lattice} $L$ of an
$n$-dimensional $F((t))$-vector space $V$ is a free
$F[[t]]$-submodule of rank $n$. In other words $L$ is a
$F[[t]]$-span of a basis of $V$. Any two lattices $L_1$, $L_2$ in
$V$ are \textit{commensurable}, which means that the quotients
$L_1/(L_1\cap L_2)$ and $L_2/(L_1\cap L_2)$ are finite-dimensional
over $F$. For example any lattice in $F((t))^n$ is commensurable
with $F[[t]]^n.$

Let $L_0$ be the lattice $F[[t]]^n,$ and $g\in GL_n(F((t))).$ Since
$gL_0/(L_0\cap gL_0)$ is finite-dimensional over $F$, we can define
the top wedge power $\wedge^{top} (gL_0/L_0\cap gL_0)$, which is a
one-dimensional vector space over $F.$ Let $\det(L_0,gL_0)$ be the
tensor product
\[
\wedge^{top} (gL_0/L_0\cap gL_0)\ot_F \wedge^{top} (L_0/L_0\cap
gL_0)^{-1},
\]
where $(-)^{-1}=\rm{Hom}_F(-,F)$ denotes the dual vector space. And
let $\det(L_0,gL_0)^\times$ be the set of nonzero vectors in
$\det(L_0,gL_0)$, which form a $F^\times$-torsor. Now define the
group
\[
\wh{GL}_n(F((t)))_{\rm{st}}=\{(g, \om_g)| g\in GL_n(F((t))),
\om_g\in \det(L_0, gL_0)^\times\}.
\]
Here the subscript ``st" stands for the standard representation of
$GL_n.$ The multiplication in the group is given by
\[
(g,\om_g)(h,\om_h)=(gh, \om_g\wedge g\om_h),
\]
where $g\om_h$ is image of $\om_h$ under the natural map $\det(L_0,
hL_0)\stackrel{g}{\ra} \det (gL_0, ghL_0)$, and $\om_g\wedge g\om_h$
is defined by the isomorphism $\det(L_0,gL_0)\wedge \det(gL_0,
ghL_0)\ra \det (L_0, ghL_0).$

$\s(\mf{q})\in \rm{Aut}_FF((t))$ with $\mf{q}\in F^\times$ preserves $L_0$,
hence induces the maps
  \[
  gL_0/L_0\cap gL_0\lra (\s(\mf{q})\c g)L_0/L_0\cap (\s(\mf{q})\c g)L_0,\q L_0/L_0\cap gL_0\lra L_0/L_0\cap (\s(\mf{q})\c g)L_0.
  \]
  Therefore $\s(\mf{q})$ induces the map
  $\det(L_0,gL_0)\ra \det(L_0,(\s(\mf{q})\c g)L_0)$,
 hence acts on the group $\wh{GL}_n(F((t)))_{\rm{st}}.$ We can form the semi-direct product group
 \[
 \wtl{GL}_n(F((t)))_{\rm{st}}=\wh{GL}_n(F((t)))_{\rm{st}}\rtimes \s(F^\times).
 \]

Suppose that $G$ is a linear reductive algebraic group over $F$, and $(\rho,V)$ is a (faithful) rational representation of $G$. Then
 $G(F((t)))$ acts on $V_{F((t))}=V\ot_F F((t))$. Let $V_0=V_{F[[t]]}=V\ot_F F[[t]]$, then $V_0$ is a lattice of $V_{F((t))}.$ Define the
 following loop group
 \[
 \wh{G}(F((t)))_\rho=\{(g,\om_g)| g\in G(F((t))), \om_g\in \det(V_0,\rho(g)V_0)^\times\}.
 \]
 The group law is defined similarly as the $GL_n$ case. It is clear that $\wh{G}(F((t)))_\rho$ is a central extension of $G(F((t)))$ by $F^\times.$
 $\s(\mf{q})$ acts on $V_{F((t))}$ as $\rm{id}\ot \s(\mf{q})$, which preserves $V_0.$ Therefore we
 can form the full loop group:
 \[
 \wtl{G}(F((t)))_\rho=\wh{G}(F((t)))_\rho\rtimes \s(F^\times).
 \]

From the construction it is seen that this notion only depends on
the equivalence class of $\rho.$ Namely, if $(\rho,V)$ and
$(\rho',V')$ are equivalent representations of $G$, then
$\wtl{G}(F((t)))_\rho$ and $\wtl{G}(F((t)))_{\rho'}$ are isomorphic.
To show this, let $f: V\ra V'$ be any intertwining linear bijection,
then $f(\rho(g)v)=\rho'(g)fv$ for all $v\in V,$ $g\in G.$ The action
of $f$ extends to $V_{F((t))}=V\ot_F F((t))$ by scalar extension.
Then $f(V_0)=V_0',$ $f(\rho(g)V_0)=\rho'(g)V_0',$ and $f$ induces an
isomorphism of $F^\times$-torsors $f_g:
\det(V_0,\rho(g)V_0)^\times\ra \det(V'_0,\rho'(g)V_0')^\times,$
which commutes with the action of $\s(F^\times).$ Let us identify
$f_g$ with a scalar in $F^\times.$ Then
\[
\wtl{f}: \wtl{G}(F((t)))_\rho\ra \wtl{G}(F((t)))_{\rho'}, \q
(g,\om_g)\rtimes \s(u)\mt (g, f_g\om_g)\rtimes \s(u)
\]
is a group isomorphism, which can be checked easily. Write $\hat{f}$
for the restriction of $\wtl{f}$ to $\wh{G}(F((t)))_\rho$, then the
the following diagram with exact rows is commutative
\begin{displaymath}
\xymatrix
{1 \ar[r] & F^\times \ar[d]^{\rm{id}} \ar[r] & \wh{G}(F((t)))_\rho \ar[d]^{\hat{f}} \ar[r]^{\pi_\rho} & G(F((t)))\ar[r] \ar[d]^{\rm{id}} & 1\\
1 \ar[r] & F^\times \ar[r] & \wh{G}(F((t)))_{\rho'}
\ar[r]^{\pi_{\rho'}} & G(F((t))) \ar[r] & 1}
\end{displaymath}
It is also clear that $\wtl{f}_2\circ \wtl{f}_1=\wtl{f_2\circ f_1}$ if
$V_1\stackrel{f_1}{\ra}V_2\stackrel{f_2}{\ra}V_3$ are equivalence of
representations of $G$.

\begin{lemma}\label{iso}
If $G$ is a connected and simply connected semisimple algebraic
group split over $F$, then the isomorphism $\hat{f}:
\wh{G}(F((t)))_\rho\ra \wh{G}(F((t)))_{\rho'}$ does not depend on
the choice of the intertwining map $f$.
\end{lemma}
\textit{Proof.} It is equivalent to prove that if $\rho=\rho'$ and $f:
V\ra V$ is an intertwining map, then $\hat{f}=\wh{\rm{id}}.$ In
other words, we have to show that $f_g:
\det(V_0,\rho(g)V_0)^\times\ra \det(V_0, \rho(g)V_0)^\times$ is the
identity map for all $g\in G(F((t))).$ We may assume that $V=m
V_\tau$ where $V_\tau$ is irreducible. Then Schur's lemma implies
that $f$ is given by an element $a_f\in GL_m(F).$ Recall that for
any two lattices $L_1, L_2$ in $V_{F((t))}$, we have the notion of
relative dimension
\[
\rm{dim}(L_1, L_2)=\rm{dim}_F\f{L_1}{L_1\cap
L_2}-\rm{dim}_F\f{L_2}{L_1\cap L_2}.
\]
It is easy to see that
\\
($i$) $\rm{dim}(hL_1, hL_2)=\rm{dim}(L_1, L_2)$ for any $h\in
GL(V_{F((t))}),$
\\
($ii$) $\rm{dim}(L_1, L_3)=\rm{dim}(L_1, L_2)+\rm{dim}(L_2, L_3)$.
\\
Since $f$ commutes with $\rho(g),$ we have $f_g=(\det
a_f)^{\rm{dim}(V_{\tau0}, \rho(g)V_{\tau0})},$ where
$V_{\tau0}=V_\tau\ot_F F[[t]].$ Therefore we only have to show that
 $\rm{dim}(V_{\tau0}, \rho(g)V_{\tau0})=0.$ We have the decomposition $G=P_\Delta Q^\vee P_\Delta$ since $G$ is simply connected. Using this together
  with the fact $\rho(P_\Delta)$ preserves $V_{\tau0},$ by ($i$) above we may assume that $g\in \wtl{W}$ and is mapped into $Q^\vee$, e.g.
   $g$ is a product of elements of the form $h_\a(t)$, $\a\in \Phi$. Using ($i$) and ($ii$) repeatedly we get the formula
 \[
 \rm{dim}(L, h_1h_2\cs h_n L)=\rm{dim}(L, h_1 L)+\cs +\rm{dim}(L, h_n L).
 \]
Hence we are reduced to prove that $\rm{dim}(V_{\tau0},
h_\a(t)V_{\tau0})=0.$ However this follows from the facts that
$h_\a(u)$ acts on the weight space $V_{\tau,\la}$ of $V_{\tau}$ of
weight $\la$ by the scalar $u^{\lg\la,\a^\vee\rg}$, $\lg r_\a \la,
\a^\vee\rg=-\lg\la,\a^\vee\rg$, and
$\rm{dim}V_{\tau,w\la}=\rm{dim}V_{\tau,\la},$ $\forall w\in W.$
\hfill$\Box$

\

Let us denote by $[\rho]$ the equivalence class of representations
of $\rho$. Under the condition of Lemma \ref{iso},
$\wh{G}(F((t)))_{\rho_1}$ and $\wh{G}(F((t)))_{\rho_2}$ are
canonically isomorphic for any $\rho_1, \rho_2\in [\rho]$, and our
loop group can be written as $\wh{G}(F((t)))_{[\rho]}.$

To see the relations with the construction in section 3.1, let us
assume that $G$ is a connected and simply connected simple linear
algebraic group split over $F$. Let $\mf{g}_F$
be the (simple) Lie algebra of $G$, and
$\mf{g}=\mf{g}_\mb{Z}\ot_\mb{Z}\mb{C}$ be the complex simple Lie
algebra, where $\mf{g}_\mb{Z}$ is the lattice spanned by a Chevalley
basis of $\mf{g}_F.$ Let $\wh{G}(F((t)))$ be the central extension
of $G(F((t)))$ constructed in section 3.1. Let $(\rho,V)$ be a rational
representation of $G$.

\begin{thm}\label{new}
There exists a group homomorphism $\phi_{\rho}: \wh{G}(F((t)))\ra
\wh{G}(F((t)))_{\rho}$ such that the following diagram is
commutative
\begin{displaymath}
\xymatrix{1 \ar[r] & F^\times \ar[d]^{d_\rho} \ar[r] & \wh{G}(F((t))) \ar[d]^{\phi_\rho} \ar[r]^{\pi} & G(F((t)))\ar[r] \ar[d]^{\mathrm{id}} & 1\\
1 \ar[r] & F^\times \ar[r] & \wh{G}(F((t)))_{\rho}
\ar[r]^{\pi_{\rho}} & G(F((t))) \ar[r] & 1}
\end{displaymath}
where $d_\rho$ is the Dynkin index of the representation $\rho$, and
$F^\times\stackrel{d_\rho}{\ra}F^\times$ is the $d_\rho$-th power.
\end{thm}

The Dynkin index of a representation, introduced to the theory of $G$-bundles over a curve by
Faltings \cite{Fa} and Kumar et. al. \cite{KNR}, is defined as follows. By abuse of notation we also write the representation
$\rho: \mf{g}\ra \mf{sl}(V).$ Let
\[
\rm{ch }V=\sum_\la n_\la e^\la
\]
be the formal character of $V$. Then the Dynkin index of $\rho$ is defined to be
\[
d_\rho=\f{1}{2}\sum_\la n_\la\lg\la,\wtl{\a}^\vee\rg^2.
\]
\cite{LS} contains a Lie algebra version of this theorem, which is much easier to prove. The minimal Dynkin index $d_{\mf{g}}$ is defined to be
$\min d_\rho$ where $\rho$ runs over all representations $\rho: \mf{g}\ra \mf{sl}(V).$ For a dominant weight $\lambda$, let $\rho_\lambda$ be the irreducible $\frak{g}$-module with highest weight $\lambda$. The following table is given in \cite{LS}.

\begin{tabular} {|c|c|c|c|c|c|c|c|c|c|}
\hline
Type of $\mf{g}$ & $A_n$ & $B_n$, $n\geq 3$ & $C_n$ & $D_n$, $n\geq 4$ & $E_6$ & $E_7$ & $E_8$ & $F_4$ & $G_2$ \\
\hline
$d_{\mf{g}}$ & 1 & 2 & 1 & 2 & 6 & 12 & 60 & 6 & 2\\
\hline
$\la$ s.t. $d_{\rho_\la}=d_{\mf{g}}$ & $\vpi_1$ & $\vpi_1$ & $\vpi_1$ &  $\vpi_1$ & $\vpi_6$ & $\vpi_7$ & $\vpi_8$ & $\vpi_4$ & $\vpi_1$ \\
\hline
\end{tabular}

\

\textit{Proof} of the Theorem. By the existence of canonical lifting (Theorem
\ref{lift})
\[
G(F[[t]])\hookrightarrow \wh{G}(F((t))),
\]
and that $\rho(G(F[[t]]))$ preserves $V_0,$ we first define
$\phi_\rho(\wtl{x}_\a(u))=(x_\a(u),1)$ for $\a\in \Phi, u\in F[[t]].$
For general $u\in F((t))$, choose $\b\in\Phi$ with
$\lg\a,\b^\vee\rg\neq 0.$ There exists $v\in F((t))$ such that
$v^{-\lg\a,\b^\vee\rg}u\in F[[t]].$
 Let $\overline{h}_\b(v)=(h_\b(v),\om)$ be an element in the preimage of $h_\b(v)$ under the projection $\wh{G}(F((t)))_\rho\ra \wh{G}(F((t))).$
We define \be\label{def} \phi_\rho(\wtl{x}_\a(u))=
\overline{h}_\b(v)\phi_\rho(\wtl{x}_\a(v^{-\lg\a,\b^\vee\rg}u))\overline{h}_\b(v)^{-1}.
\ee Let us check that this is well-defined, namely, the
right-hand-side of (\ref{def}) does not depend on the choice of
$\overline{h}_\b(v).$ But we have \bq
&&\overline{h}_\b(v)\phi_\rho(\wtl{x}_\a(v^{-\lg\a,\b^\vee\rg}u))\overline{h}_\b(v)^{-1}\label{omega}\\
&=&(h_\b(v),\om)(x_\a(v^{-\lg\a,\b^\vee\rg}u),1)(h_\b(v)^{-1}, h_\b(v)^{-1}\om^{-1})\n\\
&=& (x_\a(u), \om \wedge x_\a(u)\om^{-1}).\n \eq Then we
have to prove that for $\om'\in \det(V_0, h_{\b'}(v')V_0)^\times$
with $v'^{-\lg\a,\b'^\vee\rg}u\in F[[t]]$,
\[
\om\wedge x_\a(u)\om^{-1}=\om' \wedge x_\a(u)\om'^{-1}.
\]
Let $\eta=\om^{-1}\om'\in\det(h_\b(v)V_0, h_{\b'}(v')V_0)^\times.$
Since
\[
x_\a(u)h_\b(v)V_0=h_\b(v)x_\a(v^{-\lg\a,\b^\vee\rg} u)V_0=h_\b(v)V_0
\]
and similarly $x_\a(u)h_{\b'}(v')V_0=h_{\b'}(v')V_0,$ it is clear
that $x_\a(u)$ acts on the finite-dimensional spaces
\[
\f{h_\b(v)V_0}{h_\b(v)V_0\cap h_{\b'}(v')V_0},\q
\f{h_{\b'}(v')V_0}{h_\b(v)V_0\cap h_{\b'}(v')V_0 } \] unipotently.
Taking top wedge product we see that $x_\a(u)$ fixes $\eta.$ This
proves (\ref{def}) is well-defined.

To show that $\phi_\rho$ is the required homomorphism, we need to
verify that
\\
($i$) $\phi_\rho(\wtl{x}_\a(u))$, $\a\in\Phi,$ $u\in F((t))$ satisfy
(\ref{R1}) and (\ref{R2}) if rank $\frak{g}\geq 2$, or $(\ref{R1})$ and $(\ref{R2'})$ if $\frak{g}=\frak{sl}_2$.
\\
($ii$)
$\phi_\rho(\wtl{h}_\a(u))\phi_\rho(\wtl{h}_\a(v))\phi_\rho(\wtl{h}_\a(uv))^{-1}=\left(
u, v\right)_{tame}^{-\f{2}{\left(\a,\a\right)}d_\rho},$ $\a\in\Phi,$ $u, v\in
F((t))^\times.$
\\
(\ref{R1}) is trivial. (\ref{R2}) is clearly true when $u, v\in
F[[t]],$ and to prove the general case we need a lemma.

\begin{lemma}\label{conj}
Suppose that $\mf{g}$ is a simple complex Lie algebra with root
system $\Phi$, $\a,\b$ are positive roots in $\Phi$,  then there
exists $\g\in\Phi$ such that $\lg\a,\g^\vee\rg\lg\b,\g^\vee\rg>0.$
\end{lemma}

The lemma can be verified for the Lie algebra $\mf{g}$ of type $A, B, \ld, G$
separately. Since we cannot find this lemma in the literature, let
us sketch procedures of the proof.  If $(\a,\b)>0$ then it is obvious. If $(\a,\b)<0$, it is reduced to checking the lemma for all irreducible root
systems of rank 2, namely $A_2,$ $B_2$ and $G_2.$ For $\left(\a,\b\right)=0$ we only
have a case-by-case proof, and we omit the details here.

We continue to prove (\ref{R2}). Since $\a+\b\neq0,$ there exists
$w\in W$ such that $w\a, w\b>0.$ By Lemma \ref{conj}, we can find
$\g\in\Phi$ satisfying $\lg\a,\g^\vee\rg\lg\b,\g^\vee\rg>0.$ Then
there exists $a\in F((t))^\times$ such that
$a^{-\lg\a,\g^\vee\rg}u,$ $a^{-\lg\b,\g^\vee\rg}v\in F[[t]].$ As
before let $\overline{h}_\g(a)=(h_\g(a),\om)\in \wh{G}(F((t)))_\rho.$ It
follows that \bq
&&\phi_\rho(\wtl{x}_\a(u))\phi_\rho(\wtl{x}_\b(v))\phi_\rho(\wtl{x}_\a(u))^{-1}\phi_\rho(\wtl{x}_\b(v))^{-1}\n\\
&=&   \overline{h}_\g(a) \phi_\rho(\wtl{x}_\a(a^{-\lg\a,\g^\vee\rg}u))\phi_\rho(\wtl{x}_\b(a^{-\lg\b,\g^\vee\rg}v))\n\\
&&\quad \times\phi_\rho(\wtl{x}_\a(a^{-\lg\a,\g^\vee\rg}u))^{-1}\phi_\rho(\wtl{x}_\b(a^{-\lg\b,\g^\vee\rg}v))^{-1}\overline{h}_\g(a)^{-1}\n\\
&=& \overline{h}_\g(a)\prod_{i,j\in\mb{Z}^+,i\a+j\b\in\Phi}\phi_\rho(\wtl{x}_{i\a+j\b}(c_{ij}^{\a\b}(a^{-\lg\a,\g^\vee\rg}u)^i(a^{-\lg\b,\g^\vee\rg}v)^j)) \overline{h}_\g(a)^{-1}\n\\
&=&
\prod_{i,j\in\mb{Z}^+,i\a+j\b\in\Phi}\phi_\rho(\wtl{x}_{i\a+j\b}(c_{ij}^{\a\b}u^iv^j)).\n\eq
This proves (\ref{R2}). The same trick also applies to the proof of
(\ref{R2'}).

Let us compute the 2-cocycle and prove ($ii$) above. It suffices to
treat the case $\mf{g}=\mf{sl}_2.$ In fact, if we define
\[
e_\a=\f{1}{2}\sum_\la n_\la \lg\la,\a^\vee\rg^2
\]
for $\a\in \Phi,$ then $e_\a$ is proportional to $1/\left(\a,\a\right).$ So
assume $\mf{g}=\mf{sl}_2$, $\a$ is the simple root, and we shall
prove that
\[
\phi_\rho(\wtl{h}_\a(ut^i))\phi_\rho(\wtl{h}_\a(vt^j))\phi_\rho(\wtl{h}_\a(uvt^{i+j}))^{-1}=(ut^i,
vt^j)_{tame}^{-d_\rho},
\]
where $u, v\in F^\times,$ $i, j\in\mb{Z}.$ Let us restrict to the
case $i, j\geq 0.$ Other cases can be treated similarly. We can
further assume that $\rho$ is irreducible, say, of highest weight
$m$, then $V$ has the weight space decomposition
\[
V=\bop_{\la=m,m-2,\ld,-m}V_\la,
\]
and $d_\rho=\df{1}{2}\sum\limits_\la \la^2=\sum_{\la>0}\la^2.$ We
first study the element $\phi_\rho(\wtl{x}_{-\a}(ct^{-i}))$ where
$c\in F^\times,$ $i\in\mb{Z}_{\geq 0}$. By (\ref{def}) and (\ref{omega}) we have
\[
\phi_\rho(\wtl{x}_{-\a}(ct^{-i}))=(x_{-\a}(ct^{-i}),\om\wedge
x_{-\a}(ct^{-i})\om^{-1}),
\]
where $\om\in \det(V_0, h_\a(t^i)V_0)^\times.$ Let
$\{v_\la|\la=m,m-2,\ld,-m\}$ be a basis of $V$ such that
 $X_{-\a} v_\la=v_{\la-2}.$ Since $h_\a(t^i)$ acts on $V_{\la,F((t))}$ as the scalar $t^{i\la}$,
 we can choose  $\om=\bw\limits_{\la=m,m-2,\cs,-m} \om_\la,$ where
\be\label{om}
\om_\la=\l\{\ba{ll}\bw^{-1}_{l=-i\la}t^l v_\la,& \rm{if } \la<0,\\
\l(\bw^{i\la-1}_{l=0}t^l v_\la\r)^{-1},& \rm{if } \la>0.\ea\r. \ee
Recall that the action of $x_{-\a}(ct^{-i})$ is given by
\[
x_{-\a}(ct^{-i})v_\la=\sum^\i_{j=0}\f{c^j t^{-ij}}{j!}X_{-\a}^j
v_\la=\sum^\i_{j=0}\f{c^j t^{-ij}}{j!}v_{\la-2j},\] from which it is
seen that $x_\a(ct^{-i})$ preserves the $F$-span of
$\{t^lv_\la|\la<0, i\la\leq l\leq -1\}$ and acts unipotently. Write
$\om=\om_+\wedge \om_-,$ where $\om_+=\bw\limits_{\la>0}\om_\la,$
$\om_-=\bw\limits_{\la<0}\om_\la,$ with the $\la$'s in decreasing
order. Then $x_{-\a}(ct^{-i})$ fixes $\om_-,$ and
\[
\om\wedge x_{-\a}(ct^{-i})\om^{-1}=\om_+\wedge
x_{-\a}(ct^{-i})\om_+^{-1}.
\]                                                                                                                                                                                           Now we have
\bq
&&\phi_\rho(\wtl{h}_\a(ut^i))\n\\&=&\phi_\rho(w_\a(ut^i))\phi_\rho(w_\a(1))\n\\
&=& (x_\a(ut^i),1)\phi_\rho(\wtl{x}_{-\a}(-u^{-1}t^{-i}))(x_\a(ut^i),1)(w_\a(1),1)\n\\
&=& (x_\a(ut^i),1)(x_{-\a}(-u^{-1}t^{-i}),\om_+\wedge x_{-\a}(-u^{-1}t^{-i})\om_+^{-1})(x_\a(ut^i)w_\a(1),1)\n\\
&=& (h_\a(ut^i), \wtl{\om}),\n \eq where \bq
\wtl{\omega}&=& x_\a(ut^i)(\om_+\wedge x_{-\a}(-u^{-1}t^{-i})\om_+^{-1})\n\\
&=& x_\a(ut^i)\om_+\wedge x_\a(ut^i)x_{-\a}(-u^{-1}t^{-i})\om_+^{-1}\n\\
&=& \om_+\wedge w_\a(ut^i)x_\a(-ut^i)\om_+^{-1}\n\\
&=& \om_+\wedge w_\a(ut^i)\om_+^{-1}.\n \eq From the action of
$w_\a$ on $V$ it is seen that \bq
w_\a(ut^i)\om_+^{-1}&=&w_\a(ut^i)\bw_{\la>0}\bw^{i\la-1}_{l=0}t^l v_\la\n\\
&=& \bw_{\la>0}\bw^{i\la-1}_{l=0}(-1)^{\f{\la+m}{2}}u^{-\la}t^{l-i\la}v_{-\la}\n\\
&=& \bw_{\la>0}(-1)^{\f{i\la(\la+m)}{2}}u^{-i\la^2}\bw^{-1}_{l=-i\la}t^lv_{-\la}\n\\
&=&
\l(\prod_{\la>0}(-1)^{\f{i\la(\la+m)}{2}}u^{-i\la^2}\r)\bw_{\la>0}\om_{-\la}.\n
\eq In summary, we can write \be
\phi_\rho(\wtl{h}_\a(ut^i))=(h_\a(ut^i), \om_{u,i}), \ee where
\[
\om_{u,i}=\prod_{\la>0}(\v_\la^i
u^{-i\la^2})\bw_{\la>0}\om_{i,\la}\wedge\bw_{\la>0}\om_{i,-\la},
\]
$\v_\la=(-1)^{\f{\la(\la+m)}{2}},$ and $\om_{i,\la}$ is given by
(\ref{om}). Let us apply this to prove ($ii$). We have
\[\phi_\rho(\wtl{h}_\a(ut^i))\phi_\rho(\wtl{h}_\a(vt^j))=(h_\a(uvt^{i+j}),\om_{u,i}\wedge h_\a(ut^i)\om_{v,j})
\]
and \bq
&&\om_{u,i}\wedge h_\a(ut^i)\om_{v,j}\n\\
&=&\prod_{\la>0}(\v_\la^{i+j}u^{-i\la^2}v^{-j\la^2})\bw_{\la>0}\l(\bw^{i\la-1}_{l=0}t^lv_\la\r)^{-1}\bw_{\la>0}\l(\bw^{-1}_{l=-i\la}t^lv_{-\la}\r)\n\\
&&\q\q \wedge \bw_{\la>0}\l(\bw^{j\la-1}_{l=0}u^\la t^{i\la+l}v_\la\r)^{-1}\bw_{\la>0}\l(\bw^{-1}_{l=-j\la}u^{-\la}t^{-i\la+l}v_{-\la}\r)\n\\
&=&\prod_{\la>0}(u^{-j\la^2}v^{i\la^2})\prod_{\la>0}(\v_\la^{i+j}(xy)^{-(i+j)\la^2})\bw_{\la>0}\l(\bw^{i\la-1}_{l=0}t^lv_\la\r)^{-1}\bw_{\la>0}\l(\bw^{-1}_{l=-i\la}t^lv_{-\la}\r)\n\\
&&\q\q \wedge \bw_{\la>0}\l(\bw^{(i+j)\la-1}_{l=i\la}t^{l}v_\la\r)^{-1}\bw_{\la>0}\l(\bw^{-i\la-1}_{l=-(i+j)\la}t^{l}v_{-\la}\r)\n\\
&=& \prod_{\la>0}(u^{-j\la^2}v^{i\la^2}(-1)^{ij\la^2})\om_{uv,i+j}\n\\
&=& (xt^i, yt^j)_{tame}^{-d_\rho}\om_{uv,i+j},\n \eq where the
second last equality is obtained from a direct counting. This
finishes the proof of ($ii$),  hence $\phi_\rho$ is the required
homomorphism. \hfill$\Box$

\subsection{Adelic loop groups and arithmetic quotients}

Let $F$ be a number field, and for each place $v$ let $F_v$ be the
completion of $F$ at $v$. For each local field $F_v$, we have the
local loop groups $G(F_v((t)))$, $\wh{G}(F_v((t)))$ and
$\wtl{G}(F_v((t)))$ constructed in section 3.1 which correspond to a
complex simple Lie algebra $\mf{g}.$ We add the subscript $v$ to
indicate the corresponding local subgroups. So we have\[
\wh{B}_{0v}\hr\wh{B}_v\hr \wtl{B}_v,\quad \wh{K}_{0v}\hr\wh{K}_v\hr\wtl{K}_v,\quad
T_v\hr\wh{T}_v\hr \wtl{T}_v.
\]
 For example in
$\wtl{G}(F_v((t)))$ we have
$\wtl{T}_v=\wh{T}_v\times\s(F_v^\times).$ We form the restricted
direct product group $\prod'_v G(F_v((t)))$ (resp.
$\prod_v'\wh{G}(F_v((t))$, $\prod_v'\wtl{G}(F_v((t)))$) with respect
to $\wh{K}_{0v}$ (resp. $\wh{K}_v$, $\wtl{K}_v$)'s.

Let $\mb{A}$ and $\mb{I}$ be the adele ring and the idele group of $F$ respectively. We let
$\mb{A}\lg t\rg$ be the restricted product $\prod'_v F_v((t))$ with respect to ${\cal O}_v((t))$'s for finite places $v.$
In other words,
\[
\mb{A}\lg t\rg=\{(x_v)_v| x_v\in F_v((t)),\rm{ and }x_v\in{\cal O}_v((t))\rm{ for almost all finite places }v\}.
\]
Note that we do not require that $(x_v)$'s in $(x_v)_v$ have bounded poles, so the ring $\mb{A}\lg t \rg$ is not
a subring of $\mb{A}((t)).$ Let
\[
F\lg t \rg=F((t))\cap\mb{A}\lg t \rg,
\]
i.e. $F\lg t\rg$ is the subset of the elements $x\in F((t))$ such that $x\in {\cal O}_v((t))$ for almost all finite places $v$.
We also define
\[
F\lg t\rg_+=F\lg t\rg\cap F[[t]],\q \mb{A}\lg t\rg_+=\mb{A}\lg t\rg\cap \mb{A}[[t]].
\]

\begin{lemma}
$F\lg t\rg$ is a subfield of $F((t)).$
\end{lemma}
\textit{Proof.} The only thing we need to check is that if $x\in F\lg t\rg$ and $x\neq0$ then $x^{-1}\in F\lg t\rg.$
We can assume that $x=1+x_1t+x_2t^2+\cs,$ then the coefficients in $x^{-1}$ are polynomials of $x_n$'s. Therefore $x^{-1}$ also lies in
$F\lg t\rg.$
\hfill$\Box$

\

We shall denote $\prod'_v\wh{G}(F_v((t)))$ by $\wh{G}(\mb{A}\lg t\rg)$, and $\prod_v'G(F_v((t)))$ by $G(\mb{A}\lg t \rg).$ For $\a\in\Phi$, $u\in F((t))$, we also denote by $\wtl{x}_\a(u)$
 the element in $\prod_v\wh{G}(F_v((t)))$ whose $v$-component is $\wtl{x}_\a(u)$ in $\wh{G}(F_v((t)))$. If $u\in F\lg t \rg,$
  then $\wtl{x}_\a(u)\in \wh{G}(\mb{A}\lg t \rg).$ We denote the subgroup generated by $\wtl{x}_\a(u)$ ($\a\in\Phi, u\in F\lg t \rg$) and $G(F\lg t \rg_+)$ by
  $\wh{G}(F\lg t \rg).$ It is clear that $\wh{G}(F\lg t \rg)/F^\times$ is isomorphic to $G(F\lg t \rg).$ we have the diagram with exact rows
\begin{displaymath}
\xymatrix{1 \ar[r]  & F^\times \ar[r] \ar[d] & \wh{G}(F\lg t \rg) \ar[r] \ar[d] & G(F\lg t \rg) \ar[r] \ar[d] & 1\\
1 \ar[r]  & \mb{I} \ar[r] \ar[d] & \wh{G}(\mb{A}\lg t \rg) \ar[r] \ar[d] & G(\mb{A}\lg t \rg) \ar[r] \ar[d] & 1\\
1 \ar[r]  & \mb{I}/F^\times \ar[r]  & \wh{G}(\mb{A}\lg t \rg)/F^\times \ar[r]  & G(\mb{A}\lg t \rg) \ar[r]  & 1}
\end{displaymath}
where $F^\times \hr \mb{I}$ is the diagonal subgroup.

By abuse of notations, we also use $\wh{T}$, $\wh{B}, \ldots$ to denote the following adelic subgroups of $\wh{G}(\mb{A}\lg t \rg)$:
\bq
&& \wh{T}=\wh{T}(\mb{A})=\prod'_v \wh{T}_v,\quad \wh{U}=\wh{U}(\mb{A})=\prod'_v \wh{U}_v,\n\\
&& \wh{B}=\wh{B}(\mb{A})=\prod'_v \wh{B}_v,\quad \wh{K}=\wh{K}(\mb{A})=\prod_v \wh{K}_v,\n
\eq
where the restricted products are defined with respect to the corresponding analogues of maximal compact subgroups in the finite dimensional case.
For example, $\wh{U}$ is generated by $\wtl{x}_\a(u)$ where either $\a\in \Phi_+,$ $u\in \mb{A}\lg t \rg_+$ or $\a\in \Phi_-,$ $u\in t\mb{A}\lg t \rg_+.$
Then every element $g\in \wh{G}(\mb{A}\lg t \rg)$ can be wrtten as $g= u_g a_g k_g$ with $u_g\in\wh{U}$, $a_g\in\wh{T}$ and $k_g\in \wh{K}.$ The local actions
of $\s(F_v^\times)$ on $\wh{G}(F_v((t)))$ piece together to form an action of the group $\s(\mb{I})$ on $\wh{G}(\mb{A}\lg t \rg).$ Define the semi-direct product group
\be\label{lg}
\wtl{G}(\mb{A}\lg t \rg)=\wh{G}(\mb{A}\lg t \rg)\rtimes\s(\mb{I}).
\ee
Similarly we can set the subgroups of $\wtl{G}(\mb{A}\lg t \rg)$:
\bq
&& \wtl{T}=\wh{T}\times\s(\mb{I})=\prod'_v \wtl{T}_v,\quad \wtl{B}=\wh{B}\rtimes \s(\mb{I})=\prod_v'\wtl{B}_v,\n\\
&& \wtl{K}=\wh{K}\rtimes\prod_v\s({\cal M}_{F_v})=\prod_v \wtl{K}_v.\n
 \eq
Since $\s(\mb{I})$ normalizes $\wtl{B},$  any $g\in \wtl{G}(\mb{A}\lg t \rg)$ can be written as
$g= u_g a_g k_g$ with $u_g\in\wh{U}$, $a_g\in\wtl{T}$ and $k_g\in \wtl{K}.$

Let $\wtl{G}(F\lg t \rg)=\wh{G}(F\lg t \rg)\rtimes \s(F^\times)\hr \wtl{G}(\mb{A}\lg t \rg).$ By abuse of notation, for any subgroup $H$ of $\wtl{G}(F((t)))$, we still denote by $H$ the subgroup $H\cap \wtl{G}(F\lg t \rg)$ of $\wtl{G}(F\lg t \rg).$ For example we have the subgroup $\wh{U}(F)$
 of $\wh{G}(F\lg t \rg)$, which is generated by $\wtl{x}_\a(u)$, where either $\a\in \Phi,$ $u\in F\lg t \rg_+$ or $\a\in\Phi_-,$ $u\in t F\lg t \rg_+.$

By Lemma \ref{equiv} we have a $\wtl{W}$-equivarient map
\[
Q^\vee\times \mb{I}\ra \wtl{T}
\]
given by, for $\la=\a^\vee+ic+jd\in Q_{\rm{aff}}^\vee$ where $\a\in \Phi,$ $u\in\mb{I}$,
\[
(\la,u)\mt \la(u)= \wtl{h}_\a(u)u^i\s(u^j).
\]

For a dominant integral weight $\la,$ we have a representation of $V_{\la, F_v}$ of $\wtl{G}(F_v((t)))$ for each place $v$. Form the restricted product
\[
V_{\la,\mb{A}}=\prod_v'V_{\la,F_v}
\]
with respect to the lattices $V_{\la,{\cal O}_v}$ which are defined for all finite places. Denote by $v_\la\in V_{\la,\mb{A}}$ the element with
$v$-component $v_\la$ for each place $v.$ Note that $V_{\la, F}$ embeds diagonally into $V_{\la,\mb{A}}.$

We define a map $|\c|: V_{\la,\mb{A}}\ra \mb{R}_{\geq 0}$ as follows. Recall that we have defined a norm on $V_{\la,F_v}$ for each place $v$ in section 3.2.  For $(u_v)_v\in V_{\la,\mb{A}},$ if $v$ is real or $p$-adic let $|u_v|= \|u_v\|$; if $v$ is complex let $|u_v|=\|u_v\|^2.$ Then define $|(u_v)_v|=\prod_v|u_v|.$ Note that almost all $|u_v|$'s are less than or equal to 1, hence the product is finite. If $u\in V_{\la,\mb{A}}$ and $k\in\wtl{K}$, then $|ku|=|u|$; and if $x=(x_v)_v\in \mb{I}$, then $|xu|=|x||u|$ where $|x|=\prod_v|x_v|_v.$ In particular $|xu|=|u|$ for $x\in F^\times$ by the Artin product formula.

For $\mu\in\wtl{\mf{h}}^\ast,$ define a quasi-character $\mu: \wtl{T}(F)\bs\wtl{T}(\mb{A})\ra\mb{C}^\times$ by, for $g=\wtl{h}_\a(x)y\s(z)\in \wtl{T}$ where $x, y, z\in\mb{I},$
\[
\mu(g)= |x|^{\lg\mu,\a^\vee\rg}|y|^{\lg\mu,c\rg} |z|^{\lg\mu,d\rg}.
\]
In particular, for $\la\in Q_{\rm{aff}}^\vee$ we have
\[
\la(x)^\mu=|x|^{\lg\la,\mu\rg}.
\]
However sometimes we also use the following notation: if $a=(a_v)_v\in\wtl{T}(\mb{A})$, $\b\in X^\ast(\wtl{T})$, the character lattice of $\wtl{T}$ spanned by $\Phi$, $\de$ and $L$, then write
\[
a^\b=(a_v^{\b})_v\in \mb{I}.
\]
Interpretations of the notations we shall use depend on the situation and would not cause any confusion.

\begin{lemma}\label{hight}
For each $g\in \wtl{G}(\mb{A}\lg t \rg)$ with decomposition $g=u_g a_g k_g$, and $v_\la\in V_{\la,\mb{A}_F}$ be the highest weight vector as above,
then
\[
|g^{-1}v_\la|=a_g^{-\la}.
\]
\end{lemma}
\textit{Proof.} Note that $u_g v_\la=v_\la,$ and since $\wtl{K}$ preserves $|\c|$, it follows that
\[
|g^{-1}v_\la|=|a_g^{-1}v_\la|=a_g^{-\la}.
\]
\hfill$\Box$

In the rest of this section we collect some lemmas on the arithmetic quotients of loop groups based on \cite{Gar1, Zhu2}.

\begin{lemma}\label{fin}
For each $i=0,1,\ld,n$, and each integer $l>0,$ the set of weights of $V_\la$ of the form
\[
\la-\sum^{n}_{j=0}l_j\a_j
\]
with $l_i\leq l$, is finite.
\end{lemma}

\begin{lemma}
Let $a\in \wh{T}$ and $a\s(\mf{q})\in\wtl{T}$ with $\mf{q}\in\mb{I}$, $|\mf{q}|<1$, then there exists $w\in \wtl{W}$ such that
$(wa\s(\mf{q}))^{\a_i}\leq 1$ for all $i=0,1,\ld,n.$
\end{lemma}
\textit{Proof.} It is easy to see that the lemma can be reduced to $\wh{G}(\mb{R}((t)))$. And it follows from the well-known fact
that for each $h=h_0+ic+jd\in\wtl{\mf{h}}_{\mb{R}}$ with $j>0,$ there exists $w\in\wtl{W}$ such that $\left( wh,\a_i\right)\geq0$ for
$i=0,1,\ld, n.$
\hfill$\Box$

\begin{lemma} Assume the conditions of the previous lemma. Moreover suppose that $(a\s(\mf{q}))^{\a_i}\leq 1$ for $i=0,1,\ld,n.$ Then there exists $0\leq j\leq n$ such that
\[
(a\s(\mf{q}))^{\a_j}<1.
\]
\end{lemma}

\begin{lemma}\label{min}
For any $g\in\wh{G}(\mb{A}\lg t \rg)\rtimes\s(\mf{q})$ with $\mf{q}\in\mb{I}$, $|\mf{q}|<1$, there exists $\g_0\in \wh{G}(F\lg t \rg)$ such that
\[
|g\g_0v_\la|\leq |g\g v_\la|
\]
for all $\g\in\wh{G}(F\lg t \rg).$
\end{lemma}
\textit{Proof.} The lemma is essentially an adelic formulation of Lemma 17.15 in \cite{Gar1}. Write $g=ka\s(\mf{q})u.$ First note
that for any positive number $C$, we may choose a finite set of weights $w_C$ of $V_\la$ such that $|(a\s(\mf{q}))^\mu|>C$ for any weight $\mu$ of $V_\la$ which is not in $w_C$. Enlarge $w_C$ if necessary, we may assume that if $\mu\in w_C$, then all the weights of $V_\la$ with depth less than the depth of $\mu$ are also in $w_C.$ Recall that the depth of a weight $\la-\sum^n_{i=0}l_i\a_i$ is $\sum^n_{i=0}l_i$. Then Lemma \ref{fin} guarantees the finiteness of $w_C.$ Now
set $C=(a\s(\mf{q}))^\la$, and divide $\wh{G}(F((t)))$ into two parts
\[
\wh{G}(F((t)))=G_1\cup G_2,
\]
where $G_1$ consists of all the $\g$'s such that some component of $\g v_\la$ has weight not in $w_C,$ and $G_2$ consists of other $\g$'s, i.e. all the $\g$'s such that all components of $\g v_\la$ have weights in $w_C.$ If $\g\in G_1$, let $\mu\not\in w_C$ be a maximal weight of $\g v_\la,$ and $v'\neq 0$ be the $\mu$-component of $\g v_\la.$ It is clear that $v'\in V_{\mu,F}$, and the $\mu$-component of $a \s(\mf{q})u\g v_\la$ is $(a\s(\mf{q}))^\mu v'$. Then
\[
|g\g v_\la|=|a \s(\mf{q})u\g v_\la|\geq (a \s(\mf{q}))^\mu|v'|=(a\s(\mf{q}))^\mu>C.
\] For $\g\in G_2$, $\g v_\la$ lies in the finite dimensional $F$-space
$\sum_{\mu\in w_C} V_{\la,\mu, F}.$
Consequently $a\s(\mf{q})u\g v_\la$ lies in a finite dimensional $F$-space in $\sum_{\mu\in w_C} V_{\la,\mu,\mb{A}}$, so there exists
$\g_0\in G_2$ such that
\be \label{ineq}
|g\g_0 v_\la|\leq |g\g v_\la|
\ee
for all $\g\in G_2.$ In particular $e\in G_2$, hence by Lemma \ref{hight}
\[
|g\g_0 v_\la|\leq |gv_\la|=(a\s(\mf{q}))^\la=C.
\]
Therefore (\ref{ineq}) holds for all $\g\in \wh{G}(F\lg t \rg).$
\hfill$\Box$

\

Consider the partial order on $\Phi_+$ such that $\a<\b$ if $\b-\a$ is a sum of positive roots. Fix a total order on $\Phi_+$ which extends this partial order, and induce the corresponding order on $\Phi_-$ by
identifying $\Phi_-$ with $\Phi_+$ via $\a\mt-\a.$

\begin{lemma} \label{uadecom} We have unique factorizations
\[
\wh{U}(\mb{A})=U^+(\mb{A})D^1(\mb{A})U^-(\mb{A})=\prod_{\a\in\Phi_+}U_\a(\mb{A}\lg t \rg_+)\prod^n_{i=1}\wtl{h}_{\a_i}(1+t\mb{A}\lg t \rg_+)\prod_{\a\in\Phi_-}U_\a(t\mb{A}\lg t \rg_+),
\]
\[
\wh{U}(F)=U^+(F)D^1(F)U^-(F)=\prod_{\a\in\Phi_+}U_\a(F\lg t \rg_+)\prod^n_{i=1}\wtl{h}_{\a_i}(1+tF\lg t \rg_+)\prod_{\a\in\Phi_-}U_\a(tF\lg t \rg_+),
\]
where the product is taken with respect to above fixed orders on $\Phi_+$ and $\Phi_-$.
\end{lemma}

Let ${\cal D}\ss\mb{A}$ be a fundamental domain of $\mb{A}/F.$ We shall take
\be\label{D}
{\cal D}={\cal D}_\i\times \prod_{v<\i}{\cal O}_v,
\ee
where ${\cal D}_\i$ is a fundamental domain of
$(\prod\limits_{v|\i} F_v)/{\cal O}_F=(F\ot_{\mb{Q}}\mb{R})/{\cal O}_F$ constructed as follows.
Let $\om_1,\ld, \om_N$ be a basis of ${\cal O}_F$ over $\mb{Z},$ where $N=[F:\mb{Q}].$
The diagonal embedding ${\cal O}_F\hr \prod\limits_{v|\i}F_v$ identify ${\cal O}_F$ with a lattice in $F\ot_{\mb{Q}}\mb{R}=\mb{R}^N$. Let ${\cal D}_\i$ be the following subset of $\prod\limits_{v|\i}F_v,$
\[
{\cal D}_\i=\l\{\sum_{i=1}^N t_i\om_i| 0\leq t_i<1\r\}.
\]
For example if $F=\mb{Q}$ then ${\cal D}_\i=[0,1)$ is a fundamental domain of $\mb{R}/\mb{Z}.$
We define
\[
{\cal D}\lg t \rg=\l\{\sum_i u_i t^i\in\mb{A}\lg t \rg| u_i\in{\cal D}, \forall i\r\},\q
{\cal D}\lg t \rg_+={\cal D}\lg t \rg\cap \mb{A}\lg t \rg_+.
\]
Let $\wh{U}_{\cal D}$ be subset of $\wh{U}(\mb{A})$
\be \label{ud}
\wh{U}_{\cal D}=\prod_{\a\in\Phi_+}U_\a({\cal D}\lg t \rg_+)\prod^n_{i=1}\wtl{h}_{\a_i}(1+t{\cal D}\lg t \rg_+)\prod_{\a\in\Phi_-}U_\a(t{\cal D}\lg t \rg_+),
\ee
where the product is taken with respect to the order in Lemma \ref{uadecom}.

\begin{lemma}\label{fd}
Every $u\in\wh{U}(\mb{A})$ can be written as $u=\g_u u_{\cal D}$ $($or $u_{\cal D}\g_u)$ for some $\g_u\in\wh{U}(F)$ and $u_{\cal D}\in\wh{U}_{\cal D}.$
\end{lemma}

\section{Eisenstein Series and Their Coefficients}

Let $G$ be the Chevalley group associated to a complex simple Lie algebra $\mf{g}$, and $F$ be a number field.
In this section we construct the Eisenstein series $E(s,f,g)$
defined on $\wh{G}(\mb{A}\lg t \rg)\rtimes\s(\mf{q})$, where $s\in\mb{C}$
and $f$ is an unramified cusp form on $G(\mb{A}).$ We always make
the assumption that $\mf{q}\in \mb{I}$ and $|\mf{q}|>1.$ We establish the
absolute convergence of the constant terms of $E(s,f,g)$ along
$\wh{U}$, under the condition that $\mf{Re}s$ is large enough. The
proof makes use of the Gindikin-Karpelevich formula, which will be
used frequently. The same method gives the values of
constant terms and Fourier coefficients of $E(s,f,g)$ along
unipotent radicals of parabolic subgroups.

\subsection{Definition of Eisenstein series}

 Let $G({\mb{A}})$ be the restricted product $\prod_v'G(F_v)$ with respect to $K_v,$ where $K_v$ is a maximal compact subgroup of $G(K_v).$ If $v$ is finite we take $K_v$ to be $G({\cal O}_v)$. Then
 \[
 G(\mb{A})=\varprojlim_{S} \prod_{v\in S}G(F_v)\prod_{v\not\in S}K_v,
 \]
where the inverse limit is taken over all finite sets of places.

Let $f\in L^2(G(F)\bs G(\mb{A}))$ be an unramified cuspidal automorphic form, i.e.
\\
($i$) $f$ is invariant under right translation of $K=\prod_v K_v,$
\\
($ii$) $f$ is an eigenform for all $p$-adic Hecke operators,
\\
($iii$) $f$ is an eigenform for all invariant differential operators at all infinite places,
\\
($iv$) the constant term of $f$ along the unipotent radical of any parabolic subgroup of $G$ is zero, i.e.
\[
\int_{U_P(F)\bs U_P(\mb{A})}f(ug)du=0,
\]
where $P$ is any parabolic subgroup of $G$ and $U_P$ is the unipotent radical of $P$.

Associated to $f$ and $s\in\mb{C}$, we shall define a function $\wtl{f}_s$ on $\wtl{G}(\mb{A}\lg t \rg).$ Suppose $g\in \wtl{G}(\mb{A}\lg t \rg)$ has a
decomposition
\[
g=cp\s(\mf{q}) k,
\]
where $c, \mf{q}\in \mb{I}$, $p\in G(\mb{A}\lg t \rg_+),$ $k\in \wh{K}=\prod_v \wh{K}_v.$ Write $p_0$ for the image of $p$
under the projection
\[
G(\mb{A}\lg t \rg_+)\ra G(\mb{A}).
\] Then we define
\[
\wtl{f}_s(g)=|c|^sf(p_0).
\]
We have to check that $\wtl{f}_s$ is well-defined, namely, $\wtl{f}_s(g)$ does not depend on the choice of the decomposition of $g.$
In fact, if $cp\s(\mf{q})k=c'p'\s(\mf{q})k',$ then $c'^{-1}c \s(\mf{q})^{-1}\c(p'^{-1}p)=k'k^{-1}\in \wh{B}\cap \wh{K},$ which implies that
$|c'c^{-1}|=1,$ $p_0=p'_0k_0$ for some $k_0\in K.$ Note that $(\s(\mf{q})\c p)_0=p_0$. This proves that $\wtl{f}_s$ is well-defined since $f$ is right $K$-invariant. $\wtl{f}_s$ has the following invariance properties:

\begin{lemma}\label{invf}
$($i$)$ $\wtl{f}_s$ is right $\wtl{K}$-invariant,\\
$($ii$)$ $\wtl{f}_s$ is left $G(F\lg t \rg_+)$-invariant, and left $\s(\mb{I})$-invariant,\\
$($iii$)$ $\wtl{f}_s(cg)=|c|^s\wtl{f}_s(g),$ and $\wtl{f}_s$ is $F^\times$-invariant.
\end{lemma}
\textit{Proof.} ($i$) By definition $\wtl{f}_s$ is right $\wh{K}$-invariant. It is also $\wtl{K}$-invariant since $\s({\cal M}_{F_v})$ normalizes $\wh{K}_v$ for each place $v$, where ${\cal M}_{F_v}$ is given by (\ref{mf}). ($ii$) follows from the fact that $f$ is left $G(F)$-invariant, and that $(\s(\mf{q})\c p)_0=p_0$, as noted above. ($iii$) follows from the definition and Artin product formula.
\hfill$\Box$

\

Let $\wh{G}(F\lg t \rg_+)=G(F\lg t \rg_+)\times F^\times,$ then $\wtl{f}_s$ is left $\wh{G}(F\lg t \rg_+)$-invariant by the previous lemma. We define the Eisenstein series $E(s,f,g)$ on $\wtl{G}(\mb{A}\lg t \rg)$ by
\be \label{eis}
E(s,f,g)=\sum_{\g\in \wh{G}(F\lg t \rg_+)\bs \wh{G}(F\lg t \rg)}\wtl{f}_s(\g g).
\ee
It is clear that the right-hand-side is a countable sum, right $\wtl{K}$-invariant and left $\wh{G}(F\lg t \rg)$-invariant. $E(s,f,g)$ is also left $\s(F^\times)$-invariant, hence left $\wtl{G}(F\lg t \rg)$-invariant.
To see this, let $\mf{q}\in F^\times.$ By Lemma \ref{invf} ($ii$),
\[
E(s,f,\s(\mf{q})g)=\sum_\g \wtl{f}_s(\g \s(\mf{q})g)=\sum_\g \wtl{f}_s(\s(\mf{q}^{-1})\g \s(\mf{q}) g)=E(s,f,g).
\]
Note that in this case $\s(\mf{q})$ acts on $\wh{G}(F\lg t \rg)$ as automorphism and preserves $\wh{G}(F\lg t \rg_+).$

For completeness let us construct Eisenstein series induced from cusp forms on other parabolic subgroups. Let $P$ be a parabolic
subgroup of $\wh{G}$ with Levi decomposition $P=M_PN_P$. Then $M_P$ is a finite dimensional split reductive group. Let $f_{M_P}$ be an unramified cusp form on $M_P(\mb{A})$, and $\nu$
be an unramified quasi-character of $M_P(\mb{A})$. If $g\in\wtl{G}(\mb{A}\lg t\rg)$
decomposes as $g=mn\sigma(\mf{q})k$, where $m\in M_P(\mb{A})$, $n\in N_P(\mb{A})$ and $k\in\wh{K}$, then we define a function $\wtl{f}_{M_P,\nu}$ on $\wtl{G}(\mb{A}\lg t\rg)$ associated to $f_{M_P}$ and $\nu$ by
\be\label{nonmaxf}
\wtl{f}_{M_P,\nu}(g)=\nu(m)f_{M_P}(m).
\ee
Then we form the Eisensten series
\be\label{nonmax}
E(\nu,f_{M_P},g)=\sum_{\gamma\in P(F)\bs \wh{G}(F\lg t\rg)}\wtl{f}_{M_P,\nu}(\gamma g).
\ee
Similarly one may verify invariance properties of $\wtl{f}_{M_P,\nu}$ and $E(\nu,f_{M_P},g)$, and check that they are well-defined.

For later use, we shall specialize to the case $P$ is maximal. We follow the treatment in \cite{sha1}. Assume $P=P_\t$, and
$\a_P$ is the corresponding simple root. Let $A_P$ denote
the (split) torus in the center of $M_P$. For any group $H$ defined over $F$, let $X(H)_F$ be the group of $F$-rational characters of $H$. Set
\[
\mf{a}_P=\rm{Hom}(X(M_P)_F,\mb{R})
\]
 the real Lie algebra of $A_P$. Then
\[
\mf{a}_P^\ast=X(M_P)_F\ot_\mb{Z}\mb{R}=X(A_P)_F\ot_\mb{Z}\mb{R}.
\]
Let $H_P: M_P(\mb{A})\to \mf{a}_P$ be the homomorphism defined in \cite{sha1}. Let $\rho_{M_P}$ be the half sum of the roots in $\Phi_{\t+}$, and $\rho_P=\wtl{\rho}-\rho_{M_P}$. Then $\wtl{\a}_P=\lg\rho_P,\a_P^\vee\rg^{-1}\rho_P$ belongs to $X(\wh{T})_F$, thus by restriction to $A_P$ can be
viewed as an element in $\mf{a}_P^\ast$. We shall now identify $s\in\mb{C}$ with $s\wtl{\a}_P\in\mf{a}_{P,\mb{C}}^\ast$, and with
\be\label{nup}s\nu_P=\exp\lg s\wtl{\a}_P, H_P(\cdot)\rg\ee
which is an unramified quasi-character of $M_P(\mb{A})$. Then we set
$\wtl{f}_{M_P,s}=\wtl{f}_{M_P,s\nu_P}$ and $E(s,f_{M_P},g)=E(s\nu_P,f_{M_P},g)$, which are defined by (\ref{nonmaxf}) and (\ref{nonmax}). Finally we remark that
this definition is compatible with previous definition of $E(s,f,g)$, but different from the usual one which uses the quasi-character
\[
\exp\lg s\wtl{\a}_P+\rho_P, H_P(\cdot)\rg.
\]
For the case $P=P_\D$ the latter corresponds to shifting $s$ by $h^\vee$.

\subsection{Absolute convergence of constant terms}

The constant term of $E(s,f,g)$ along the unipotent radical $\wh{U}$ of $\wh{B}$ is defined to be
the following integral
\be \label{const}
E_{\wh{B}}(s,f,g)=\int_{\wh{U}(F)\bs\wh{U}(\mb{A})}E(s,f,ug)du.
\ee
We have to specify the topology of $\wh{U}(F)\bs\wh{U}(\mb{A})$ and the measure $du.$ By Lemma \ref{uadecom} it suffices to define topologies and measures on $\mb{A}\lg t \rg_+/F\lg t \rg_+$ and $(1+t\mb{A}\lg t \rg_+)/(1+tF\lg t \rg_+).$

\begin{lemma}\label{quo}
The natural map $\vp: \mb{A}\lg t \rg_+/F\lg t \rg_+\lra\mb{A}[[t]]/F[[t]]\is\prod\limits^\i_{i=0}(\mb{A}/F)_i$ is an isomorphism of abelian groups.
\end{lemma}
\textit{Proof.} $\vp$ is clearly injective. Let $S_\i$ be the set of all infinite places of $F$, and let
\[
\mb{A}_{S_\i}=\prod_{v\in S_\i} F_v \times\prod_{v\not\in S_\i}{\cal O}_v,
\]
then $F+\mb{A}_{S_\i}=\mb{A}.$ Hence for any $u=\sum\limits^\i_{i=0}u_it^i \mod F[[t]]\in \mb{A}[[t]]/F[[t]],$ we may assume that
$u_i\in\mb{A}_{S_\i}$ for each $i$. Then $u\in \mb{A}\lg t \rg_+$ and therefore $\vp$ is surjective.
\hfill$\Box$

\begin{lemma}\label{quo2}
The natural map
\[
\tau: (1+t\mb{A}\lg t \rg_+)/(1+tF\lg t \rg_+)\lra (1+t\mb{A}[[t]])/(1+tF[[t]])\stackrel{\log}{\lra}t\mb{A}[[t]]/tF[[t]]
 \]
 is an isomorphism of abelian groups.
\end{lemma}
\textit{Proof.} Again it is clear that $\tau$ is injective. To prove $\tau$ is surjective, we have to show that for any $x\in 1+t\mb{A}[[t]]$
there exists $y\in 1+tF[[t]]$ such that $z=xy\in 1+t\mb{A}\lg t \rg_+.$ Write $x=1+\sum^\i_{i=1}x_i t^i,$ $y=1+\sum^\i_{i=1}y_it^i$, then
$z=1+\sum^\i_{i=1}z_it^i$ with
\[
z_i=x_i+x_{i-1}y_1+\cs+ x_1y_{i-1}+y_i.
\]
Apply $F+\mb{A}_{S_\i}=\mb{A}$ repeatedly, we can find a sequence of $y_i\in F$ such that $z_i\in \mb{A}_{S_\i}$ for each $i$. This finishes the proof. \hfill$\Box$

\

Since $\mb{A}/F$ is compact, Lemma \ref{quo} and Lemma \ref{quo2} imply that $\mb{A}\lg t \rg_+/F\lg t \rg_+$ and $(1+t\mb{A}\lg t \rg_+)/(1+tF\lg t \rg_+)$ are compact and inherit the product measure from that of $\mb{A}/F$, which will be defined as follows.

We first specify the self-dual Haar measure on the local field $F_v$ with respect to a non-trivial additive character $\psi_{F_v}$ of $F_v$.
If $F_v=\mb{R}$, we take
\be\label{psir}
\psi_{\mb{R}}(x)=e^{2\pi ix}
\ee
 and $dx$ is the usual Lebesgue measure on $\mb{R}$; if $F_v=\mb{C}$, we take
 \be\label{psic}
 \psi_{\mb{C}}(z)=e^{2\pi i\rm{tr}z}
=e^{4\pi i \mf{Re} z}
 \ee
 and $dzd\overline{z}=2dxdy$ is twice the usual Lebesgue measure on $\mb{C}$; if $F_v$ is a finite extension of $\mb{Q}_p$, then we first take the character
 $\psi_p$ of $\mb{Q}_p$ given by
 \be\label{psip}
 \psi_p(x)=e^{-2\pi i(\rm{fractional part of }x)},
  \ee
  and define $\psi_{F_v}$ by
  \be\label{psiv}
  \psi_{F_v}(x)=\psi_p(\rm{tr}_{F_v/\mb{Q}_p}x).
   \ee
   The self-dual Haar measure on $\mb{Q}_p$ satisfies $\rm{vol}(\mb{Z}_p)=1.$ Let ${\cal O}_{v}$ be the ring of integers of $F_v$, and $\de_{F_v}^{-1}$ be the inverse different
\[
\de^{-1}_{F_v}=\{x\in F_v|\psi(x{\cal O}_{F_v})=1\}=\{x\in F_v|\rm{tr}(x{\cal O}_v)\in\mb{Z}_p\}.
\]
Then the self-dual Haar measure on $F_v$ satisfies $\rm{vol}({\cal O}_v)=N(\de_{F_v})^{-\f{1}{2}}.$ If $\varpi_v\in{\cal O}_v$ is a uniformizer, then $\de_{F_v}=\varpi_v^e{\cal O}_v$ for some non-negative integer $e$, and $N(\de_{F_v})=q_v^e$, where $q_v$ is the cardinality of the residue field ${\cal O}_v/\mf{p}_v.$

Let us compute $\rm{vol}(\mb{A}/F)$ under above self-dual measures. Recall from (\ref{D}) that $\mb{A}/F$ has a fundamental domain
${\cal D}={\cal D}_\i\times\prod\limits_{i<\i}{\cal O}_v.$ Therefore
\[
\rm{vol}(\mb{A}/F)=\rm{vol}((F\ot_{\mb{Q}}\mb{R})/{\cal O}_F)\times\prod_{v<\i}\rm{vol}({\cal O}_v).
\]
$F$ is unramified at almost all places $v$, hence the right-hand-side is a finite product.
It is known that $\rm{vol}((F\ot_\mb{Q}\mb{R})/{\cal O}_F)=|\D_F|^{\f{1}{2}},$ where $\D_F$ is the discriminant of $F$.
On the other hand
\[
\prod_{v<\i}\rm{vol}({\cal O}_v)=\prod_{v<\i}N(\de_{F_v})^{-\f{1}{2}}=\prod_{v<\i}|\D_{F_v}|^{-\f{1}{2}}=|\D_F|^{-\f{1}{2}},
\]
where $\D_{F_v}$ is the relative discriminant $\D_{F_v/\mb{Q}_p}$ with $v|p.$ We conclude that $\rm{vol}(\mb{A}/F)=1.$
It follows that the quotient spaces $U_\a(F)\bs U_\a(\mb{A})$ ($\a\in\wtl{\Phi}_+$) and $\wh{U}(F)\bs\wh{U}(\mb{A})$ are compact with volume equal to 1.

The main result of this section is the following theorem.

\begin{thm}\label{conv}
$($i$)$ Suppose that $g\in\wh{G}(\mb{A}\lg t \rg)\rtimes \s(\mf{q})$ with $\mf{q}\in\mb{I}$ and $|\mf{q}|>1,$ $s\in H=\{z\in\mb{C}|\mathrm{\mf{Re} }z>h+h^\vee\},$ where $h$ $($resp. $h^\vee)$ is the Coxeter $($resp. dual Coxeter$)$ number. Then
$E(s,f,ug)$, as a function on $\wh{U}(F)\bs\wh{U}(\mb{A})$,
converges absolutely outside a subset of measure zero and is
measurable,
\\
$($ii$)$ For any $\v, \eta>0,$ let $H_\v=\{z\in\mb{C}|\mathrm{\mf{Re}}z\geq h+h^\vee+\v\},$ $\s_\eta=\{\s(\mf{q})|\mf{q}\in\mb{I}, |\mf{q}|\geq 1+\eta\}.$ The integral $($\ref{const}$)$ defining
$E_{\wh{B}}(s,f,g)$ converges absolutely and uniformly for $s\in H_\v,$ $g\in \wh{U}(\mb{A})\Om \s_\eta\wh{K},$ where $\Om$ is a compact
subset of $T(\mb{A}).$
\end{thm}

Define the height function $h_s$, $s\in\mb{C}$, on $\wtl{G}(\mb{A}\lg t \rg)$ by
\be\label{htf}
h_s(cp\s(\mf{q})k)=|c|^s
\ee
if $c, \mf{q}\in\mb{I},$ $p\in G(\mb{A}\lg t \rg_+),$ $k\in\wh{K}.$ Then $h_s$ has the same invariance properties as those of $\wtl{f}_s.$ Also note that
the restriction of $h_s$ on $\wtl{T}$ can be expressed as $h_s(a)=a^{sL}.$
Let us define
\be\label{eish}
E(s,h,g)=\sum_{\g\in \wh{G}(F\lg t \rg_+)\bs \wh{G}(F\lg t \rg)}h_s(\g g),
\ee
and
\be\label{consth}
E_{\wh{B}}(s,h,g)=\int_{\wh{U}(F)\bs\wh{U}(\mb{A})}E(s,h,ug)du.
\ee

\begin{lemma}\label{convh}
Theorem \ref{conv} is true for $E(s,h,g).$ Moreover, for $s$, $\mf{q}$
satisfying the conditions of the theorem and $a\in \wh{T}$ one has
\be\label{zeta} E_{\wh{B}}(s,h,a\s(\mf{q}))=\sum_{w\in
W(\D,\emptyset)}(a\s(\mf{q}))^{\wtl{\rho}+w^{-1}(sL-\wtl{\rho})}c_w(s),
\ee where \be\label{cfun} c_w(s)=\prod_{\b\in \wtl{\Phi}_+\cap
w\wtl{\Phi}_-}\f{\La_F(\lg sL-\wtl{\rho},\b^\vee\rg)}{\La_F(\lg
sL-\wtl{\rho},\b^\vee\rg+1)}, \ee with $\La_F$ the normalized Dedekind zeta
function of $F$ defined below.
\end{lemma}

Let $r_1$ (resp. $r_2$) be the number of real  (resp. complex) places of $F$, and let
\[
\G_{\mb{R}}(s)=\pi^{-s/2}\G(s/2),\q
\G_{\mb{C}}(s)=2(2\pi)^{-s}\G(s)
\]
where $\G(s)$ is the Gamma function. Then $\La_F$ is given by
\[
\La_F(s)=|\D_F|^{\f{s}{2}}\G_{\mb{R}}(s)^{r_1}\G_{\mb{C}}(s)^{r_2}\zeta_F(s),
\]
where $\zeta_F$ is the Dedekind zeta function of $F$ and has an Euler product over all prime ideals ${\cal P}$ of ${\cal O}_F$
\[
\zeta_F(s)=\prod_{{\cal P}\ss{\cal O}_F}\f{1}{1-N_{F/\mb{Q}}({\cal P})^{-s}},\q \mf{Re}s>1.
\]
E. Hecke proved that $\zeta_F$ has a meromorphic continuation to the complex plane with only a simple pole
at $s=1.$ Moreover $\La_F$ satisfies the functional equation
\[
\La_F(s)=\La_F(1-s).
\]

\textit{Proof} of Lemma \ref{convh}: From the Bruhat decomposition (Theorem \ref{bruhat}), it follows
that
\[
\wh{G}(F\lg t \rg)=\bigcup_{w\in W(\D,\emptyset)}\wh{G}(F\lg t \rg_+)wU_w(F).
\]
Note that $U_w=\prod\limits_{\a>0,w\a<0}U_\a$ is finite-dimensional (of dimension $l(w)$), hence lies in $G(F\lg t \rg_+).$
Then for $u\in \wh{U}(\mb{A})$,
\[
E(s,h,ug)=\sum_{w\in W(\D,\emptyset)}\sum_{\g\in U_w(F)}h_s(w\g ug)=:\sum_{w\in W(\D,\emptyset)}H_w(s,ug).
\]
We first prove that each inner sum $H_w(s,ug)$ is a measurable function on $\wh{U}(F)\bs\wh{U}(\mb{A}).$
Let us introduce
\be \label{uwa}
U_w'(\mb{A})=G(\mb{A}\lg t \rg_+)\cap \prod_{\a>0, w\a>0} U_\a(\mb{A}),
\ee
and
\be \label{uwk}
U_w'(F)=G(F\lg t \rg_+)\cap \prod_{\a>0, w\a>0}U_\a(F).
\ee
Then $\wh{U}=U_w'U_w$ and $wU_w'w^{-1}\ss\wh{U}$, therefore
\[
H_w(s,ug)=\sum_{\g\in U_w'(F)\bs \wh{U}(F)}h_s(w\g ug),
\]
which is left $\wh{U}(F)$-invariant. Since $\rm{dim }U_w=l(w)<\i$, applying Lemma \ref{ma} ($c$) and Corollary \ref{r} ($a$)
it is easy to see that there exists $i_w\in\mb{N}$ such that the commutator $[U_w, U_\b]\ss U_w'$ for each $\b=\a+i\de$ with $\a\in\Phi\cup\{0\}$ and $i\geq i_w$.
If we define
\bq
\label{uw''a}U_w''(\mb{A})&=&G(\mb{A}\lg t \rg_+)\cap\prod_{\a\in\Phi\cup\{0\}, i\geq i_w}U_{\a+i\de}(\mb{A}),\\
\label{uw''F}U_w''(F)&=&G(F\lg t \rg_+)\cap\prod_{\a\in\Phi\cup\{0\}, i\geq i_w}U_{\a+i\de}(F),
\eq
then $U_w''$ is of finite codimension in $\wh{U}$ and $H_w(s,ug)$ is left $U_w''(F)\bs U_w''(\mb{A})$-invariant. This proves that $H_w(s,ug)$ is
measurable.

The lemma can be reduced to the case $s\in\mb{R}.$ Indeed, since $|h_s|=h_{\mf{Re}s}$
we have $|H_w(s,ug)|\leq H_w(\mf{Re}s,ug).$ By Fubini's theorem,
\[
E_{\wh{B}}(\mf{Re}s,h,g)=\int_{\wh{U}(F)\bs\wh{U}(\mb{A})}E(\mf{Re}s,h,ug)du=\sum_{w\in W(\D,\emptyset)}\int_{\wh{U}(F)\bs\wh{U}(\mb{A})}H_w(\mf{Re}s,ug)du.
\]
If we can show that $E_{\wh{B}}(\mf{Re}s,h,g)$ is finite, then $E(\mf{Re}s,h,g)$ converges almost everywhere, hence is measurable. It follows that
$E(s,h,ug)$ converges absolutely almost everywhere and is measurable. Using Lebesgue dominated convergence theorem we get
 \[
  E_{\wh{B}}(s,h,g)=\int_{\wh{U}(F)\bs\wh{U}(\mb{A})}E(s,h,ug)du=\sum_{w\in W(\D,\emptyset)}\int_{\wh{U}(F)\bs\wh{U}(\mb{A})}H_w(s,ug)du.
 \]
 So we may assume that $s\in\mb{R}$. Let us evaluate $E_{\wh{B}}(s,h,g)$, prove its finiteness and uniform convergence. The computation for general $s$ is the same.

 Since $E_{\wh{B}}$ is left $\wh{U}$-invariant
and right $\wh{K}$-invariant, by Iwasawa decomposition we may
also assume that $g=a\s(\mf{q})$ with $a\in T$ and $|\mf{q}|>1.$ Let us prove
(\ref{zeta}), and show that the summation is finite. By previous discussion we have
\bq
E_{\wh{B}}(s,h,g)&=&\sum_{w\in W(\D,\emptyset)} \int_{\wh{U}(F)\bs \wh{U}(\mb{A})}H_w(s,ug)du\label{int} \\
&=&\sum_{w\in W(\D,\emptyset)} \int_{\wh{U}(F)\bs \wh{U}(\mb{A})}\sum_{\g\in U_w'(F)\bs\wh{U}(F)}h_s(w\g u g)du\n\\
&=& \sum_{w\in W(\D,\emptyset)} \int_{U_{w}'(F)\bs U_{w}'(\mb{A})}\int_{U_w(\mb{A})}h_s(w u' u g)du du'\n\\\
&=& \sum_{w\in W(\D,\emptyset)} \int_{U_w(\mb{A})} h_s(wug)du,\n
\eq
where the last equality follows from the facts that $\rm{vol}(U_w'(F)\bs U_w'(\mb{A}))=1$ and that $h_s$ is left $G(\mb{A}\lg t \rg_+)$-invariant.

To evaluate (\ref{int}) let us introduce some notations. Let
$w=r_{i_1}\cs r_{i_l}$ be the reduced expression of $w$, where
$l=l(w).$ Let \[\wtl{\Phi}_w=\wtl{\Phi}_+\cap
w\wtl{\Phi}_-=\{\b_1,\ld,\b_l\},\] where $\b_j= r_{i_1}\cs
r_{i_{j-1}}\a_{i_j}.$ Then
\[
\wtl{\Phi}_{w^{-1}}=\wtl{\Phi}_+\cap w^{-1}\wtl{\Phi}_-=\{\g_1,\ld,\g_l\},
\]
where
$\g_j=-w^{-1}\b_j=r_{i_l}\cs r_{i_{j+1}}\a_{i_j}.$ Note that
\[
\b_1+\cs+\b_l=\wtl{\rho}-w\wtl{\rho},\q \g_1+\cs+\g_l=\wtl{\rho}-w^{-1}\wtl{\rho}.
\]

Recall that we have assumed $s\in\mb{R}$, and $g=a\s(\mf{q}),$ $a\in T.$
We have \bq
&&\int_{U_w(\mb{A})}h_s(wua\s(\mf{q}))du \label{int2} \\
&= & \int_{U_w(\mb{A})}h_s(wa\s(\mf{q})\rm{Ad}(a\s(\mf{q}))^{-1}(u))du\n\\
&=& \int_{\mb{A}^l}h_s(wa\s(\mf{q})\wtl{x}_{\g_1}((a\s(\mf{q}))^{-\g_1}u_1)\cs \wtl{x}_{\g_l}((a\s(\mf{q}))^{-\g_l}u_l))du_1\cs du_l\n\\
&=& (a\s(\mf{q}))^{\g_1+\cs+\g_l+w^{-1}sL}\int_{\mb{A}^l}h_s (w\wtl{x}_{\g_1}(u_1)\cs \wtl{x}_{\g_l}(u_l))du_1\cs du_l\n\\
&=& (a\s(\mf{q}))^{\wtl{\rho}-w^{-1}\wtl{\rho}+w^{-1}sL}\int_{\mb{A}^l}h_s (\wtl{x}_{-\b_1}(u_1)\cs \wtl{x}_{-\b_l}(u_l))du_1\cs du_l.\n
\eq
By Iwasawa decomposition we have
\[
\wtl{x}_{-\b_l}(u_l)=n(u_l)a(u_l)k(u_l),
\]
where $a(u_l)\in T_{\b_l}=\{\wtl{h}_{\b_l}(u)|u\in\mb{I}\},$ $n(u_l)\in U_{\b_l},$ $k(u_l)\in K_{\b_l}.$
Let $w'=r_{i_1}\cs r_{i_{l-1}},$  then $\{\b_1,\ld,\b_{l-1}\}=\wtl{\Phi}_+\cap w'\wtl{\Phi}_-.$ Consider the decomposition
\be\label{dec}
\wh{U}=U_{w'}'U_{w'}
\ee
where $U_{w'}'$ is given by (\ref{uwa}), (\ref{uwk}) with $w$ replaced by $w'.$ Then
\[
U_{-\b_1}\cs U_{-\b_{l-1}}=w'U_{w'}w'^{-1}.
\]
Use (\ref{dec}) we can define the projection
\[
\pi: w'\wh{U}w'^{-1}\lra w'U_{w'}w'^{-1}.
\]
Since $U_{-\b_1},\ld, U_{-\b_{l-1}},$ $U_{\b_l}\ss w'\wh{U}w'^{-1},$ we have the map
\[
\pi\circ\rm{Ad}(n(u_l)): w'U_{w'}w'^{-1}\lra w'U_{w'}w'^{-1}
\]
which is unimodular. From this fact, together with invariance properties of $h_s$ and Corollary \ref{r}
($f$) we get
\bq
&&\int_{\mb{A}^l}h_s(\wtl{x}_{-\b_1}(u_1)\cs \wtl{x}_{-\b_l}(u_l))du_1\cs du_l\label{int3}\\
&=&\int_{\mb{A}^l}h_s\l(\wtl{x}_{-\b_1}(u_1)\cs \wtl{x}_{-\b_{l-1}}(u_{l-1})a(u_l)\r)du_1\cs du_l\n\\
&=&\int_{\mb{A}^l} h_s\l(a(u_l)\wtl{x}_{-\b_1}(a(u_l)^{\b_1}u_1)\cs \wtl{x}_{-\b_{l-1}}(a(u_l)^{\b_{l-1}}u_{l-1})\r)du_1\cs du_l\n\\
&=&\int_{\mb{A}}a(u_l)^{sL-\b_1-\cs-\b_{l-1}}du_l \int_{\mb{A}^{l-1}}h_s(\wtl{x}_{-\b_1}(u_1)\cs \wtl{x}_{-\b_{l-1}}(u_{l-1}))du_1\cs du_{l-1}\n\\
&=&\int_{\mb{A}}a(u_l)^{sL-\wtl{\rho}+w'\wtl{\rho}}du_l \int_{\mb{A}^{l-1}}h_s(\wtl{x}_{-\b_1}(u_1)\cs \wtl{x}_{-\b_{l-1}}(u_{l-1}))du_1\cs du_{l-1}.\n
\eq
From the Gindikin-Karpelevich formula \cite{GGP, Lan}, the first integral in (\ref{int3}) equals
\[
\f{\La_F( z_l)}{\La_F(z_l+1)},
\]
where
\[
z_l=\lg sL-\wtl{\rho}+w'\wtl{\rho},\b_l^\vee\rg-1=\lg sL-\wtl{\rho},\b_l^\vee\rg.
\]
Note that $\lg w'\wtl{\rho},\b_l^\vee\rg=\lg \wtl{\rho}, w'^{-1}\b_l\rg= \lg \wtl{\rho},\a_{i_l}^\vee\rg=1.$ By induction on
$l$ it is clear that (\ref{int3}) equals $c_w(s).$

Now we prove that the right-hand-side of (\ref{zeta}) is finite. We
first prove that $s>h+h^\vee$ implies that $\lg
sL-\wtl{\rho},\b^\vee\rg>1$ for each $\b\in \Phi_w= \wtl{\Phi}_+\cap
w\wtl{\Phi}_-.$ In fact, since $w\ss W(\D,\emptyset),$
$w^{-1}\D\ss\wtl{\Phi}_+.$ It follows that $\b=i\de+\a$ with $i>0,$
$\a\in\Phi.$ By Lemma \ref{coroot}, $\b^\vee=jc+\a^\vee$ with $j>0.$
Then
\[
\lg sL-\wtl{\rho},\b^\vee\rg =j(s-h^\vee)-\lg\rho,\a^\vee\rg\geq
s-h^\vee-h+1>1.
\]
By standard results on zeta functions, for every $\v>0,$ there
exists a constant $c_\v>0$ such that whenever $\mf{Re}z\geq 1+\v$ we have
\[
\l|\f{\La_F(z)}{\La_F(z+1)}\r|<c_\v.
\]
It follows that $c_w(s)\leq c_\v^{l(w)}$ for $s>h+h^\vee+\v.$

Next consider $(a\s(\mf{q}))^{\wtl{\rho}-w^{-1}\wtl{\rho}+w^{-1}sL}.$ Write $
w^{-1}=T_\la w_0$ where $\la\in Q^\vee,$ $w_0\in W.$ By Lemma \ref{weyl2},
\[
w^{-1}L= L+\la-\f{1}{2}\left( \la, \la\right) \de,
\]
\[
w^{-1}\wtl{\rho}=w_0\rho-\lg w_0\rho, \la\rg\de+h^\vee w^{-1}L.
\]
Let $\|\la\|=\left( \la,\la\right)^{\f{1}{2}}.$ If we write $\la=\sum_{i=1}^n l_i\a_i^\vee$, then
there exists a constant $c_1>0$ which does not depend on $\la$ such that $\sum^{n}_{i=1}|l_i|\leq c_1\|\la\|.$ Then
\[
|\lg w_0\rho,\la\rg|=\l|\sum_{i=1}^n l_i\lg\rho,w_0^{-1}\a_i^\vee\rg\r|\leq c_1h\|\la\|.
\]
Combine above equations we obtain
\be\label{esti}
\l\{\ba{l}\s(\mf{q})^{w^{-1}L}=|\mf{q}|^{-\f{1}{2}\|\la\|^2},\\ \s(\mf{q})^{-w^{-1}\wtl{\rho}}=|\mf{q}|^{\lg w_0\rho,\la\rg+\f{h^\vee}{2}\|\la\|^2}\leq|\mf{q}|^{c_1h\|\la\|+\f{h^\vee}{2}\|\la\|^2}.\ea\r.
\ee Let $c_a=\max\limits_{\a\in\Phi}|a^\a|,$ then \be \label{esta}
a^{\wtl{\rho}-w^{-1}\wtl{\rho}+w^{-1}sL}\leq c_a^{l(w)}a^{s\la} \leq
c_a^{l(w)+c_1s\| \la\|}. \ee But we have
\begin{equation}\label{estl}   l(w)\leq  l(T_\la)+l(w_0)\leq  \sum^n_{i=1}|l_i| l(T_{\a_i^\vee})+|\Phi_+|\leq  c_2\|\la\|+|\Phi_+|,
\end{equation}
where $c_2=c_1\max\limits_{1\leq i\leq n}l(T_{\a_i^\vee})$. In summary we
obtain \bq
E_{\wh{B}}(s,h,a\s(\mf{q}))&=&\sum_{w\in W(\D,\emptyset)}(a\s(\mf{q}))^{\wtl{\rho}+w^{-1}(sL-\wtl{\rho})}c_w(s)\n\\
\label{est}&\leq & |W|(c_\v c_a)^{|\Phi_+|}\sum_{\la\in Q^\vee}(c_\v^{c_2}c_a^{c_1s+c_2})^{\|\la\|}|\mf{q}|^{c_1h\|\la\|-\f{s-h^\vee}{2}\|\la\|^2}
\eq
for $s\geq h+h^\vee+\v.$  It is clear that the last series in (\ref{est}) is finite and satisfies the required uniform convergence properties.
\hfill$\Box$

\

The proof of Theorem \ref{conv} is similar to that of Lemma
\ref{convh}, and we only need the following two observations: (1) $f$ is
bounded by cuspidality, (2) Consider
\[
E(s,f,ug)=\sum_{w\in W(\D,\emptyset)} F_w(s,ug)=\sum_{w\in W(\D,\emptyset)}\sum_{\g\in U_w(F)}\wtl{f}_s(w\g ug).
\]
We can modify the definition of $U_w''$ in (\ref{uw''a}) and (\ref{uw''F}) by taking $i_w$ large enough such that
$[U_w,U_\b]\ss w^{-1}N_\D w\ss U_w'$ whenever $\b=\a+i\de$ with $\a\in\Phi\cup\{0\}$ and $i\geq i_w.$
Here by our convention for any $\t\ss\wtl{\D}$ we let
\be \label{rad1}
N_\t(\mb{A})=G(\mb{A}\lg t \rg_+)\cap \prod_{\a\in\wtl{\Phi}_+-\Phi_\t}U_\a(\mb{A}),
\ee
\be \label{rad2}
N_\t(F)=G(F\lg t \rg_+)\cap \prod_{\a\in\wtl{\Phi}_+-\Phi_\t}U_\a(F).
\ee
Then $U_w''$ is again of finite codimension in $\wh{U}$, and $F_w(s,ug)$ is left $U_w''(F)\bs U_w''(\mb{A})$-invariant.

\begin{cor}
If $\mathrm{\mf{Re}}s>h+h^\vee$ and $|\mf{q}|>1,$ then $E(s,f,ua\s(\mf{q}))$ and
$E(s,h,ua\s(\mf{q}))$ defined by (\ref{eis}) and (\ref{eish}), as
functions on $\wh{U}(F)\bs \wh{U}(\mb{A})\times \wh{T}$, are measurable and converge
absolutely outside a subset of measure zero.
\end{cor}

Let us state the Gindikin-Karpelevich formula in the case of $SL_2(F_v)$ where $F_v$ is the completion of $F$ at place $v$, whose proof is well-known. Let $\chi_v$ be an unramified character of $F_v^\times$. Let $f_{s,\chi_v}\in\rm{Ind}_B^{SL_2}(|\c|^s\ot\chi_v)$, $s\in\mb{C},$ $\mf{Re}s>0$ be the unique spherical function satisfying
\be \label{fchi}
f_{s,\chi_v}\l(\bpm a & x \\ 0 & a^{-1} \epm k\r)=\chi_v(a)|a|_v^{s+1},
\ee
where $k\in K$, the standard maximal compact subgroup of $SL_2(F_v).$ Let $w=\bpm 0 & 1\\ -1 & 0\epm.$ Define
\be\label{cx}
c(s,\chi_v)=\int_{F_v}f_{s,\chi_v}\l(w^{-1}\bpm 1 &x \\ 0 &1 \epm \r)dx.
\ee

\begin{prop} \label{gk} Suppose $F_v$ is $p$-adic. Let $\vpi_v$ be a uniformizer of $\mf{p}_v\ss{\cal O}_v$, $q_v$ be the cardinality of the residue field ${\cal O}_v/\mf{p}_v.$ Then
  \[
  c(s,\chi_v)=\mathrm{vol}({\cal O}_v)\f{1-\chi_v(\vpi_v)q_v^{-s-1}}{1-\chi_v(\vpi_v)q_v^{-s}} =\f{|\D_v|^{\f{s}{2}}L(s,\chi_v)}{|\D_v|^{\f{s+1}{2}}L(s+1,\chi_v)}.
  \]
  \end{prop}

Recall that $\D_v$ is the relative discriminant $\D_{F_v/\mb{Q}_p}$ satisfying
$\mathrm{vol}({\cal O}_v)=|\D_v|^{-\f{1}{2}}.$ In particular if $\chi_v$ is trivial, then $c_v(s,\chi_v)$ contributes the local factor of $\df{\La_F(s)}{\La_F(s+1)}$ at $v$.
In the case that $F_v=\mb{R}$ or $\mb{C}$, unramified characters of $F_v$ are of the form $|\c|_v^{s_0}$ for $s_0$ purely imaginary. Then one has the following

\begin{prop} \label{gkreal}
Suppose $F_v=\mb{R}$ or $\mb{C}.$ Then
\[
c(s, |\c|_v^{s_0})=\f{\G_{F_v}(s+s_0)}{\G_{F_v}(s+s_0+1)}.
\]
\end{prop}

\subsection{Constant terms and Fourier coefficients}

In this section we shall compute the constant terms of Eisenstein
series along unipotent radicals of parabolic subgroups of $\wh{G}.$
In the classical theory \cite{Lan}, if $P=MN$ is a parabolic
subgroup of $G$, $f$ is a cusp form on $M(\mb{A})$, then the
$L$-functions associated to $f$ and the representations of $^LM$ on
$^L\mf{n}$ will appear in certain constant terms of the Eisenstein
series\ induced by $f$. Unfortunately for loop groups the constant
terms are all trivial. However we will obtain certain non-trivial higher Fourier coefficients
of the Eisenstein series.

Let $P=P_{\theta}$ and $Q=P_{\theta'}$ be two maximal parabolics of $\wh{G}$ with Levi decompositions
$P=M_PN_P$ and $Q=M_QN_Q$. Let $\a_P$ and $\a_Q$ be the corresponding simple roots. Let $f_{M_P}$ be an unramified cusp form on $M_P(\mb{A})$ and $E(s,f_{M_P},g)$ be the Eisenstein series
defined at the end of section 4.1.
The constant term of $E(s,f_{M_P},g)$ along $N_Q$
is given by
\be \label{constp}
E_Q(s,f_{M_P},g)=\int_{N_Q(F)\bs N_Q(\mb{A})}E(s,f_{M_P},ng)dn.
\ee
Using similar method as
in the proof of Theorem \ref{conv}, one can show that the integral (\ref{constp}) converges absolutely for $|\mf{q}|>1$ and $\mf{Re}s\gg0$. One
can easily find the precise range of convergence for $\mf{Re}(s)$ and we will not give it here.

\begin{thm}\label{triv}
$E_Q(s,f_{M_P},g)=0$ unless $P=Q$, in which case
\[
E_P(s,f_{M_P},g)=(s\nu_P)(m)f_{M_P}(m),
\]
where $g=mn\s(\mf{q})k$, $m\in M_P(\mb{A})$, $n\in N_P(\mb{A})$ and $k\in\wh{K}$.
\end{thm}
\textit{Proof.} We follow the arguments of Langlands \cite{Lan}. We write
\[
E_Q(s,f_{M_P},g)=\sum_{\gamma\in P(F)\bs \wh{G}(F\lg t\rg)/N_Q(F)}\int_{\gamma^{-1}P(F)\gamma\cap N_Q(F)\bs N_Q(\mb{A})}\wtl{f}_{M_P,s}(\gamma n g)dn.
\]
From the Bruhat decomposition we may assume that each $\gamma$ is of the form $\gamma=w\gamma'$ with
$w\in W(\theta_1,\theta_2)$ and $\gamma'\in M_Q(F)$. Then up to a scalar depending on $\gamma'$, a typical integral equals
\[
\int_{w^{-1}P(F)w\cap N_Q(F)\bs N_Q(\mb{A})}\wtl{f}_{M_P,s}(w n \gamma'g)dn.
\]
Since $w\in W(\t_1,\t_2)$ we have
\[
U_w\h w^{-1}P w\cap N_Q\bs N_Q.
\]
Let $w^{-1}Pw\cap N_Q= N_1 N_2,$ where $N_1=w^{-1}M_P w\cap N_Q,$ $N_2=w^{-1}N_P w\cap N_Q$, then above integral equals
\[
\int_{U_w(\mb{A})}\int_{N_2(F)\bs N_2(\mb{A})}\int_{N_1(F)\bs N_1(\mb{A})}\wtl{f}_{M_P,s}(wn_1n_2 u\gamma'g) dn_1 dn_2 du.
\]
Since $f_{M_P}$ is a cusp form on $M_P$, the most inner integral vanishes unless $N_1=1$. Then $w^{-1}\theta\subset\Phi_{\theta'+}$, which forces that $w^{-1}\theta=\theta'$.
By the following lemma, which is essentially Shahidi's Lemma \cite[Lemma 4.1]{Gar6}, we have $w=1$ and $P=Q$.
In this case we can take $\gamma'=1$ and thus $E_P(s,f_{M_P},g)=\wtl{f}_{M_P,s}(g)=(s\nu_P)(m)f_{M_P}(m).$
\hfill$\Box$

\begin{lemma}\label{ext}
In above settings if there exists $w\in\wtl{W}$ such that $w\theta=\theta'$, then $w=1$ and $P=Q$.
\end{lemma}
\textit{Proof.} For convenience let us enumerate $\a_P=\a_i$, $\a_Q=\a_j$. It is enough to show that $w\a_i>0$. For the contrary
suppose $w\a_i<0.$ Let $w_0^\t$ be the longest element in $W_{\t}$, then $ww_0^\t\t=-\t'$. We may write
\[
w_0^\t\a_i=\a_i+\sum_{k\neq i}n_k\a_k.
\]
Then $ww_0^\t\a_i$ has an expression
\[
ww_0^\t\a_i=w\a_i+\sum_{k\neq j}n_k'\a_k
\]
since $w\t=\t'$. Let us write $w\a_i=\sum^{n}_{k=0}b_k\a_k$. If $b_j=0$ then $w\a_i\in \Phi_{\t'}$, which is impossible. Therefore $b_j<0$, which further implies $ww_0^\t\a_i<0$. Thus we obtain $ww_0^\t\wtl{\D}<0$, a contradiction.
\hfill$\Box$

\

In general, for a connected reductive algebraic group $G$ which is
split over $F$ (e.g. $G=GL_n$), the theory of generalized Tits
system \cite{iwa, IM} implies that $G(F((t)))=\wh{B}_0
\wtl{W}'\wh{B}_0$ where $\wtl{W}'=W\ltimes X_\ast(T).$ The proof of
Theorem \ref{triv} together with certain variants of Lemma \ref{ext} suggests that the triviality
of constant terms should also hold for $\wh{G}.$

Now let us define and compute the Fourier coefficients of $E(s,f_{M_P},g)$.
Let $\psi$ be a character of $\wh{U}(F)\bs\wh{U}(\mb{A})$, then $\psi=\prod\limits_{\a\in\wtl{\D}}\psi_\a,$
where $\psi_\a$ is a character of $U_\a(F)\bs U_\a(\mb{A}).$ This follows from the fact
\[
\wh{U}/[\wh{U},\wh{U}]\h\prod_{\a\in\wtl{\D}}U_\a.
\]
Define the $\psi$-th Fourier
coefficient of $E(s,f_{M_P},g)$ along $\wh{B}$ by
\[
E_{\wh{B},\psi}(s,f_{M_P},g)=\int_{\wh{U}(F)\bs\wh{U}(\mb{A})}E(s,f_{M_P},g)\o{\psi}(u)du.
\]
Then $E_{\wh{B},\psi}(s,f_{M_P},g)$ is a Whittaker function on $\wh{G}(\mb{A}\lg t \rg)\rtimes\s(\mf{q})$, i.e. a function $W$ satisfying
the relation $W(ug)=\psi(u)W(g),$ $\forall u\in \wh{U}(\mb{A}).$ Let $U_P=M_P\cap\wh{U}$ be the unipotent radical of $M_P\cap\wh{B}\subset M_P$. We say that $\psi|_{U_P}$ is generic if $\psi_\a$ is non-trivial for each $\a\in\t$,
and that $\psi$ is generic if $\psi_{\a}$ is non-trivial for each $\a\in\wtl{\D}$.

\begin{thm}\label{fc}
$(i)$ Assume that $\psi$ is generic. Then $E_{\wh{B},\psi}(s,f_{M_P},g)=0$.

$(ii)$
Assume that $\psi_{\a_P}$ is trivial. Then $E_{\wh{B},\psi}(s,f_{M_P},g)=0$ unless $\psi|_{U_P}$ is generic, in which case
\[
E_{\wh{B},\psi}(s,f_{M_P},g)=(s\nu_P)(m)\int_{U_P(F)\bs U_P(\mb{A})}f_{M_P}(um)\o{\psi(u)}du,
\]
where $g=mn\s(\mf{q})k$, $m\in M_P(\mb{A})$, $n\in N_P(\mb{A})$ and $k\in\wh{K}$.
\end{thm}
\textit{Proof.} $(i)$ Since $\psi$ is $\wh{U}(F)$-invariant, similarly as
before we have \bq E_{\wh{B},\psi}(s,f_{M_P},g)&=&\sum_{w\in W(\t,\emptyset)}
\int_{\wh{U}(F)\bs\wh{U}(\mb{A})}\sum_{\g\in U_w(F)}\wtl{f}_{M_P,s}(w\g
ug)\o{\psi}(u)du\n\\
&=&\sum_{w\in W(\t,\emptyset)}\int_{w^{-1}P(F)w\cap
\wh{U}(F)\bs\wh{U}(\mb{A})}\wtl{f}_{M_P,s}(wug)\o{\psi}(u)du.\n \eq
For each
$w\in W(\t,\emptyset),$ $w^{-1}\t\ss\wtl{\Phi}_+,$ therefore
$w^{-1}U(F)w\ss w^{-1}P(F)w\cap \wh{U}(F).$ If we
write $\wh{U}=N_{2w} N_{1w} U_w$, where $N_{1w}=w^{-1}U_Pw,$ $N_{2w}=w^{-1}N_P w\cap \wh{U}$, then a typical integral
equals
\[
\int_{U_w(\mb{A})}\int_{N_{1w}(F)\bs
N_{1w}(\mb{A})}\int_{N_{2w}(F)\bs
N_{2w}(\mb{A})}\wtl{f}_{M_P,s}(wn_2n_1ug)\o{\psi}(n_2n_1u)dn_2 dn_1 du.
\]
By definition of $\wtl{f}_{M_P,s}$ the function $\wtl{f}_{M_P,s}(wn_2n_1ug)$ is constant on $N_{2w}$. Hence we are reduced
to show that $\psi|_{N_{2w}}$ is non-trivial, provided that $\psi$ is generic. In this case $\psi|_{N_{2w}}$ is trivial if and only
if $U_{\a}\not\subset w^{-1}N_P w$ for any $\a\in\wtl{\D}$. Equivalently we are looking for $w$ satisfying the conditions
$w\widetilde{\Delta}\subset \widetilde{\Phi}_- \cup \Phi_{\t+}$ and $w^{-1}\t>0$. For such $w$ in fact one has
\[
w\wtl{\D}\ss (\wtl{\Phi}_--\Phi_{\t-})\cup\Phi_{\t+}.
\]
From this it is easy to deduce that $w\wtl{\D}\cap\Phi_{\t+}$ forms a system of simple roots for $\Phi_\t$, hence $w\wtl{\D}\cap\Phi_{\t+}=\t$.
By Lemma \ref{ext} we must have $w=1$, which is obviously a contradiction. Therefore we have proved that such $w$ does not exist, hence
$\psi|_{N_{2w}}$ is non-trivial.

$(ii)$ We reverse the order of $N_{1w}$ and $N_{2w}$ to rewrite the integral as
\[
\int_{N_{2w}(F)\bs
N_{2w}(\mb{A})}\int_{N_{1w}(F)\bs
N_{1w}(\mb{A})}\wtl{f}_{M_P,s}(wn_1n_2ug)\o{\psi}(n_1n_2u)dn_1 dn_2 du.
\]
Consider the subset $\t_w=\t\cap w\Phi_{\t+}$ of $\t$ and let $P_w$ be the corresponding parabolic subgroup of $M_P$. If $\t_w\neq\t,$ then
$P_w$ has non-trivial unipotent radical $U_{P_w}\ss U_P.$ From the definition of $\t_w$ it follows that
$w^{-1}U_{P_w}w\ss N_{\t}$ and therefore by our assumption $\o{\psi}$ is trivial on $w^{-1}U_{P_w}w.$ Since $f_{M_P}$ is a cusp form
on $M_P(\mb{A})$, it follows that the most inner integral vanishes.
If $\t_w=\t$, then $w^{-1}\t\subset\Phi_{\t+}$, from which we deduce that $w^{-1}\t=\t$ and hence $w=1$ by Lemma \ref{ext}. It follows that
\begin{eqnarray*}
E_{\wh{B},\psi}(s,f_{M_P},g)&=&\int_{\wh{U}(F)\bs\wh{U}(\mb{A})}\wtl{f}_{M_P,s}(ug)\o{\psi}(u)du\\
&=&(s\nu_P)(m)\int_{U_P(F)\bs U_P(\mb{A})}f_{M_P}(um)\o{\psi}(u)du,
\end{eqnarray*}
which vanishes unless $\psi|_{U_P}$ is generic, again by cuspidality of $f_{M_P}.$
\hfill$\Box$

\

In the rest of this section we consider the Eisenstein series induced from a quasi-character on $\wh{T}$. More precisely,
let $\chi_{\wh{T}}=\bot\limits_v\chi_{\wh{T}_v}: \wh{T}(\mb{A})/\wh{T}(F)\ra\mb{C}^\times$ be a quasi-character such that $\chi_{\wh{T}_v}$ is unramified for each place $v$. Extend $\chi_{\wh{T}}$ to $\wtl{T}=\wh{T}\times \s(\mb{I})$ such that $\chi_{\wtl{T}}|_{\s(\mb{I})}$ is trivial. For each $\a\in\wtl{\Phi}_{re},$ $\chi_{\wh{T}}\a^\vee: \mb{I}/F^\times\ra\mb{C}^\times$ is the Hecke quasi-character such that $\chi_{\wh{T}}\a^\vee(x)=\chi_{\wh{T}}(\wtl{h}_\a(x)).$ In particular $\chi_{\wh{T}_v}\a^\vee$ is unramified for each place $v.$

In general if $\chi: \mb{I}/F^\ti\ra\mb{C}^\times$ is a Hecke quasi-character, we may write $\chi$ as $\mu|\c|^{s_0}$ where $\mu$ is unitary and $s_0\in\mb{C}.$ Define $\mf{Re}\chi=\mf{Re}s_0$, which is called the exponent of $\chi.$ Recall that $L(s,\chi)=\prod\limits_v L(s,\chi_v)$ is the Hecke $L$-function of $\chi.$  We may twist $\chi$ to make it unitary, then $L(s,\chi)$ is holomorphic in $\{s\in\mathbb{C}: \mf{Re}s>1\},$ admits meromorphic continuation to the entire complex plane, and satisfies the functional equation
\[
L(1-s,\chi^\vee)=\varepsilon(s,\chi)L(s,\chi),
\]
where $\chi^\vee=\chi^{-1}|\c|$ is the shifted dual of $\chi.$

We define the Eisenstein series on $\wtl{G}(\mb{A}\lg t \rg)$ induced
from $\chi_{\wh{T}}$ by
\be\label{chi}
E(\chi_{\wh{T}},g)=\sum_{\g\in \wh{B}(F)\bs\wh{G}(F\lg t \rg)}\wtl{\chi}_{\wh{T}}(\g g),
\ee
where $\wtl{\chi}_{\wh{T}}$ is given by
\be\label{chis}
\wtl{\chi}_{\wh{T}}(cb\s(\mf{q})k)=\chi_{\wh{T}}(cb_0),
\ee
where $c,$ $\mf{q}\in\mb{I},$ $b\in \wh{B}_0(\mb{A})$, $b_0\in B(\mb{A})$
is the image of $b$ under the projection $\wh{B}_0(\mb{A})\ra
B(\mb{A})$, and $k\in\wh{K}.$ Then $\wtl{\chi}_{\wh{T}}$ is well-defined,
right $\wh{K}$-invariant and left $\wh{B}(F)$-invariant. We define
the constant term and Fourier coefficients of $E(\chi,g)$ along
$\wh{B}$ by \bq
 E_{\wh{B}}(\chi_{\wh{T}},g)&=&\int_{\wh{U}(F)\bs\wh{U}(\mb{A})}E(\chi_{\wh{T}},ug)du,\n\\
E_{\wh{B},\psi}(\chi_{\wh{T}},g)&=&\int_{\wh{U}(F)\bs\wh{U}(\mb{A})}E(\chi_{\wh{T}},ug)\o{\psi}(u)du.\n
\eq

\begin{thm}\label{cha}
$($i$)$ Suppose that $g\in\wh{G}(\mb{A}\lg t \rg)\rtimes \s(\mf{q})$ with $|\mf{q}|>1,$ and $\mf{Re}(\chi_{\wh{T}}\a_i^\vee)>2,$ $i=0,1,\ld,n.$ Then
$E(\chi_{\wh{T}},ug)$, as a function on $\wh{U}(F)\bs\wh{U}(\mb{A})$,
converges absolutely outside a subset of measure zero and is
measurable.
\\
$($ii$)$ For any $\v>0$, let
\[
{\cal H}_\v=\l\{\chi_{\wh{T}}:\wh{T}(\mb{A})/\wh{T}(F)\ra\mb{C}^\times\l|\chi_{\wh{T}}\rm{ unramified}, \mf{Re}(\chi_{\wh{T}}\a_i^\vee)>2+\v, i=0,1,\ld,n\r.\r\}.
\]
The integral defining
$E_{\wh{B}}(\chi_{\wh{T}},g)$ converges absolutely and uniformly for $\chi_{\wh{T}}\in {\cal H}_\v,$ $g\in \wh{U}(\mb{A})\Om \s_\eta\wh{K},$ where $\eta>0$ and $\Om$ is a compact subset of $T(\mb{A}).$ More explicitly for $a\in\wh{T}$ one has
\be
E_{\wh{B}}(\chi_{\wh{T}},a\s(\mf{q}))=\sum_{w\in\wtl{W}}(a\s(\mf{q}))^{\wtl{\rho}+w^{-1}(\chi_{\wh{T}}-\wtl{\rho})}c_{w}(\chi_{\wh{T}})
\ee where \be\label{xfun} c_{w}(\chi_{\wh{T}})=\prod_{\b\in \wtl{\Phi}_+\cap
w\wtl{\Phi}_-}|\D_F|^{-\f{1}{2}}\f{L(-\lg\wtl{\rho},\b^\vee\rg, \chi_{\wh{T}} \b^\vee)}{L(1-\lg\wtl{\rho},\b^\vee\rg,\chi_{\wh{T}}\b^\vee)}.
\ee
\end{thm}
\textit{Proof.} The proof follows exactly the line for proving Lemma \ref{convh} if we apply analytic properties of Hecke $L$-functions together with Propositions \ref{gk} and \ref{gkreal}. With the analog of Godement's criterion
\be\label{gode}
\mf{Re}(\chi_{\wh{T}}\a_i^\vee)>2,\q i=0,1,\ld,n,
\ee
we have the following two observations, which suffice for the convergence of Eisenstein series:
(1) $\mf{Re}(\chi_{\wh{T}}\b^\vee)-\lg\wtl{\rho},\b^\vee\rg>1$ for any $\b\in\wtl{\Phi}_{re+}$, which is precisely (\ref{gode}) when $\b$ is a simple root.
(2) Consider the factor $\s(\mf{q})^{w^{-1}(\chi_{\wh{T}}-\wtl{\rho})}=(w\c\s(\mf{q}))^{\chi_{\wh{T}}-\wtl{\rho}}$. Write $w=T_\la w_0\in\wtl{W}$ such that $\la\in Q^\vee,$ $w_0\in W$. Then from
\[
wd=T_\la d=d+\la-\f{\|\la\|^2}{2}c,
\]
we see that
\[
|(w\c\s(\mf{q}))^{\chi_{\wh{T}}-\wtl{\rho}}|=|\mf{q}|^{\mf{Re}(\chi_{\wh{T}}\la)-\lg\rho,\la\rg+\l(h^\vee-\mf{Re}(\chi_{\wh{T}} c)\r)\f{\|\la\|^2}{2}}.
\]
The coefficient of the quadratic term is negative. Indeed,
\[
\mf{Re}(\chi_{\wh{T}} c)= \mf{Re}\l(\chi_{\wh{T}} (\a_0^\vee+\wtl{\a}^\vee)\r)>  2+\lg\wtl{\rho},\a_0^\vee+\wtl{\a}^\vee\rg= 2+h^\vee.
\]
Let us remark that due to the last equation, for our second consideration Godement's criterion is much stronger than required.\hfill$\Box$

\subsection{Explicit computations for $\wtl{SL}_2$}

It would be interesting to investigate the Fourier coefficients of our Eisenstein series in full generality. To obtain an explicit formula would be a quite difficult and non-trivial problem, as suggested by the paper of W. Casselman and J. Shalika \cite{CS} where they gave a formula for finite-dimensional groups. The reason their method does not work here is that we do not have a longest element in the affine Weyl group, as opposed to the classical case. However in the case of $SL_2$ we do have explicit computations and everything is known.

Let $\a_0=\de-\a, \a_1=\a$ be the simple roots of $\wtl{\mf{sl}}_2,$ $\psi=\psi_0\psi_1$ be a character on $\wh{U}(F)\bs\wh{U}(\mb{A})$ where $\psi_i$ corresponds to $\a_i,$ $i=0,1.$ Let $f$ be an unramified cusp form on $SL_2(\mb{A}).$

\begin{prop}
Assume that $\psi$ is non-trivial. Then $E_{\wh{B},\psi}(s,f,g)$ vanishes unless $\psi$ equals $\psi_1$ and is non-trivial, in which case it equals the Fourier coefficient
of $f$
\[
E_{\wh{B},\psi}(s,f,g)=|c|^s\int_{F\bs\mb{A}}f\l(\bpm 1 & u \\ 0 & 1\epm p_0\r)\overline{\psi}_1(u)du,
\]
where $g=c p\s(\mf{q})k$, $c\in\mb{I}$, $p\in SL_2(\mb{A}\lg t\rg_+)$, $k\in\wh{K}$.
\end{prop}
\textit{Proof.} The case $\psi=\psi_1$ is known due to Theorem \ref{fc} hence we only need to show the vanishing of $E_{\wh{B},\psi}$ in other cases. As usual let us write
\[
E_{\wh{B},\psi}(s,f,g)=\sum_{w\in\wtl{W}, w^{-1}\a>0}\int_{w^{-1}\wh{G}(F\lg t \rg_+)w\cap
\wh{U}(F)\bs\wh{U}(\mb{A})}\wtl{f}_s(wug)\o{\psi}(u)du.
\]
For $w\in\wtl{W}$ such that $w^{-1}\a>0,$ if $w\neq 1$ then $w^{-1}\a\neq \a_0, \a_1$, which implies that $\psi$ is trivial on the root subgroup $U_{w^{-1}\a}=w^{-1}U_\a w.$ On the other hand the integration of $\wtl{f}_s$ along $U_\a(F)\bs U_\a(\mb{A})$ is zero since $f$ is cuspidal.
Therefore we arrive at
\[
E_{\wh{B},\psi}(s,f,g)=\int_{\wh{U}(F)\bs \wh{U}(\mb{A})}\wtl{f}_s(ug)\overline{\psi}(u)du.
\]
Since by definition $\wtl{f}_s$ is left invariant under the root subgroup $U_{\a_0}$, $E_{\wh{B},\psi}$ vanishes unless
$\psi_0$ is trivial, i.e. $\psi$ equals $\psi_1$ and is non-trivial.\hfill$\Box$

\

Now suppose that $\chi_{\wh{T}}$ is an unramified character on $\wh{T}(\mb{A})/\wh{T}(F)$ such that
\be\label{chisl}
\chi_{\wh{T}}(c\c a)=\chi_0(c)\chi_1(a),\q c\in \mb{I},\q a\in T(\mb{A})\h\mb{I},
\ee
where $\chi_0$ and $\chi_1$ are unramified Hecke quasi-characters on $\mb{I}$.
We define $E(\chi_{\wh{T}},g)$ by (\ref{chi}) and the Fourier coefficient $E_{\wh{B},\psi}(\chi,g)$. Similarly as before we have
\[
E_{\wh{B},\psi}(\chi_{\wh{T}},g)=\sum_{w\in\wtl{W}}\int_{w^{-1}\wh{G}(F\lg t \rg_+)w\cap
\wh{U}(F)\bs\wh{U}(\mb{A})}\wtl{\chi}_{\wh{T}}(wug)\o{\psi}(u)du.
\]
For any $w\in\wtl{W}$, at least one of $w\a_i$, $i=0,1$ is positive. Since $\wtl{\chi}_{\wh{T}}$ is left invariant under any root subgroup
of a positive root, we see that $E_{\wh{B},\psi}$ does not vanish unless one of $\psi_i$, $i=0,1$ is trivial, i.e. $\psi=\psi_0$ or $\psi_1.$
Let us work out the explicit formula for $E_{\wh{B},\psi_i}(\chi_{\wh{T}},g)$ at $g=\s(q),$ $i=0,1.$ For simplicity let us assume that
\be\label{simp}
\psi_i(\wtl{x}_{\a_i}(u))=\psi_F(u),\q u\in\mb{A},\q i=0,1
\ee
where $\psi_F=\bot\limits_v\psi_{F_v}$ and $\psi_{F_v}$ is the standard character of $F_v$ defined by (\ref{psir})-(\ref{psiv}). As preliminary computations let us give the local coefficients
for $SL_2(F_v).$ Define
\be \label{wchi}
W(s,\chi_v)=\int_{F_v}f_{s,\chi_v}\l(w^{-1}\bpm 1 &x \\ 0 &1 \epm \r)\o{\psi}_{F_v}(x)dx
\ee
where $f_{s,\chi_v}$ is given by (\ref{fchi}). The following proposition is well-known.

\begin{prop}\label{wtv}
Suppose $F_v$ is $p$-adic, then
\[
W(s,\chi_v)=\mathrm{vol}({\cal O}_v)L(s+1,\chi_v)^{-1}.
\]
\end{prop}

We also have the local coefficients at archimedean places. Assume $\chi_v=|\c|_v^{s_v}$, $v|\i$.
If $F_v=\mb{R}$, then
\bq
W(s,\chi_v)&=&2\int^\i_0 (1+x^2)^{-\f{s+s_v+1}{2}}e^{-2\pi i x}dx\label{wtr}\\
&=&2\G_\mb{R}(s+1,\chi_v)^{-1}K_{\f{s+s_v}{2}}(2\pi).\n
\eq
If $F_v=\mb{C}$, then
\bq
W(s, |\c|_{\mb{C}}^{s_0})&=&2\int_{\mb{R}^2}(1+x^2+y^2)^{-s-s_0-1}e^{-4\pi i x}dx dy\label{wtc}\\
&=& \rm{B}\l(s+s_0+\f{1}{2},\f{1}{2}\r)\int^\i_0(1+x^2)^{-s-s_0}e^{-4\pi ix}dx\n\\
&=& \f{1}{\sqrt{2}}\G_\mb{C}(s+1,\chi_v)^{-1}\rm{B}(s+s_v,\f{1}{2})^{-1}K_{s+s_v-\f{1}{2}}(4\pi).\n
\eq
In the above, $\rm{B}(\cdot,\cdot)$ is the Beta function, and $K_s(y)$ is the $K$\textit{-Bessel function}, also known as the \textit{Macdonald Bessel function}, defined by \[
K_s(y)=\f{1}{2}\int^\i_0 e^{-y(t+t^{-1})/2}t^s\f{dt}{t}.
\]
We have used the formulas in \cite[pp.66-67]{bump} to obtain (\ref{wtr}) and (\ref{wtc}). Define
\be\label{glow}
W(s,\chi)=\prod_v W(s,\chi_v).
\ee
Then $W(s,\chi)$ can be written as $|\D_F|^{-\f{1}{2}}W'_\i(s,\chi)L(s+1,\chi)^{-1}$, where
$W'_\i(s,\chi)$ is a product involving Bessel functions and Beta functions. Now we are ready to give the formula for the Fourier coefficients $E_{\wh{B},\psi_i}(\chi_{\wh{T}},\s(q))$. We do the
case $\psi=\psi_1$. The other case is similar, which is only a little bit more involved.

\begin{prop}
Assume $|\mf{q}|>1$ and the Godement's criterion
\[
\mf{Re}\chi_1>2,\q \mf{Re}(\chi_0-\chi_1)>2.
\]
Then
\bq
&&E_{\wh{B},\psi_1}(\chi_{\wh{T}},\s(\mf{q}))\n\\&=& \sum^\i_{n=1}|\D_F|^{-n+\f{1}{2}}|\mf{q}|^{n-2n^2}\l(\chi_1^n\chi_0^{-n^2}\r)(\mf{q})W(1-4n,\chi_0^{2n}\chi_1^{-1})\prod^{2n-1}_{i=1}\f{L(1-2i,\chi_0^i\chi_1^{-1})}
{L(2-2i,\chi_0^i\chi_1^{-1})}\n\\
&+&\sum^\i_{n=0}|\D_F|^{-n}|\mf{q}|^{-n-2n^2}\l(\chi_1^{-n}\chi_0^{-n^2}\r)(\mf{q})W(-1-4n,\chi_0^{2n}\chi_1)\prod^{2n}_{i=1}\f{L(1-2i,\chi_0^{i-1}\chi_1)}
{L(2-2i,\chi_0^{i-1}\chi_1)}.\n
\eq
\end{prop}
\textit{Proof.} By previous reasonings, we have
\[
E_{\wh{B},\psi_1}(\chi_{\wh{T}},g)=\sum_{w\in\wtl{W}, w\a<0}\int_{w^{-1}\wh{G}(F\lg t \rg_+)w\cap
\wh{U}(F)\bs\wh{U}(\mb{A})}\wtl{\chi}_{\wh{T}}(wu\s(q))\o{\psi}(u)du.
\]
For $w\in\wtl{W}$, recall our notations
\[
\wtl{\Phi}_w=\wtl{\Phi}_+\cap w\wtl{\Phi}_-=\{\b_1,\ld,\b_{l}\},
\]
\[\wtl{\Phi}_{w^{-1}}=\wtl{\Phi}_+\cap w^{-1}\wtl{\Phi}_-=\{\g_1,\ld,\g_{l}\}.
\]
It is clear from the formula for $\b_i$, $\g_i$ that if $w\a<0,$ then $\g_l=\a.$
Following the arguments of (\ref{int}) (\ref{int2}) we obtain
\bq
&&\int_{w^{-1}\wh{G}(F\lg t \rg_+)w\cap
\wh{U}(F)\bs\wh{U}(\mb{A})}\wtl{\chi}_{\wh{T}}(wu\s(q))\o{\psi}(u)du\label{intf}\\
&=&\int_{U_w(\mb{A})}\wtl{\chi}_{\wh{T}}(wu\s(q))\o{\psi}(u)du\n\\
&=&\int_{\mb{A}^{l}}\wtl{\chi}_{\wh{T}}\l(w\s(\mf{q})\wtl{x}_{\g_1}(\s(\mf{q})^{-\g_{1}}u_{1})\cs \wtl{x}_{\g_{l}}(\s(\mf{q})^{-\g_{l}}u_{l})\r)\o{\psi}_F(u_{l})du_1\cs du_{l}\n\\
&=& \s(\mf{q})^{\g_1+\cs+\g_{l}+w^{-1}\chi_{\wh{T}}}\int_{\mb{A}^{l}}\wtl{\chi}_{\wh{T}} \l(w\wtl{x}_{\g_1}(u_1)\cs \wtl{x}_{\g_{l}}(u_{l})\r)\o{\psi}_F(u_{l})du_1\cs du_{l}\n\\
&=& \s(\mf{q})^{w^{-1}(\chi_{\wh{T}}-\wtl{\rho})}\int_{\mb{A}^{l}}\wtl{\chi}_{\wh{T}}\l(\wtl{x}_{-\b_1}(u_1)\cs \wtl{x}_{-\b_{l}}(u_{l})\r)\o{\psi}_F(u_{l})du_1\cs du_{l}.\n
\eq
Notice that $\s(\mf{q})^{\g_{l}}=\s(\mf{q})^\a=1.$ Similarly as (\ref{int3}),
\bq
&&\int_{\mb{A}^{l}}\wtl{\chi}_{\wh{T}}\l(\wtl{x}_{-\b_1}(u_1)\cs \wtl{x}_{-\b_{l}}(u_{l})\r)\o{\psi}_F(u_{l})du_1\cs du_{l}\n\\
&=&\int_{\mb{A}}a(u_{l})^{\chi_{\wh{T}}-\wtl{\rho}+w'\wtl{\rho}}\o{\psi}_F(u_l)du_l\int_{\mb{A}^{l-1}}\wtl{\chi}_{\wh{T}}\l(\wtl{x}_{-\b_1}(u_1)\cs \wtl{x}_{-\b_{l-1}}(u_{l-1})\r)du_1\cs du_{l-1}.\n
\eq
Applying the formula of Fourier coefficients for $SL_2$, together with Gindikin-Karpelevich formula, we see that the last equation equals
\[
|\D_F|^{-\f{l-1}{2}}W(-\lg\wtl{\rho},\b_l^\vee\rg,\chi_{\wh{T}}\b_l^\vee)\prod^{l-1}_{i=1}\f{L(-\lg\wtl{\rho},\b_i^\vee\rg,\chi_{\wh{T}}\b_i^\vee)}
{L(1-\lg\wtl{\rho},\b_i^\vee\rg,\chi_{\wh{T}}\b_i^\vee)}.
\]
There are two cases for $w\in\wtl{W}$ such that $w\a<0.$

\textit{Case} 1: $w=T_{n\a^\vee},$ $n>0.$  Then $l=2n$, $\b_i=i\de-\a$ and $\g_i=(2n-i)\de+\a.$ In this case we obtain
\[
\lg\wtl{\rho},\b_i^\vee\rg=\lg \f{1}{2}\a+2L, ic-\a^\vee\rg=2i-1,
\]
\[
wd=T_{n\a^\vee}d=d+n\a^\vee-n^2c,
\]
\[
\chi_{\wh{T}}\b_i^\vee=\chi_0^i\chi_1^{-1},\quad\chi_{\wh{T}} wd=\chi_1^n\chi_0^{-n^2},\q \lg\wtl{\rho},wd\rg=n-2n^2.
\]

\textit{Case} 2: $w=T_{-n\a^\vee}r_\a$, $n\geq 0.$ Then $l=2n+1$, $\b_i=(i-1)\de+\a,$ $\g_i=(2n+1-i)\de+\a.$ Similarly we obtain
\[
\lg\wtl{\rho},\b_i^\vee\rg=\lg \f{1}{2}\a+2L, (i-1)c+\a^\vee\rg=2i-1,
\]
\[
wd=T_{-n\a^\vee}r_\a d= T_{-n\a^\vee}d=d-n\a^\vee-n^2c,
\]
\[
\chi_{\wh{T}}\b_i^\vee=\chi_0^{i-1}\chi_1,\quad\chi_{\wh{T}}wd=\chi_1^{-n}\chi_0^{-n^2},\q \lg\wtl{\rho},wd\rg=-n-2n^2.
\]
Combining contributions from these two cases, we get the formula.\hfill$\Box$

\section{Absolute Convergence of Eisenstein Series}

Under the conditions of Theorem \ref{conv} we have proved that $E(s,f,ug)$ converges absolutely almost everywhere on $\wh{U}(F)\bs\wh{U}(\mb{A})$, by proving the finiteness of the constant term $E_{\wh{B}}(s,h,g).$
The main result of this section is the absolute and uniform convergence of $E(s,f,g)$ for $g$ in certain Siegel set and $\mf{\mf{Re}}s$ large enough. By boundedness of the cusp form $f$, it is enough to prove the absolute convergence
of $E(s,h,g).$ The main ingredient of the proof is the systematic use of Demazure modules together with some technical estimations.
We follow Garland's idea in \cite{Gar5}. Our arguments in the adelic settings also involve a property
of algebraic number fields, which is analogous to the Riemann-Roch theorem for algebraic curves.

\subsection{Demazure modules}

Recall that for any dominant integral weight $\la,$ we have the
irreducible highest weight module $V_{\la},$ and the highest weight
vector $v_\la$ such that $V_{\la,\la,\mb{Z}}=\mb{Z}v_\la.$ There is
a map $|\c|: V_{\la,\mb{A}}\ra\mb{R}_{\geq 0}$ defined in
section 3.4. The highest weight vector $v_\la$ embeds into
$V_{\la,\mb{A}}$ diagonally. Recall that $L$ is the fundamental weight such that $\lg
L,\a_i^\vee\rg=0,$ $i=1,\ld, n$, and $\lg L, \a_0^\vee\rg=\lg L,
c-\wtl{\a}^\vee\rg=1.$ Then by Lemma \ref{hight} it follows that the
height function $h_s$ can be defined as \be \label{hs}
h_s(g)=|g^{-1}v_L|^{-s}, \ee where $v_L$ is the highest weight
vector in $V_{L,\mb{A}}$ as above.

The simple equation (\ref{hs}) plays a crucial role in the proof of the absolute convergence of Eisenstein series. In this section we shall collect some basic results on the Demazure module $V_{\la}(w)$, which is associated with a dominant integral weight $\la$ and $w\in\wtl{W}$, and is a submodule of $V_{\la}$ defined below.

Recall that $\wtl{\mf{n}}_+=\bop\limits_{\a\in\wtl{\Phi}_+}\mf{g}_\a$ is the Lie algebra of $\wh{U}.$ Let ${\cal
U}^{\mb{Z}}(\wtl{\mf{n}}_+)$ be a $\mb{Z}$-form of the universal
enveloping algebra ${\cal U}(\wtl{\mf{n}}_+)$ of $\wtl{\mf{n}}_+$. We
define
\be\label{dem}
V_{\la}(w)_\mb{Z}={\cal U}^{\mb{Z}}(\wtl{\mf{n}}_+)\c w\c v_{\la}.
\ee
Let $\La(w)$ be the subset of all weights $\mu$ of $V_{\la}$
such that $\mu\geq w\la,$ i.e. $\mu-w\la=\sum\limits_{i=0}^n
l_i\a_i,$ $l_i\in\mb{Z}_{\geq 0}.$  Then it is clear that
\be\label{wti}
V_\la(w)_\mb{Z}\ss \bop_{\mu\in\La(w)}V_{\la,\mu,\mb{Z}}.
\ee

For any field $F$ we define the Demazure module over $F$
corresponding to $\la$ and $w$ by $V_{\la}(w)_F=V_{\la}(w)_\mb{Z}\ot
F.$ Let $V_\la(w)_{\mb{A}}\ss V_{\la,\mb{A}}$ be the restricted product
\[
\prod'_v V_\la(w)_{F_v}
\]
with respect to $V_\la(w)_{{\cal O}_v}=V_\la(w)_{\mb{Z}}\ot{\cal O}_v$ defined for all finite places $v.$ If $\phi=\ot_v\phi_v$ such that $\phi_v$
is a linear operator on $V_{\la,F_v}$ and $\phi_v$ preserves $V_{\la,{\cal O}_v}$ for almost all finite places $v$, then we define
\be\label{glon}
\|\phi_v\|=\sup_{x\in V_{\la,F_v}, \|x\|=1}\|\phi_v x\|,\q \|\phi\|=\prod_{v}\|\phi_v\|,
\ee
and
\be\label{gloa}
\q |\phi_v|=\sup_{x\in V_{\la,F_v}, \|x\|=1}|\phi_vx|,\q |\phi|=\prod_v|\phi_v|.
\ee
Note that $|\phi_v|=\|\phi_v\|^2$ or $\|\phi_v\|$ according to $v$ is complex or not. If $\phi_v$ preserves $V_\la(w)_{F_v}$ for each $v$, then we define $\|\c\|_w$ and $|\c|_w$ similarly by restriction. In particular
$\|\phi_v\|_w\leq \|\phi_v\|,$ $|\phi_v|_w\leq |\phi_v|.$

Let $F$ be a local field. We shall estimate the norm of $\wtl{x}_\b(u)$ acting on
$V_{\la}(w)_F,$ where $\b\in\wtl{\Phi}_{re+},$ $u\in F.$ By Lemma \ref{sl2} we have the group homomorphism
$\vp_\b: SL_2(F)\ra \wh{G}(F((t)))$ such that
\[
\vp_\b\l(\bpm 1 & u\\ 0 &1\epm\r)=\wtl{x}_\b(u).
\]
Then it follows that
\[
\vp_\b\l(\bpm u & 0\\ 0 &u^{-1}\epm\r)=\wtl{h}_\b(u),\q u\in
F^\times,
\]
and $\vp_\b$ maps the maximal compact subgroup of $SL_2(F)$ into
$\wh{K}.$

If $F=\mb{R}$ or $\mb{C}$, consider the Cartan decomposition
\be\label{car}
\bpm 1 & u\\0 & 1\epm\bpm 1 & 0\\ \overline{u} & 1\epm=k\bpm a &0\\ 0
&a^{-1}\epm k^{-1},
\ee
where $k$ lies in the maximal compact subgroup of $SL_2(F),$ and
$a\geq a^{-1}>0.$ Compare the trace and determinant of both sides it
is easy to obtain
\be \label{au}
a=a(u)=\f{|u|^2+2+|u|\sqrt{|u|^2+4}}{2}.
\ee

\begin{lemma} \label{norm} $($\cite{Gar5}~Lemma~$4.1)$
Suppose $F=\mb{R}$ or $\mb{C}$, $\b\in\wtl{\Phi}_{re+},$ $u\in F$. The square norm $\|\wtl{x}_\b(u)\|_w^2$ of $\wtl{x}_\b(u)$ restricted on $V_\la(w)_F$
, is bounded by
\[
\sup_{\mu\in\rm{weights of }V_\la(w)} a(u)^{\lg\mu,\b^\vee\rg},
\]
where $a(u)$ is given by $(\ref{au})$.
\end{lemma}
\textit{Proof.} The adjoint operator of $\wtl{x}_\b(u)$ is $\wtl{x}_{-\b}(\overline{u})$, which follows from
the facts that $X_\a\ot t^i$ and $X_{-\a}\ot t^{-i}$ ($\a\in\Phi$) are adjoint operators, and that the inner product is hermitian.
Then by above discussions we get
\bq
\|\wtl{x}_\b(u)\|^2_w&=&\sup_{v\in V_\la(w)_F, \|v\|=1}\|\wtl{x}_\b(u)v\|^2\n \\
&=& \sup_{v\in V_\la(w)_F, \|v\|=1} (\wtl{x}_{-\b}(\overline{u})\wtl{x}_\b(u)v,v)\n\\
&=& \sup_{v\in V_\la(w)_F, \|v\|=1}(k_\b\wtl{h}_\b(a)k_\b^{-1}v,v)\n\\
&=& \sup_{v\in V_\la(w)_F, \|v\|=1}(\wtl{h}_\b(a)v,v)\n\\
&\leq & \sup_{\mu\in \rm{weights of }V_\la(w)}a^{\lg\mu,\b^\vee\rg},\n
\eq
where $a=a(u)$ and $k_\b=\vp_\b(k)\in\wh{K}$, with $k$ the element appearing in (\ref{car}).
\hfill$\Box$

\begin{lemma}\label{pnorm}
Suppose $F$ is a $p$-adic field, $\b=\a+i\de\in\wtl{\Phi}_{re+},$ $u\in F$. Then $|\wtl{x}_\b(u)|_w=1$ if $|u|\leq 1$, and is bounded
by \[\sup_{\mu\in\rm{weights of }V_\la(w)} |u|^{2|\lg\mu,\a^\vee\rg|}\]
if $|u|>1.$
\end{lemma}
\textit{Proof.} If $|u|\leq 1$ the lemma is clear since $\wtl{x}_\b(u)\in\wh{K}$ which preserves the norm on $V_{\la,F}.$ Assume
that $|u|>1.$ Then
\bq
|\wtl{x}_\b(u)|_w&=&|\wtl{h}_\a(u)\wtl{x}_\b(u^{-1})\wtl{h}_\a(u^{-1})|_w\n\\
&\leq & |\wtl{h}_\a(u)|_w |\wtl{h}_\a(u^{-1})|_w\n\\
&\leq & \sup_{\mu\in\rm{weights of }V_\la(w)}|u|^{2|\lg\mu,\a^\vee\rg|}.\n
\eq
\hfill$\Box$

Consider the real Cartan subalgebra and its dual,
\[
\wtl{\mf{h}}_\mb{R}=\mf{h}_\mb{R}\op \mb{R}c\op \mb{R}d,\quad
\wtl{\mf{h}}_\mb{R}^\ast=\mf{h}_\mb{R}^\ast\op \mb{R}\de\op \mb{R}L.
\]
For any $\mu\in \wtl{\mf{h}}_\mb{R}^\ast,$ write $\mu$ as the
decomposition \be \label{decom} \mu=\mu_0-\k_\mu\de+\s_\mu L, \ee
where $\mu_0\in\mf{h}_{\mb{R}}^\ast,$ $\k_\mu,$ $\s_\mu\in\mb{R}.$
If $\mu$ is a weight of $V_\la(w),$ then $\mu\geq w\la$ and
$\k_\mu\leq \k_{w\la}.$ We impose the condition $\k_\la=0$ in order
that $\la$ be dominant integral. The following lemmas are due to
Garland \cite{Gar5}. For reader's convenience we shall sketch the
proof.

\begin{lemma}\label{mu0}
Let $\la\in \wtl{\mf{h}}_\mb{R}^\ast $ be a dominant integral weight $\la.$ Then there exists a constant $\k_0>0$ such that
for all $w\in\wtl{W},$
\[
\k_{w\la}\leq \k_0l(w)^2
\]
and for any weight $\mu$ of $V_\la(w)$, $w\neq 1,$
\[
\left(\mu_0,\mu_0\right)\leq \k_0 l(w)^2.
\]
\end{lemma}
\textit{Proof.} We only prove the first inequality. See \cite{Gar5} for the proof of the second one. Let us write
$\la=\la_0+\s_\la L,$ where $\la_0\in\mf{h}^\ast,$ $\s_\la\in\mb{N}.$ Write $w=T_\g w_0$ where $\g\in Q^\vee,$ $w_0\in W.$
Then by Lemma \ref{weyl2}
\[
w\la= w_0\la_0+\s_\la L+\s_\la \g-\l(\lg\la_0,\g\rg+\f{\s_\la}{2}\left(\g,\g\right)\r)\de.
\]
It follows from Lemma \ref{len} below that
\[
\k_{w\la}=\lg\la_0,\g\rg+\f{\s_\la}{2}\left(\g,\g\right)= O(1)\|\g\|^2= O(1)l(T_\g)^2\leq  O(1)(l(w)+|\Phi_+|)^2.
\]
From this the first inequality is clear. Here we denote by $O(1)$ a bounded term.
\hfill$\Box$

\

As a consequence of this lemma we obtain the following result.

\begin{cor}\label{F1}
Given $\k>0,$ there exists $\k_1>0$ such that for all $w\in\wtl{W},$ $w\neq 1,$ and for all
$\b\in\wtl{\Phi}$ of the form
\[
\b=\a+i\de,\q 0\leq i\leq \k l(w),~\a\in\Phi,
\]
we have $\lg\mu,\b^\vee\rg\leq \k_1l(w)$ for any weight $\mu$ of $V_\la(w).$
\end{cor}
\textit{Proof.} Write $\mu=\mu_0-\k_\mu\de+\s_\la L.$ By Lemma \ref{coroot},
$\b^\vee=\f{i}{2}\left(\a^\vee,\a^\vee\right) c+\a^\vee.$ Then from Lemma \ref{mu0} it follows
that
\[
\lg\mu,\b^\vee\rg=\lg\mu_0,\a^\vee\rg+\f{i}{2}\left(\a^\vee,\a^\vee\right)\s_\la\leq \l(\k_0^{\f{1}{2}}\|\a^\vee\|+\f{\k}{2}\|\a^\vee\|^2\s_\la\r)l(w).
\]
Let $\k_1$ be the maximum of above coefficient of $l(w)$ over $\a\in\Phi.$
\hfill$\Box$

\

The condition on $\b$ in Corollary \ref{F1} is satisfied in the following case.

\begin{lemma}\label{kappa}
There exists $\k\in\mb{N}$ such that for all $w\in\wtl{W}$ and
$\b=\a+i\de\in \wtl{\Phi}_w=\wtl{\Phi}_+\cap w\wtl{\Phi}_-$
$(\a\in\Phi, i\in\mb{Z}_{\geq 0})$, we have $0\leq i<\k l(w).$
\end{lemma}
\textit{Proof.} Write $w^{-1}=T_\g w_0$ where $\g\in Q^\vee,$ $w_0\in W.$ Then
\[
w^{-1}\b=w_0\a+(i-\lg w_0\a, \g\rg)\de<0,
\]
which implies that $i\leq \lg w_0\a, \g\rg\leq \|\g\|\|\a\|.$ The proof follows from the
following lemma, together with the inequality $l(T_\g)\leq l(w)+l(w_0)\leq l(w)+|\Phi_+|.$
\hfill$\Box$

\begin{lemma}\label{len}
There exists $\wtl{\k}>0$ such that $\|\g\|\leq \wtl{\k}l(T_\g)$ for all $\g\in Q^\vee$.
\end{lemma}
\textit{Proof.} From \cite{IM} Proposition 1.23 we get
\[
l(T_\g)=\sum_{\a\in\Phi_+}|\lg\a,\g\rg|\geq\sum^{n}_{i=1}|\lg\a_i,\g\rg|.
\]
$\g$ can be written as a linear combination of the fundamental coweights, with coefficients $\lg\a_i,\g\rg$, $1\leq i\leq n$. The lemma is clear from this observation.\hfill$\Box$

\

The Demazure module $V_\la(w)_F$ is preserved by the action of the elements $\wtl{x}_\a(u)$, $\a\in\Phi,$ $u\in tF[[t]],$ and
$\wtl{h}_\a(u),$ $\a\in\Phi,$ $u\equiv 1\mod tF[[t]].$

\begin{lemma}\label{tri}
Let $\k$ be given as in Lemma \ref{kappa}. For any $w\in\wtl{W},$ the elements $\wtl{x}_\a(u),$ $\a\in\Phi,$ $u\in t^{\k l(w)}F[[t]]$
and $\wtl{h}_\a(u)$, $\a\in \Phi,$ $u\equiv 1\mod t^{\k l(w)}$ act as identity on the Demazure module $V_\la(w)_F.$
\end{lemma}
\textit{Proof.} Let $\mf{u}_{w^{-1}}=\bop\limits_{\a\in
\wtl{\Phi}_w}\mf{g}_\a$ be the Lie algebra of $U_{w^{-1}},$ then
\[
V_\la(w)_\mb{Z}={\cal U}^\mb{Z}(\mf{u}_{w^{-1}})wv_\la.
\]
Let $\wtl{\Phi}_w=\{\b_1,\ld,\b_l\}$ with $l=l(w).$ The PBW Theorem
implies that the monomials $X_{\b_1}^{i_1}\cs X_{\b_l}^{i_l},$
$i_1,\ld,i_l\in\mb{Z}_{\geq 0}$, form a basis of ${\cal
U}(\mf{u}_{w^{-1}}).$ To prove the lemma, it suffices to show
that $U_\b(F)$ acts on $V_\la(w)_F$ trivially
for each $\b=\a+i\de$ with $\a\in\Phi\cup\{0\}$, $i\geq\k l(w).$  Then
we are further reduced to show that $\mf{g}_\b=\mf{g}_{\a+i\de}$ with
$\a\in\Phi\cup\{0\}$, $i\geq \k l(w)$ acts on $V_\la(w)_F$ as zero. We prove
by induction on $i_1+\cs+i_l$ that \be \label{ind} \mf{g}_\b
X_{\b_1}^{i_1}\cs X_{\b_l}^{i_l}wv_\la=0. \ee Since $w^{-1}\mf{g}_\b
w=\mf{g}_{w^{-1}\b}$ and $w^{-1}\b>0$ by Lemma \ref{kappa}, we have
$\mf{g}_\b wv_\la=0.$ Consider $[\mf{g}_\b,\mf{g}_{\b_i}],$ which is zero if
$\b+\b_i\not\in \wtl{\Phi}$ and equals $\mf{g}_{\b+\b_i}$
otherwise. In each case the induction follows and (\ref{ind}) is
proved.
\hfill$\Box$

\subsection{Estimations of some norms}

Let $F$ be a local field. In this section we shall apply the results of the previous section to estimate
the norms of elements in $\wh{U}(F)$ acting on $V_\la(w)_F$, under certain conditions.

\begin{lemma} \label{est1} $($\cite[pp.228\rm{-}232]{Gar5}$)$
Suppose that $F=\mb{R}$ or $\mb{C}$, $\wtl{x}_\a(u)\in\wh{U}(F)$, where $\a\in\Phi,$ $u=\sum\limits^\i_{i=0}u_i t^i\in F[[t]]$ $(\in tF[[t]]$ if $\a\in\Phi_-)$ such that $|u_i|\leq M \tau^i$, $i=0,1,\ld$ for some $M>0$
and $0<\tau<1.$ Then $\|\wtl{x}_\a(u)\|_w\leq\exp(\k_{M,\tau}l(w))$ for some constant
$\k_{M,\tau}$ only depending on $M$ and $\tau.$
\end{lemma}
\textit{Proof.} Consider $\b=\a+i\de\in\wtl{\Phi}_{re+},$ $u_i\in F.$ Let $\k\in\mb{N}$ be the constant in Lemma \ref{kappa}. If $i\geq \k l(w)$ then $\wtl{x}_\b(u_i)$ acts on $V_\la(w)_F$ trivially; if $i<\k l(w),$ by Lemma \ref{norm} and Corollary \ref{F1} we have
\[
\|\wtl{x}_\b(u_i)\|_w^2  \leq \sup_{\mu\in\rm{weights of }V_\la(w)}a(u)^{\lg\mu,\b^\vee\rg}\leq  a(u_i)^{\k_1 l(w)},
\]
where $a(u_i)\geq 1$ is given by (\ref{au}). It is easy to show that there exists $c_M>0$ depending on $M$ such that
$a(u_i)\leq 1+ c_M \tau^i.$ Then by Lemma \ref{tri} we have
\bq
\|\wtl{x}_\a(u)\|_w^2&=&\|\prod^{\k l(w)-1}_{i=0}\wtl{x}_{\a+i\de}(u_i)\|_w^2\leq  \prod^{\k l(w)-1}_{i=0}\|\wtl{x}_{\a+i\de}(u_i)\|_w^2\n\\
&\leq & \prod^{\k l(w)-1}_{i=0}(1+ c_M \tau^i)^{\k_1 l(w)}\leq  \exp\l(\k_1c_M l(w) \sum ^{\k l(w)-1}_{i=0}\tau^i\r)\n\\
&\leq & \exp\l(\f{\k_1c_M}{1-\tau}l(w)\r).\n
\eq
The lemma follows if we set $\k_{M,\tau}=\df{\k_1c_M}{2(1-\tau)}.$
\hfill$\Box$

\begin{lemma} \label{est2} $($\cite[pp.233\rm{-}240]{Gar5}$)$
 Suppose that $F=\mb{R}$ or $\mb{C}$, $u=1+\sum\limits^\i_{j=1}u_j t^j\in 1+tF[[t]]$ such that $|u_j|<M\tau^j$, $j=1,2,\ld$ for some $M>0$ and $0<\tau<1$. Then $\|\wtl{h}_{\a_i}(u)\|_w\leq \exp(\wtl{\k}_{M,\tau} l(w)),$ $i=0,1,\ld, n,$
for some constant $\wtl{\k}_{M,\tau}$ only depending on $M$ and $\tau.$
\end{lemma}
\textit{Proof.} Let us consider the following two cases:

\textit{Case} 1: $w^{-1}\a_i<0.$ Let $w'=\wtl{w}_{\a_i}w,$ then $l(w)=l(w')+1.$ Moreover if we write $\wtl{\Phi}_{w'}=\wtl{\Phi}_+\cap w'\wtl{\Phi}_-=\{\b_1,\ld,\b_{l-1}\},$ then $\wtl{\Phi}_w=\{\a_i, r_i\b_1, \ld, r_i\b_{l-1}\}.$
We shall prove the following:
\be\label{xinv}
X_{-\a_i}V_\la(w)_F\ss V_\la(w)_F,
\ee
\be\label{ss}
V_\la(w')_F\ss V_\la(w)_F.
\ee
To prove (\ref{xinv}) it suffices to prove that for $j, j_1,\ld, j_{l-1}\in\mb{Z}_{\geq 0}$,
\[
X_{-\a_i}X_{\a_i}^j X_{r_i\b_1}^{j_1}\cs X_{r_i \b_{l-1}}^{j_{l-1}}wv_\la\in V_\la(w)_F.
\]
We use induction on $j+j_1+\ld+j_{l-1}.$ It is clear by assumption that $X_{-\a_i}wv_\la=0.$
The induction follows since $[X_{-\a_i}, X_{\a_i}], [X_{-\a_i}, X_{r_i\b_1}],\ld, [X_{-\a_i}, X_{r_i\b_{l-1}}]\in \wtl{\mf{b}}
=\wtl{\mf{h}}\op \wtl{\mf{n}}_+.$
$V_\la(w)_F$ is also preserved by $X_{\a_i},$ hence it is preserved by $w_{\a_i}.$ Then (\ref{ss}) follows from
\[X_{\b_1}^{j_1}\cs X_{\b_{l-1}}^{j_{l-1}}w'v_\la= w_{\a_i}(w_{\a_i}X_{\b_1}w_{\a_i})^{j_1}\cs (w_{\a_i}X_{\b_{l-1}}w_{\a_i})^{j_{l-1}}wv_\la\in V_\la(w)_F.\]
Since $w_{\a_i}=w_{\a_i}(1)\in\wh{K}$ which preserves the norm, by Lemma \ref{est1} and its proof we obtain
\bq
\|\wtl{h}_{\a_i}(u)\|_w&=& \| w_{\a_i}(u)w_{\a_i}(1)^{-1}\|_w= \| w_{\a_i}(u)\|_w\n\\
&\leq & \|\wtl{x}_{\a_i}(u)\|_w\|\wtl{x}_{-\a_i}(-u^{-1})\|_w\|\wtl{x}_{\a_i}(u)\|_w\leq \exp(3\k_{M,\tau}l(w)).\n
\eq

 \textit{Case }2: $w^{-1}\a_i>0.$ Then thanks to (\ref{xinv}) and (\ref{ss}), $V_\la(w')_F$ is preserved by $X_{-\a_i}$ and $\wtl{w}_{\a_i}$, and we have
 $V_\la(w)_F\ss V_\la(w')_F.$ Apply the result in Case 1 we get
 \[
 \|\wtl{h}_{\a_i}(u)\|_w\leq \|\wtl{h}_{\a_i}(u)\|_{w'}\leq \exp(3\k_{M,\tau}(l(w)+1)).
 \]
 Combine the two cases, the lemma holds for $\wtl{\k}_{M,\tau}=3\k_{M,\tau}+1.$
 \hfill$\Box$

 \

We have the following $p$-adic analog of Lemma \ref{est1}.

\begin{lemma}\label{estp}
 Suppose that $F$ is a $p$-adic field, $\wtl{x}_\a(u)\in\wh{U}$, where $\a\in\Phi,$ $u=\sum\limits^\i_{i=0}u_i t^i\in F[[t]]$ $(\in tF[[t]]$ if $\a\in\Phi_-)$ such that $|u_i|<M$, $i=0,1,\ld$ for some $M>0.$ Then
$|\wtl{x}_\a(u)|_w\leq\exp(\k_M l(w))$ for a constant $\k_M$ only depending on $M$.
\end{lemma}
\textit{Proof.} We may assume that $M>1.$ By Lemma \ref{mu0} we have
\[
|\lg\mu,\a^\vee\rg|=|\lg \mu_0, \a^\vee\rg|\leq \k_0^{\f{1}{2}}\|\a^\vee\|l(w)\leq\k_1l(w)
\]
for any weights $\mu$ of $V_\la(w).$ Following the proof of Lemma \ref{pnorm} we obtain
\bq
|\wtl{x}_\a(u)|_w&=& |\wtl{h}_\a(M)\wtl{x}_\a(M^{-2}u)\wtl{h}_\a(M^{-1})|_w\leq  |\wtl{h}_\a(M)|_w |\wtl{h}_\a(M^{-1})|_w\n\\
&\leq & \sup_{\mu\in \rm{weights of }V_\la(w)}|M|^{2|\lg\mu,\a^\vee\rg|}\leq  |M|^{2\k_1l(w)},\n
\eq
note that $\wtl{x}_\a(M^{-2}u)\in \wh{K}.$ The lemma follows by setting $\k_M=2\k_1\log M.$
\hfill$\Box$

\

Let $F$ be any local field. Recall (\ref{sac}) that $\s(\mf{q}),$ $\mf{q}\in F^\times$ acts on $V_{\la,F}$ by
\[
\s(\mf{q})v=\mf{q}^{\lg\mu,d\rg}v
\]
for each $v\in V_{\la,\mu,F}.$ Then
\[
\|\s(\mf{q})\|_w\leq \sup_{\mu\in\rm{weights of }V_\la(w)}|\mf{q}|^{\lg\mu,d\rg}\leq\sup_{\mu\in\La(w)}|\mf{q}|^{\lg\mu,d\rg}.
\]
If $\la=L,$ then $\lg\mu,d\rg\leq 0$ for any $\mu\in \La(w).$ In this case we obtain the following:

\begin{lemma}\label{sn}
Assume $\la=L,$ $\mf{q}\in F^\times.$ Then
\[
\|\s(\mf{q})\|_w=\l\{\ba{ll} 1,&\rm{ if }|\mf{q}|\geq 1,\\ |\s(\mf{q})^{wL}|, & \rm{ if }|\mf{q}|\leq 1.\ea\r.
\]
\end{lemma}

From this lemma we can get the following global result.

\begin{cor}\label{gsn}
Assume $\la=L,$ $\mf{q}=(\mf{q}_v)_v\in\mb{I} $. Then
$|\s(\mf{q})|_w=1$ if $|\mf{q}_v|_v\geq 1$ for each $v$, and $|\s(\mf{q})|_w=\s(\mf{q})^{wL}$ if $|\mf{q}_v|_v\leq 1$ for each $v.$
\end{cor}
\textit{Proof.} This is clear from the fact $|\s(\mf{q})|_w=\prod_v |\s(\mf{q}_v)|_w.$ \hfill$\Box$

\subsection{Convergence of Eisenstein series}

In this section we prove the absolute convergence of Eisenstein series everywhere, whenever the conditions
of Theorem \ref{conv} are satisfied. The uniform convergence of Eisenstein series over certain analog of Siegel set will
be established. The main result is the following:

\begin{thm}\label{convae}
Fix $\mf{q}\in\mb{I},$ $|\mf{q}|>1.$ There exists a constant $c_\mf{q}>0$
depending on $\mf{q}$, such that for any $\v>0$ and compact subset
$\Om$ of $T(\mb{A})$, $E(s,f,g)$ and $E(s,h,g)$ converge absolutely
and uniformly for $s\in \{z\in\mb{C}|\mathrm{\mf{Re}}z>\max(h+h^\vee+\v,
c_\mf{q})\}$ and $g\in \wh{U}(\mb{A})\Om \s(\mf{q}) \wh{K}.$
\end{thm}
\textit{Proof.} We only need to prove the theorem for $E(s,h,g)$. Let
us write $g=u_ga_g\s(\mf{q})k_g,$ where $u_g\in\wh{U}(\mb{A}),$
$a_g\in\Om,$ $|\mf{q}|>1$ and $k_g\in \wh{K}.$ Since $E(s,h,g)$ is left
$\wh{G}(F\lg t \rg)$-invariant, by Lemma \ref{fd} we may assume that
$u_g\in \wh{U}_{\cal D}.$ We may also assume that $k_g=1.$ Recall
(\ref{hs}) that $h_s$ is the height function such that
$h_s(g)=|g^{-1}v_L|^{-s}$, and
\[
E(s,h,g)=\sum_{w\in W(\D,\emptyset)}\sum_{\g\in U_w(F)}h_s(w\g g).
\]

Let $C>1$ be a constant which will be determined later. Write $\mf{q}=\mf{q}_1\mf{q}_2^{-1}$ such that \\
($i$) $\mf{q}_1, \mf{q}_2\in\mb{I}_+:=\{x=(x_v)_v\in\mb{I}|~|x_v|_v\geq1, \forall v\},$\\
($ii$) $|\mf{q}_{1v}|\geq C$ for each $v|\i.$
\\
By Lemma \ref{gsn}, $|\s(\mf{q}_2^{-1})|_{w^{-1}}=\s(\mf{q}_2^{-1})^{w^{-1}L}$.
Assume $\g\in U_w(F)$. Let $g_1=u_ga_g\s(\mf{q}_1)=g\s(\mf{q}_2)$. Since $(w\g g_1)^{-1}v_L\in V_L(w^{-1})_{\mb{A}}$, we have
\bq
\label{trick} h_1(w\g g)&=& |(w\g g)^{-1}v_L|^{-1}=|\s(\mf{q}_2)(w\g g_1)^{-1}v_L|^{-1}\\
&\leq& |\s(\mf{q}_2^{-1})|_{w^{-1}}|(w\g g_1)^{-1}v_L|^{-1}=\s(\mf{q}_2^{-1})^{w^{-1}L}h_1(w\g g_1).\n
\eq Similarly for any $u\in\wh{U}(\mb{A})$ we have
$(w\g u g_1)^{-1}v_L\in V_L(w^{-1})_{\mb{A}}.$ Note that $g_1$ acts on $V_L(w^{-1})_{\mb{A}}$, therefore
\bq
\label{comp} h_1(w\g g_1)&=& |g_1^{-1}\g^{-1}w^{-1}v_L|^{-1}\\
&\leq & |g_1^{-1}u^{-1}g_1|_{w^{-1}}\c |(w\g ug_1)^{-1}v_L|^{-1}\n\\
&=& |g_1^{-1}u^{-1}g_1|_{w^{-1}} h_1(w\g ug_1).\n
\eq
 Assume $u\in\wh{U}_{\cal D}$, and we shall estimate
$|g_1^{-1}u^{-1}g_1|_{w^{-1}},$ which is bounded by
\[
|(a_g\s(\mf{q}_1))^{-1}u_g^{-1}(a_g\s(\mf{q}_1))|_{w^{-1}}|(a_g\s(\mf{q}_1))^{-1}u^{-1}(a_g\s(\mf{q}_1))|_{w^{-1}}|(a_g\s(\mf{q}_1))^{-1}u_g(a_g\s(\mf{q}_1))|_{w^{-1}}.
\]
Let us estimate $|(a_g\s(\mf{q}_1))^{-1}u^{-1}(a_g\s(\mf{q}_1))|_{w^{-1}}$.
Other factors can be treated similarly. Since $u^{-1}\in \wh{U}_{\cal D},$ $u^{-1}$ is a product of the elements $\wtl{x}_{\a}(u_\a)$ where either $\a\in\Phi_+,$
$u_\a\in -{\cal D}\lg t \rg_+$ or $\a\in\Phi_-$, $u_\a\in -t{\cal D}\lg t \rg_+,$ and the elements $\wtl{h}_{\a_i}(u_i),$ $u_i\in
(1+t{\cal D}\lg t \rg_+)^{-1}$, $i=1,\ld, n.$

By our choice of ${\cal D}$, there exists $M_{\cal D}>0$ such that for any $x=(x_v)_v
\in{\cal D}$ we have $|x_v|\leq M_{\cal D}$ for each $v|\i.$ Then for any $\eta_{\cal D}>2M_{\cal D}$, there exists $M>0$
such that if $x=1+\sum\limits_{j=1}^\i x_jt^j\in (1+t{\cal D}\lg t \rg_+)^{-1}$, then $|x_{jv}|\leq M\eta_{\cal D}^j$, $j=1,2,\ld$ for
 each $v|\i$. In fact $x\in (1+t{\cal D}\lg t \rg)^{-1}$ implies that $x^{-1}_v$ defines a non-vanishing series absolutely convergent in the range
 $|t|<(2M_{\cal D})^{-1}$, hence so is $x_v$ itself. This implies
 \[
 \f{1}{\lim\sup\sqrt[j]{|x_{jv}|}}\geq \f{1}{2M_{\cal D}}>\f{1}{\eta_{\cal D}},
 \]
whence the assertion follows.

Consider
 \[(a_g\s(\mf{q}_1))^{-1}\wtl{h}_{\a_i}(u_i)(a_g\s(\mf{q}_1))=\wtl{h}_{\a_i}(\s(\mf{q}_1^{-1})\c u_i).\]
 Let $C>1$ be any constant such that
 \be\label{cst}
 C>\eta_{\cal D}>2M_{\cal D}.
 \ee
Applying Lemma \ref{est2} with $M$ chosen as above and $\tau=C^{-1}\eta_{\cal D}$, we get a constant $\k_{C, {\cal D}}:=2\wtl{\k}_{M,\tau}$ such that
\[
|\wtl{h}_{\a_i}(\s(\mf{q}_{1}^{-1})\c u_{i})_v|_{w^{-1}}\leq \exp\l(\k_{C,{\cal D}}l(w)\r)
\]
for each $v|\i.$ On the other hand it is clear that $\wtl{h}_{\a_i}(\s(\mf{q}_{1}^{-1})\c u_{i})_v\in \wh{K}_v$ for $v<\i.$ Therefore
\be\label{esth}
|\wtl{h}_{\a_i}(\s(\mf{q}_{1}^{-1})\c u_{i})|_{w^{-1}}\leq \exp\l(|S_\i|\k_{C,{\cal D}}l(w)\r),
\ee
where $S_\i$ is the set of infinite places of $F.$

Now consider
\[(a_g\s(\mf{q}_1))^{-1}\wtl{x}_\a(u_\a)(a_g\s(\mf{q}_1))=\wtl{x}_\a(a_g^\a \s(\mf{q}_1^{-1})\c u_\a).\]
Since $\Om$ is a compact subset of $T(\mb{A}),$ we may assume that
\[
\Om\ss \prod_{v\in S_\Om}T(F_v)\times \prod_{v\not\in S_\Om}T({\cal O}_v)
\]
where $S_\Om\supset S_\i$ is a finite set of places, and we can find $M_\Om>0$ such that $|a_v^\a|<M_\Om$ for any $a\in\Om,$ $v\in S_\Om$
and $\a\in\Phi.$ Applying Lemma \ref{est1} with $M'=M_\Om M_{\cal D}$ and $\tau'=C^{-1}$, we get a constant $\k_{C,{\cal D},\Om}:=2\k_{M',\tau'}$
such that
\[
|\wtl{x}_\a(a_g^\a \s(\mf{q}_1^{-1})\c u_\a)_v|_{w^{-1}}\leq \exp\l(\k_{C,{\cal D},\Om}l(w)\r)
\]
for $v\in S_\i.$ Applying Lemma \ref{estp} we get
\[
|\wtl{x}_\a(a_g^\a \s(\mf{q}_1^{-1})\c u_\a)_v|_{w^{-1}}\leq \exp\l(\k_{M_\Om}l(w)\r)
\]
for $v\in S_\Om\bs S_\i.$ Therefore
\be\label{estx}
|\wtl{x}_\a(a_g^\a \s(\mf{q}_1^{-1})\c u_\a)|_{w^{-1}}\leq \exp\l((|S_\i|\k_{C,{\cal D},\Om}+|S_\Om\bs S_\i|\k_{M_\Om})l(w)\r).
\ee
Combining (\ref{esth}) and (\ref{estx}),  there exists a constant $\wtl{\k}_{C, {\cal D}, \Om}$ such that
\[
|g_1^{-1}ug_1|_{w^{-1}}\leq \exp\l(\wtl{\k}_{C, {\cal D}, \Om}l(w)\r)
\]
for any $u\in\wh{U}_{\cal D},$ hence for any $u\in\wh{U}(F)\bs\wh{U}(\mb{A}).$ From (\ref{trick}) and (\ref{comp}), it follows that
\bq \label{comp2}
h_1(w\g g)\leq \exp\l(\wtl{\k}_{C, {\cal D}, \Om}l(w)\r)\s(\mf{q}_2^{-1})^{w^{-1}L}h_1(w\g u g_1).
\eq

Now we are ready to finish the proof of the theorem. We may assume that $s\in H_\v$ is a real number. Taking $s$th power of both sides of (\ref{comp2}) and integrating over $\wh{U}(F)\bs \wh{U}(\mb{A})$,
we obtain
\bq
E(s,h,g)&\leq & \sum_{w\in W(\D,\emptyset)}\exp\l(s\wtl{\k}_{C, {\cal D}, \Om}l(w)\r)\s(\mf{q}_2^{-1})^{sw^{-1}L}
\int_{\wh{U}(F)\bs \wh{U}(\mb{A})}\sum_{\g\in U_w(F)}h_s(w\g ug_1)du\n\\
&=& \sum_{w\in W(\D,\emptyset)}\exp\l(s\wtl{\k}_{C, {\cal D}, \Om}l(w)\r)\s(\mf{q}_2^{-1})^{sw^{-1}L}(a_g\s(\mf{q}_1))^{\wtl{\rho}+w^{-1}(sL-\wtl{\rho})}c_w(s)\n\\
&=& \sum_{w\in W(\D,\emptyset)}\exp\l(s\wtl{\k}_{C, {\cal D},
\Om}l(w)\r)\s(\mf{q}_2)^{\wtl{\rho}-w^{-1}\wtl{\rho}}(a_g\s(\mf{q}))^{\wtl{\rho}+w^{-1}(sL-\wtl{\rho})}c_w(s).\n
\eq Let us keep track of the proof of Lemma \ref{convh}. Let
\[
c_\Om=\max_{a\in\Om}c_a=\max_{a\in\Om, \a\in\Phi}|a^\a|\leq
M_\Om^{|S_\Om|}.
\]
Then (\ref{esta}) reads
\[
a_g^{\wtl{\rho}-w^{-1}\wtl{\rho}+w^{-1}sL}\leq
c_\Om^{l(w)+c_1s\|\la\|}.
\]
Plugging in (\ref{esti}), (\ref{estl}) and $c_w(s)\leq c_\v^{l(w)}$, we
see that \bq E(s,h,g) &\leq & |W|\exp(s\wtl{\k}_{C,{\cal
D},\Om}|\Phi_+|)(c_\v
c_\Om)^{|\Phi_+|}\n\\
&& \times\sum_{\la\in
Q^\vee}(c_\v^{c_2}c_\Om^{c_1s+c_2})^{\|\la\|}|\exp(s\wtl{\k}_{C,{\cal
D},\Om}c_2\|\la\|)|\mf{q}_2|^{c_1h\|\la\|+\f{h^\vee}{2}\|\la\|^2}|\mf{q}|^{c_1h\|\la\|-\f{s-h^\vee}{2}\|\la\|^2}.
\n\eq The last summation converges if and only if
\[
|\mf{q}_2|^{\f{h^\vee}{2}}|\mf{q}|^{-\f{s-h^\vee}{2}}<1,
\]
i.e. \be\label{conrg} s>c_\mf{q}:=h^\vee\l(1+\f{\log |\mf{q}_2|}{\log
|\mf{q}|}\r).\ee It is clear that convergence is uniform for all $s$
satisfying (\ref{conrg}).\hfill$\Box$

\

When $|\mf{q}|$ is large enough, we can replace $c_\mf{q}$ by a constant
which does not depend on $\mf{q}.$ We need the following lemma
\cite[p.143]{Lang}.

\begin{lemma}\label{lang}
There exists a constant $c_F>1$ depending on the number field $F$
such that, for all $\mf{q}\in\mb{I}$, $|\mf{q}|\geq c_F$ there exists $x\in
F^\times$ such that $1\leq |x\mf{q}_v|_v\leq |\mf{q}|$ for each place $v.$
\end{lemma}

Since the Eisenstein series are left $\s(F^\times)$-invariant, and
$\s(F^\times)$ normalizes $\wh{U}(\mb{A})T(\mb{A})$, we may replace
$\mf{q}$ in Theorem \ref{convae} by $x\mf{q}$  for any $x\in F^\times.$
Therefore we may assume that $\mf{q}$ satisfies the conclusion of Lemma
\ref{lang} whenever $|\mf{q}|\geq c_F.$ In this case $\mf{q}_2$ can be chosen
such that $|\mf{q}_{2v}|=1$ for $v<\i,$ and $|\mf{q}_{2v}|\leq C$ for
each $v|\i.$ Then we have \be\label{Ck} |\mf{q}_2|\leq C^{2|S_\i|},\q
c_{\mf{q}}\leq c'_F:=h^\vee\l(1+2|S_\i|\f{\log C}{\log c_F}\r). \ee In
summary we obtain the following result.

\begin{thm}
For any $\v>0$ and compact subset $\Om$ of $T(\mb{A})$, $E(s,f,g)$
and $E(s,h,g)$ converge absolutely and uniformly for
$s\in\{z\in\mb{C}|\mathrm{\mf{Re}}z\geq \max(h+h^\vee+\v, c'_F)\}$ and $g\in
\wh{U}(\mb{A})\Om\s_{c_F-1}\wh{K}$, where $c_F$ and $c'_F$ are given
by Lemma \ref{lang} and (\ref{Ck}) respectively.
\end{thm}

The constants in the theorem can be made explicit for the case $F=\mb{Q}.$ In fact
we get the same range of convergence as that of constant term of the Eisenstein series. Namely Theorem \ref{conv}
also holds for Eisenstein series itself.

\begin{cor}\label{corq}
Let $F=\mb{Q}.$ Then for any $\v, \eta>0$, $E(s,f,g)$
and $E(s,h,g)$ converge absolutely and uniformly for $s\in H_\v$ and
$g\in \wh{U}(\mb{A})\s_\eta\wh{K}.$
\end{cor}
\textit{Proof.} It is clear that $c_\mb{Q}$ in Lemma \ref{lang} can be chosen to be arbitrary constant greater than 1.
Fix $c_\mb{Q}=1+\eta$ with $\eta>0.$ We may choose
\[
{\cal D}=[-\f{1}{2},\f{1}{2})\times \prod_p\mb{Z}_p
\] and
$M_{\cal D}=\df{1}{2}.$ It follows from (\ref{cst}) and (\ref{Ck}) that we can choose $C$ to be close enough to $1$ such that
$c_\mb{Q}'\leq h^\vee+\v.$ \hfill$\Box$

\

As mentioned in the introduction, we conjecture that Corollary \ref{corq} holds for arbitrary number field, or in general for a global field $F$. For the
geometric case that $F$ is the function field of a smooth projective curve $X$ over a finite field $\mb{F}_q$, Lemma \ref{lang}
boils down to the Riemann-Roch theorem. Namely, let $D_\mf{q}$ be the divisor corresponding to $\mf{q}\in\mb{I}$, in order that
$H^0(-D_\mf{q})\neq 0$ it is sufficient that $-\deg(D_\mf{q})+1-g>0,$ i.e.
\[
|\mf{q}|=q^{-\deg D_\mf{q}}\geq q^g,
\]
where $g$ is the genus of $X.$ Hence we set $c_F=q^g+\eta$ and $c_F'=h^\vee,$ since $S_\i=\emptyset.$ In particular, the condition for $s$
reduces to $s\in H_\v$, and Corollary \ref{corq} is true for $X=\mb{P}^1_{\mb{F}_q}.$

\section*{Acknowledgement}

This paper is based on the author's Ph.D thesis. The author would like to thank his advisor, Prof. Yongchang Zhu, for many enlightening and helpful discussions during his Ph.D study at the Hong Kong University of Science and Technology.

\bibliographystyle{amsplain}

\end{document}